\documentclass[11pt,a4paper]{article}
\usepackage{amsfonts}
\usepackage{amssymb}
\usepackage{amsbsy}
\usepackage{amsthm}
\usepackage{graphicx}

\newtheorem{thm}{Theorem}[section]
\newtheorem{cor}[thm]{Corollary}
\newtheorem{lem}[thm]{Lemma}
\newtheorem{prop}[thm]{Proposition}
\theoremstyle{definition}
\newtheorem{defn}[thm]{Definition}
\theoremstyle{remark}
\newtheorem{rem}[thm]{\bf Remark}
\newtheorem{ex}[thm]{\bf Example}
\newcommand{\bt}{\begin{thm}}
\newcommand{\et}{\end{thm}}
\newcommand{\bc}{\begin{cor}}
\newcommand{\ec}{\end{cor}}
\newcommand{\bl}{\begin{lem}}
\newcommand{\el}{\end{lem}}
\newcommand{\bp}{\begin{prop}}
\newcommand{\ep}{\end{prop}}
\newcommand{\bd}{\begin{defn}}
\newcommand{\ed}{\end{defn}}
\newcommand{\br}{\begin{rem}}
\newcommand{\er}{\end{rem}}
\newcommand{\bpr}{\begin{proof}}
\newcommand{\epr}{\end{proof}}
\newcommand{\bex}{\begin{ex}}
\newcommand{\eex}{\end{ex}}

\newcommand{\bi}{\begin{itemize}}
\newcommand{\ei}{\end{itemize}}
\newcommand{\be}{\begin{enumerate}}
\newcommand{\ee}{\end{enumerate}}
\newcommand{\ds}{\displaystyle}
\newcommand{\ba}{\begin{array}}
\newcommand{\ea}{\end{array}}
\newcommand{\beq}{\begin{equation}}
\newcommand{\eeq}{\end{equation}}
\newcommand{\beqa}{\begin{eqnarray}}
\newcommand{\eeqa}{\end{eqnarray}}

\newcommand{\N}{{\mathbb N}}
\newcommand{\Z}{{\mathbb Z}}
\newcommand{\R}{{\mathbb R}}
\newcommand{\C}{{\mathbb C}}
\newcommand{\T}{{\mathbb T}}
\newcommand{\D}{{\mathbb D}}

\newcommand{\PP}{{\mathbb P}}

\newcommand{\cB}{{\mathcal  B}}

\newcommand{\cD}{{\mathcal  D}}

\newcommand{\cI}{{\mathcal  I}}
\newcommand{\cJ}{{\mathcal  J}}
\newcommand{\cK}{{\mathcal  K}}

\newcommand{\frB}{{\mathfrak B}}

\newcommand{\bsa}{{\boldsymbol a}}
\newcommand{\bsb}{{\boldsymbol b}}
\newcommand{\bsu}{{\boldsymbol u}}

\newcommand{\bsc}{{\boldsymbol \chi}}

\newcommand{\bsE}{{\boldsymbol E}}

\newcommand{\spn}{\mathrm{span}}
\newcommand{\supp}{\mathrm{supp}}
\newcommand{\spec}{\sigma}

\newcommand{\re}{\mathrm{Re}}

\newcommand{\lims}{\mathop{\overline{\lim}}}
\newcommand{\limi}{\mathop{\underline{\lim}}}
\newcommand{\limss}{\mathop{\overline{\overline{\lim}}}}
\newcommand{\rank}{\mathrm{rank}}

\newcommand{\even}{\mathrm{even}}
\newcommand{\odd}{\mathrm{odd}}
\newcommand{\sg}{\mathrm{sign}}
\newcommand{\co}{\mathrm{Co}}

\begin{document}

\title{\bf Measures on the unit circle and unitary truncations of unitary operators}
\author{M.J. Cantero $^a$, L. Moral $^b$, L. Vel\'azquez $^c$
\thanks{The work of the authors was supported by Project E-12/25 of Diputaci\'on
General de Arag\'on (Spain) and by Ibercaja under grant IBE2002-CIEN-07.} \\
\small{Departamento de Matem\'atica Aplicada} \\
\small{Universidad de Zaragoza} \\
\small{50009 Zaragoza, Spain} \\
\small{(a) \texttt{mjcante@unizar.es}} \\
\kern-5pt\small{(b) \texttt{lmoral@unizar.es}} \\
\kern5pt\small{(c) \texttt{velazque@unizar.es}}
}
\date{March 12, 2005}
\maketitle

\begin{abstract}

In this paper we obtain new results about the orthogonality measure of orthogonal
polynomials on the unit circle, through the study of unitary truncations of the
corresponding unitary multiplication operator, and the use of the five-diagonal
representation of this operator.

Unitary truncations on subspaces with finite co-dimension give information about the
derived set of the support of the measure under very general assumptions for the related
Schur parameters $(a_n)$. Among other cases, we study the derived set of the support of
the measure when $\lim_n|a_{n+1}/a_n|=1$, obtaining a natural generalization of the known
result for the L\'opez class $\lim_na_{n+1}/a_n\in\T$, $\lim_n|a_n|\in(0,1)$.

On the other hand, unitary truncations on subspaces with finite dimension provide
sequences of unitary five-diagonal matrices whose spectra asymptotically approach the
support of the measure. This answers a conjecture of L. Golinskii concerning the relation
between the support of the measure and the strong limit points of the zeros of the
para-orthogonal polynomials.

Finally, we use the previous results to discuss the domain of convergence of rational
approximants of Carath\'eodory functions, including the convergence on the unit circle.

\end{abstract}


\medskip

\noindent{\it (2000) AMS Mathematics Subject Classification}: 42C05, 47B36.

%%%%%%%%%%%%%%%%%%%%%%%%%%%%%%%%%%%%%%%%%%%%%%%%%%%%%%%%%%%%%%%%%%%%%%%%%%%%%%%%%%%%%
\section{Introduction}
%%%%%%%%%%%%%%%%%%%%%%%%%%%%%%%%%%%%%%%%%%%%%%%%%%%%%%%%%%%%%%%%%%%%%%%%%%%%%%%%%%%%%

In \cite{FIVE} and \cite{MIN}, a new operator theoretic approach for the orthogonal
polynomials with respect to a measure on the unit circle $\T:=\{z\in\C:|z|=1\}$ was
established. The five-diagonal representation of unitary operators introduced there has
proved to be a powerful tool for the study of such orthogonal polynomials, as it has been
shown in \cite{Si04-1,Si04-2}, where numerous new results have been obtained (for a summary
of some of the main new results in \cite{Si04-1,Si04-2}, see \cite{SiNEW}). Let us summarize
the main facts concerning this five-diagonal representation, since it is the starting point
of this paper.

In what follows $\mu$ denotes a probability measure on $\T$ with an infinite support
$\supp\,\mu$. Then,
$$
U^\mu \colon \mathop{L^2_\mu \longrightarrow L^2_\mu} \limits_{f(z) \; \to \; zf(z)}
$$
is a unitary operator on the Hilbert space $L^2_\mu$ of $\mu$-square-integrable functions
with the inner product
$$
(f,g) := \int f(z)\overline{g(z)}\,d\mu(z), \qquad \forall f,g \in L^2_\mu.
$$
The associated five-diagonal representation is just the matrix representation of $U^\mu$
with respect to an orthonormal Laurent polynomial basis $(\chi_n)_{n\geq0}$ of $L^2_\mu$.
This basis is related to the usual orthonormal polynomials
$(\varphi_n)_{n\geq0}$ in $L^2_\mu$, defined by
$$
\varphi_n(z)=\kappa_n(z^n+\cdots+a_n), \quad \kappa_n>0, \quad (\varphi_n,\varphi_m)=\delta_{n,m},
\quad n,m\geq0,
$$
through the relations \cite{Th88,FIVE}
\beq \label{OP-OLP}
\chi_{2j}(z) = z^{-j}\varphi_{2j}^*(z), \quad
\chi_{2j+1} = z^{-j} \varphi_{2j+1}(z), \quad j\geq0,
\eeq
where, for every polynomial $p$ of degree $n$, $p^*(z) := z^n \overline p(z^{-1})$ is called
the reversed polynomial of $p$. The five-diagonal representation has the form \cite{FIVE}
\beq \label{C}
C(\bsa):=\pmatrix{
-a_1 & \kern-7pt -\rho_1 a_2 & \rho_1 \rho_2
\cr
\kern7pt \rho_1 & \kern-7pt -\overline a_1 a_2 & \overline a_1 \rho_2 & 0
\cr
\kern5pt 0 & \kern-7pt -\rho_2 a_3 & \kern-7pt -\overline a_2 a_3 &
\kern-7pt -\rho_3 a_4 & \rho_3 \rho_4
\cr
& \rho_2 \rho_3 & \overline a_2 \rho_3 & \kern-7pt -\overline a_3 a_4  &
\overline a_3 \rho_4 & 0
\cr
&& 0 & \kern-7pt -\rho_4 a_5 & \kern-7pt -\overline a_4 a_5 &
\kern-7pt -\rho_5 a_6 & \rho_5 \rho_6
\cr
&&& \rho_4 \rho_5 & \overline a_4 \rho_5 & \kern-7pt -\overline a_5 a_6 &
\overline a_5 \rho_6 & 0
\cr
&&&& \hskip-35pt\ddots & \hskip-35pt\ddots &
\hskip-35pt\ddots & \hskip-20pt\ddots & \ddots},
\eeq
where $\bsa:=(a_n)_{n\geq1}$ satisfies $|a_n|<1$ and
$\rho_n:=\kappa_{n-1}/\kappa_n=\sqrt{1-|a_n|^2}$.
The transposed matrix $C(\bsa)^t$ of $C(\bsa)$ is also a representation of $U^\mu$, but
with respect to the orthonormal Laurent polynomial basis $(\chi_{n*})_{n\geq0}$,
where $f_*(z):=\overline f(z^{-1})$ for any Laurent polynomial $f$.

$C(\bsa)$ and $C(\bsa)^t$ can be identified with unitary operators on the Hilbert space
$\ell^2$ of square-sumable sequences in $\C$, these operators being unitarily equivalent
to $U^\mu$. Due to the properties of the multiplication operator, the spectrum of
$C(\bsa)$ and $C(\bsa)^t$ coincides with $\supp\,\mu$, the mass points being the
corresponding eigenvalues. Since the eigenvalues are simple, the essential spectrum of
$C(\bsa)$ and $C(\bsa)^t$ (that is, the spectrum except the isolated eigenvalues with
finite multiplicity) is the derived set $\{\supp\,\mu\}'$ of $\supp\,\mu$.

$\bsa$ is called the sequence of Schur parameters of $\mu$. The Schur parameters
establish a one to one correspondence between sequences in the open unit disk
$\D:=\{z\in\C : |z|<1\}$ and probability measures on $\T$ with infinite support. The
Schur parameters also appear in the forward recurrence relation
\beq \label{RR-OP}
\rho_n\varphi_n(z) = z\varphi_{n-1}(z) + a_n \varphi_{n-1}^*(z), \quad n\geq1,
\eeq
that generates the orthonormal polynomials, which is also equivalent to the backward
recurrence relation
\beq \label{RR-OP2}
\rho_nz\varphi_{n-1}(z) = \varphi_n(z) - a_n \varphi_n^*(z), \quad n\geq1.
\eeq

Therefore, the matrix $C(\bsa)$ provides a connection between this practical way of
constructing sequences of orthonormal polynomials on $\T$, and the properties of the
related orthogonality measure $\mu$, lost in such construction. In particular, the
spectral analysis of $C(\bsa)$ permits us to recover features of $\supp\,\mu$ from
properties of the Schur parameters $\bsa$. Only in the case $|a_n|<1$ the matrices
$C(\bsa)$ are related to a measure on $\T$ with infinite support. However, they are
well defined unitary matrices even if $|a_n|\leq1$. As we will see, this is important
when using perturbative arguments for the analysis of the measure.

In this paper we use unitary truncations of $C(\bsa)$ as a source of information for the
spectrum and the essential spectrum of $C(\bsa)$, that is, for $\supp\,\mu$ and
$\{\supp\,\mu\}'$. Let us denote by $\{e_n\}_{n\geq1}$ the canonical basis of $\ell^2$,
and let $\ell^2_n:=\spn\{e_1,e_2,\dots,e_n\}$. As we will see in the following section,
for any infinite bounded normal band matrix, the normal truncations on $\ell^2_n$ or
$\ell^{2\bot}_n$ have a spectrum closely related to the spectrum of the full matrix. This
justifies the study of this kind of normal truncations for $C(\bsa)$, which is the aim of
Section 2. We find that all these truncations are indeed unitary and can be parameterized
by the points in $\T$, leading to the para-orthogonal polynomials \cite{JoNjTh89} and to
the family of Aleksandrov measures \cite{GoNe01} related to the associated polynomials
\cite{Pe96}.

In Section 3 we use the unitary truncations of $C(\bsa)$ on $\ell^{2\bot}_n$ to obtain
general relations between the asymptotic behaviour of the Schur parameters and the
location of $\{\supp\,\mu\}'$. A well known result is that, under the conditions
$\lim_na_{n+1}/a_n\in\T$, $\lim_n|a_n|\in(0,1)$, which define the so-called L\'opez
class, $\{\supp\,\mu\}'$ is a closed arc centred at $-\lim_na_{n+1}/a_n$ with angular
radius $2\arccos(\lim_n|a_n|)$ \cite{BaLo99} (see also \cite[Chapter 4]{Si04-1} for an
approach using the five-diagonal representation $C(\bsa)$). It is of interest to extend
these results to a bigger class than the L\'opez one. With our techniques we can get
information about $\{\supp\,\mu\}'$ when only one of the two L\'opez conditions is
satisfied, or, even, under the more general condition $\lim_n|a_{n+1}/a_n|=1$. Among other
results, we prove that, if $\lim_n|a_{n+1}/a_n|=1$, $\{\supp\,\mu\}'$ lies inside the
union of closed arcs with centre at the limit points of the sequence
$(-a_{n+1}/a_n)_{n\geq1}$ and angular radius $2\arccos(\limi_n|a_n|)$, which is a natural
generalization of the result for the L\'opez class.

Section 4 is devoted to the study of the approximation of $\supp\,\mu$ by means of the
spectra of the unitary truncations of $C(\bsa)$ on $\ell^2_n$, which means the approximation
of the spectrum of an infinite unitary matrix by the spectra of finite unitary matrices. The
proofs now include methods, not only from operator theory, but also from the theory of
analytical functions. We prove that, for any measure $\mu$ on $\T$, there exist infinitely
many sequences of unitary truncations whose spectra exactly converge to $\supp\,\mu$, in a
strong sense that we will specify later on. This result proves a conjecture formulated by L.
Golinskii in \cite{Go02}, concerning the coincidence of the support of a measure on $\T$ and
the strong limit points of the zeros of the related para-orthogonal polynomials. We also
present some other results that deal with weaker notions of convergence of the finite
spectra, which are of interest in the following section.

Finally, in Section 5, we consider an application of the previous results, that is, the
study of the convergence of rational approximants of the Carath\'eodory function of a
measure on $\T$. It is known that the standard rational approximants constructed with the
related orthogonal polynomials, or their reversed ones, converge on
$\C\backslash\overline\D$ and $\D$ respectively. We focus our analysis on the study of
the rational approximants related to the unitary truncations of $C(\bsa)$ on $\ell^2_n$,
that always converge on $\C\backslash\T$, and are just the rational approximants
constructed with the para-orthogonal polynomials \cite{JoNjTh89}. The domain of
convergence of these approximants is closely related to the asymptotic behaviour of the
finite spectra of the above unitary truncations. Therefore, the results of the previous
sections give information about the convergence of these approximants on $\T$, where the
situation is more delicate. Some results in this direction for the standard rational
approximants can be found in \cite{Kh02}.

%%%%%%%%%%%%%%%%%%%%%%%%%%%%%%%%%%%%%%%%%%%%%%%%%%%%%%%%%%%%%%%%%%%%%%%%%%%%%%%%%%%%%
\section{Normal truncations of $C(\bsa)$}
%%%%%%%%%%%%%%%%%%%%%%%%%%%%%%%%%%%%%%%%%%%%%%%%%%%%%%%%%%%%%%%%%%%%%%%%%%%%%%%%%%%%%

In what follows, given a Hilbert space $H$, $(\cdot,\cdot)$ is the corresponding inner
product and $\|\cdot\|$ the related norm. We will deal with the set $\frB(H)$ of bounded
linear operators on $H$. $\|\cdot\|$ also denotes the standard operator norm in $\frB(H)$
while, for any operator $T$ on $H$,
$$
\|T\|_S:=\sup_{x\in S\backslash\{0\}} \frac{\|Tx\|}{\|x\|}, \qquad
\gamma(T;S) := \inf_{x\in S\backslash\{0\}} \frac{\|Tx\|}{\|x\|},
\qquad \forall S \subset H,
$$
and $\gamma(T):=\gamma (T;H)$. Given a sequence $(T_n)_{n\geq1}$ of operators on $H$,
$T_n \to T$ means that $\lim_n\|T_nx-Tx\|=0$, $\forall x \in H$ ($T$ is the strong limit
of $(T_n)_{n\geq1}$).

Let $S$ be a subspace of $H$. If an operator $T$ leaves $S$ invariant, the operator on
$S$ defined by the restriction of $T$ to $S$ is denoted by $T_{\upharpoonright S}$. In
particular, $0_{\upharpoonright S}$ and $1_{\upharpoonright S}$ are the null and
identity operators on $S$ respectively (the identity operator will be omitted when it
is clear from the context). Also, if $T$ is an operator on $S$, we define an operator
on $H$ by $\hat T := T \oplus 0_{\upharpoonright S^\bot}$.

If $T\in\frB(H)$, $\spec(T)$ is its spectrum and $\spec_e(T)$ its essential spectrum.
When $T$ is normal it is known that $\spec(T)=\{z\in\C : \gamma(z-T)=0\}$. In fact, in
this case, denoting by $d(\cdot,\cdot)$ the distance between points and sets in $\C$,
we have $\gamma(z-T)=d(z,\spec(T))$ for any $z\in\C$.

Let $T\in\frB(H)$ and $Q$ be a projection on $S \subset H$ along $S'\subset H$. The
operator $T[Q]:=QT_{\upharpoonright S}$ is called the truncation of $T$ associated with
$Q$, or the truncation of $T$ on $S$ along $S'$. $T[Q]$ is finite (co-finite) when $S$
has finite dimension (co-dimension). If $Q$ is an orthogonal projection we say that $T[Q]$
is an orthogonal truncation. To compare the operator with its truncation, it is
convenient to consider $\hat T[Q] = T[Q] \oplus 0_{\upharpoonright S^\bot} = QTP$, where
$P$ is the orthogonal projection on $S$. Notice that, if $T[Q]$ is bounded (for example,
this is the case of a finite truncation), $\|\hat T[Q]\|=\|T[Q]\|$. Also, $\hat T[Q]^* =
T[Q]^* \oplus 0_{\upharpoonright S^\bot}$, so, any orthogonal truncation of a
self-adjoint operator is self-adjoint too. However, in general, to get normal truncations
of a normal operator can require non-orthogonal truncations.

In what follows, any infinite bounded matrix $M$ is identified with the operator
$T\in\frB(\ell^2)$ defined by $Tx=Mx$, $\forall x\in\ell^2$. $T$ is called a band
operator if $M$ is a band matrix. The following result shows the interest in finding
normal finite and co-finite truncations of a normal band operator.

\bp \label{GENTRUNC}

Let $T\in\frB(\ell^2)$ be a normal band operator.
\be
\item If $T_n$ is a normal truncation of $T$ on $\ell^2_n$ for $n\geq1$ and
$(\|T_n\|)_{n\geq1}$ is bounded, then $\hat T_n \to T$ and
$$
\spec(T) \subset \{z\in\C : \lim_nd(z,\spec(T_n))=0\}.
$$
\item For any bounded normal truncation $T_n$ of $T$ on $\ell^{2\bot}_n$,
$$
\spec_e(T_n) = \spec_e(T).
$$
\ee

\ep

\bpr

Let $T_n$ be a normal truncation of $T$ on $\ell^2_n$. If $T$ is a $2N+1$-band operator,
$T\ell^2_n\subset\ell^2_{n+N}$ for $n\geq1$. Therefore, $\hat T_nx=Tx$ if
$x\in\ell^2_{n-N}$, $n>N$, and we get
$$
\|\hat T_nx-Tx\| \leq (\|T_n\|+\|T\|)\|x-P_{n-N}x\|,
\quad \forall n>N, \quad \forall x\in\ell^2.
$$
Since $P_n\to1$ and $(\|T_n\|)_{n\geq1}$ is bounded we find that $\hat T_n \to T$.

If $\lims_nd(z,\spec(T_n))>0$, there exist $\delta>0$
and an infinite set $\cI\subset\N$ such that $d(z,\spec(T_n))\geq\delta$, $\forall
n\in\cI$. Since $T_n$ is normal, $\gamma(z-T_n) \geq \delta$, $\forall n\in\cI$. Hence,
if $P_n$ is the orthogonal projection on $\ell^2_n$,
\beq \label{FINITETRUNC}
\|(zP_n-\hat T_n)x\| = \|(z-T_n)P_nx\| \geq \delta \|P_nx\|,
\quad \forall x\in\ell^2, \quad \forall n\in\cI.
\eeq
Taking limits in (\ref{FINITETRUNC}) we obtain
$\|(z-T)x\| \geq \delta \|x\|$, $\forall x\in\ell^2$,
which, taking into account that $T$ is normal, implies that $z \notin \spec(T)$.
This proves 1.

As for the second statement, let $T_n$ be a normal truncation of $T$ on $\ell^{2\bot}_n$.
If $T$ is $2N+1$-band, $T\ell^{2\bot}_n\subset\ell^{2\bot}_{n-N}$ for $n>N$. Hence, if
$n\geq1$, $\hat T_nx=Tx$ for $x\in\ell^{2\bot}_{n+N}$, and, thus, $\rank\,(T-\hat T_n)
\leq n+N$. Since $T-\hat T_n$ has finite rank, Weyl's theorem implies that
$\spec_e(T)=\spec_e(\hat T_n)=\spec_e(T_n)$.

\epr

The first statement of the above proposition is the first step in establishing a numerical
method for the approximation of the spectrum of an infinite bounded normal band operator.
This statement will be improved in the case of finite normal truncations of $C(\bsa)$,
getting an equality instead of an inclusion (see Section 4), which can be used for the
numerical approximation of the support of the related measure on $\T$.

The importance of the second assertion of Proposition \ref{GENTRUNC} is that it can be
used to extract properties of the essential spectrum of an infinite bounded normal band
operator from the asymptotic behaviour of the coefficients of its diagonals.
When applying Proposition \ref{GENTRUNC} to co-finite truncations of $C(\bsa)$, we can obtain
properties of the derived set of the support of a measure on $\T$ from the asymptotic
behaviour of the related Schur parameters (see Section 3).

Our next step is to study the normal truncations of $C(\bsa)$ for an arbitrary sequence
$\bsa$ in $\D$. This is equivalent to studying the normal truncations of $U^\mu$, where $\mu$ is
the measure on $\T$ related to $\bsa$. Concerning this problem we have the following result.

\bt \label{NORMALTRUNC}

Let $\mu$ be a measure on $\T$ with infinite support and
$\PP_{m,n}:=\spn\{z^m,z^{m+1},\dots,z^{m+n-1}\}$, $m\in\Z$, $n\in\N$. The normal
truncations of $U^\mu$ on $\PP_{m,n}$ $(\PP_{m,n}^\bot)$ are unitary, and they are
parameterized by the points in $\T$. The normal truncation on $\PP_{m,n}$
$(\PP_{m,n}^\bot)$ corresponding to a parameter $u\in\T$ is
$U^\mu_{m,n}(u):=U^\mu[Q^\mu_{m,n}(u)]$, where $Q^\mu_{m,n}(u)$ is the projection on
$\PP_{m,n}$ $(\PP_{m,n}^\bot)$ along $\spn\{z^mp_n^u\}\oplus\PP_{m,n+1}^\bot$
$(\spn\{z^mq_n^u\} \oplus \PP_{m+1,n-1})$, and
$$
p_n^u(z) := z\varphi_{n-1}(z)+u\varphi_{n-1}^*(z), \qquad
q_n^u(z) := \varphi_n^*(z)-\overline u\varphi_n(z),
$$
$(\varphi_n)_{n\geq0}$ being the orthonormal polynomials in $L^2_\mu$. If $(a_n)_{n\geq1}$
are the Schur parameters of $\mu$,
$$
\|Q^\mu_{m,n}(u)\|=\sqrt{1+|u-a_n|^2/\rho_n^2}.
$$
The spectrum of the truncation $U^\mu_{m,n}(u)$ on $\PP_{m,n}$ is simple and coincides with
the zeros of $p_n^u$.

\et

\bpr

The problem can be reduced to the study of the normal truncations
$U^\mu_n=U^\mu[Q^\mu_n]$ of $U^\mu$ on $\PP_n:=\PP_{0,n}$ $(\PP_n^\bot)$, since the normal
truncations $U^\mu_{m,n}=U^\mu[Q^\mu_{m,n}]$ on $\PP_{m,n}$ $(\PP_{m,n}^\bot)$ are related
to the previous ones by $Q^\mu_{m,n}=z^mQ^\mu_nz^{-m}$ and $U^\mu_{m,n}=z^mU^\mu_nz^{-m}$.

Let $U^\mu_n$ be a truncation of $U^\mu$ on $\PP_n$. For any $f\in\PP_n$, the decomposition
$U^\mu_nf = z(f-(f,\varphi_{n-1})\varphi_{n-1})+(f,\varphi_{n-1})U^\mu_n\varphi_{n-1}$ gives
\beq \label{FINTRUNC}
U^\mu_nf = zf-(f,\varphi_{n-1})p_n, \qquad
p_n=z\varphi_{n-1}-f_n, \qquad f_n=U^\mu_n\varphi_{n-1}.
\eeq
Thus, $U^\mu_n=U^\mu[Q^\mu_n]$, $Q^\mu_n$ being the projection on $\PP_n$ along
$\spn\{p_n\}\oplus\PP_{n+1}^\bot$.

From (\ref{FINTRUNC}), for an arbitrary $f\in\PP_n$, we get
$$
U^{\mu*}_nf = z^{-1}(f-(f,\varphi_n^*)\varphi_n^*)-(f,p_n)\varphi_{n-1},
$$
and, therefore,
\beqa \label{UU*}
& \kern-58pt U^\mu_n U^{\mu*}_nf = f-(f,\varphi_n^*)\varphi_n^*-(f,z\varphi_{n-1})z\varphi_{n-1}
+(f,f_n)f_n,
\smallskip \\ \label{U*U}
& \ba{l}
U^{\mu*}_nU^\mu_n f = f+(zf,f_n)\varphi_{n-1} \, +
\smallskip \\ \kern63pt
+ \, (f,\varphi_{n-1})
\left\{z^{-1}(f_n-(f_n,\varphi_n^*)\varphi_n^*)+(\|p_n\|^2-2)\varphi_{n-1}\right\}.
\ea
\eeqa

Let us suppose that $U^\mu_n$ is normal, that is
$(U^\mu_nU^{\mu*}_n-U^{\mu*}_nU^\mu_n)f=0$ for any $f\in\PP_n$. Using (\ref{UU*}) and
(\ref{U*U}) we find that
$$
(f,f_n)f_n = (zf,f_n)\varphi_{n-1}, \quad \forall f \in z\PP_{n-2}.
$$
If $f \in z\PP_{n-2}$ is such that $(f,f_n)\neq0$, the above equality implies that
$f_n$ is proportional to $\varphi_{n-1}$, which gives a contradiction since
$\varphi_{n-1}\bot z\PP_{n-2}$.
So, $(f,f_n)=(zf,f_n)=0$ for any $f \in z\PP_{n-2}$, that is,
$$
f_n\in(z\PP_{n-2}+z^2\PP_{n-2})^{\bot\PP_n}=z\PP_{n-1}^{\bot\PP_n}
=\spn\{\varphi_{n-1}^*\}.
$$
Therefore, $f_n=-u\varphi_{n-1}^*$, $u\in\C$. Then, if we take $f=1$ in (\ref{UU*}) and
(\ref{U*U}), the condition $(U^\mu_nU^{\mu*}_n-U^{\mu*}_nU^\mu_n)1=0$ gives
$|u|^2\varphi_{n-1}^* = \rho_n\varphi_n^* - \overline a_n z\varphi_{n-1}$, and the
reversed form of (\ref{RR-OP}) shows that $u\in\T$.

Moreover, if $f_n=-u\varphi_{n-1}^*$, $u\in\T$, we get from (\ref{UU*})
$$
U^\mu_nU^{\mu*}_nf-f =
(f,\varphi_{n-1}^*)(\varphi_{n-1}^*-\rho_n\varphi_n^*+\overline a_nz\varphi_{n-1}) = 0,
$$
and, hence, the finite truncation $U^\mu_n$ is unitary.

\medskip

Let us consider now a truncation $U^\mu_n$ of $U^\mu$ on $\PP_n^\bot$. The orthogonal
decomposition
$\PP_n^\bot=z^{-1}\PP_{n+1}^\bot\oplus\spn\{z^{-1}\varphi_n^*\}$
gives, for any $f\in\PP_n^\bot$, the equality
$U^\mu_nf = z(f-(f,z^{-1}\varphi_n^*)z^{-1}\varphi_n^*)
+(f,z^{-1}\varphi_n^*)U^\mu_nz^{-1}\varphi_n^*$,
and, thus,
\beq \label{COFINTRUNC}
U^\mu_nf = zf-(zf,\varphi_n^*)q_n, \qquad
q_n=\varphi_n^*-g_n, \qquad g_n=U^\mu_nz^{-1}\varphi_n^*.
\eeq
So, $U^\mu_n=U^\mu[Q^\mu_n]$, $Q^\mu_n$ being the projection on $\PP_n^\bot$ along
$\spn\{q_n\}\oplus z\PP_{n-1}$.

From (\ref{COFINTRUNC}) we find that, for an arbitrary $f\in\PP_n^\bot$,
$$
U^{\mu*}_nf = z^{-1}(f-(f,z\varphi_{n-1})z\varphi_{n-1})-(f,q_n)z^{-1}\varphi_n^*,
$$
and
$$
\ba{l}
U^\mu_n U^{\mu*}_nf = f-(f,\varphi_n^*)\varphi_n^*-(f,z\varphi_{n-1})z\varphi_{n-1}
+(f,g_n)g_n,
\smallskip \\
U^{\mu*}_nU^\mu_n f =f+(zf,g_n)z^{-1}\varphi_n^* \, +
\smallskip \\ \kern55pt
+ \, (zf,\varphi_n^*)
\left\{z^{-1}(g_n-(g_n,z\varphi_{n-1})z\varphi_{n-1})
+(\|q_n\|^2-2)z^{-1}\varphi_n^*\right\}.
\ea
$$

When $(U^\mu_nU^{\mu*}_n-U^{\mu*}_nU^\mu_n)f=0$ for any $f\in\PP_n^\bot$, we find that
$$
(f,g_n)g_n = (zf,g_n)z^{-1}\varphi_n^*, \quad \forall f \in z^{-1}\PP_{n+2}^\bot,
$$
and, in a similar way to the previous case, we obtain that
$$
g_n\in\PP_n^{\bot(\PP_{n+2} \cap z^{-1}\PP_{n+2})}=\PP_n^{\bot\PP_{n+1}}
=\spn\{\varphi_n\}.
$$
Therefore, $g_n=\overline u\varphi_n$, $u\in\C$, and
$(U^\mu_nU^{\mu*}_n-U^{\mu*}_nU^\mu_n)z^{-1}\varphi_{n+1}=0$ gives
$|u|^2\varphi_n = \rho_nz\varphi_{n-1} + a_n \varphi_n^*$. From (\ref{RR-OP2}) we conclude
that $u\in\T$.

Moreover, if $g_n=\overline u\varphi_n$, $u\in\T$,
$$
\ba{l}
U^\mu_nU^{\mu*}_nf-f =
(f,\varphi_n)(\varphi_n-\rho_nz\varphi_{n-1}-a_n\varphi_n^*) = 0,
\smallskip \\
U^{\mu*}_nU^\mu_nf-f =
u\rho_n(f,\varphi_{n-1})z^{-1}\varphi_n^* = 0,
\ea
$$
which shows that $U^\mu_n$ is unitary.

\medskip

Notice that $p_n=\rho_n\varphi_n+(u-a_n)\varphi_{n-1}^*$ and
$q_n=\rho_n\varphi_{n-1}^*+(\overline a_n-\overline u)\varphi_n$. Therefore,
$$
Q^\mu_n \varphi_{n-1}^* = \alpha_n \varphi_n,
\qquad Q^\mu_n \varphi_n = \beta_n \varphi_{n-1}^*,
$$
where $\alpha_n=1$, $\beta_n=(a_n-u)/\rho_n$ in the finite case and $\alpha_n=(\overline
u-\overline a_n)/\rho_n$, $\beta_n=1$ in the co-finite case. In any of the two cases,
$Q^\mu_{n \upharpoonright z\PP_{n-1}}$ and $Q^\mu_{n \upharpoonright \PP_{n+1}^\bot}$ are
the unit or null operators. Since
$L^2_\mu=z\PP_{n-1}\oplus\spn\{\varphi_{n-1}^*,\varphi_n\}\oplus\PP_{n+1}^\bot$ we find
that $\|Q^\mu_n\| = \|Q^\mu_{n \upharpoonright \spn\{\varphi_{n-1}^*,\varphi_n\}}\| =
\sqrt{1+|u-a_n|^2/\rho_n^2}$.

\medskip

Finally, $1$ is a cyclic vector for any finite normal truncation $U^\mu_n$ since
$\spn\{U^{\mu k}_n1\}_{k=0}^{n-1}=\PP_n$, so, the spectrum of $U^\mu_n$ is simple.
The identity (\ref{FINTRUNC}) implies that any eigenvalue of $U^\mu_n$ is a zero
of $p_n$, hence, the $n$ different eigenvalues of $U^\mu_n$ must fulfill the $n$
(simple) zeros of $p_n$.

\epr

\br \label{POP}

The polynomials $p_n^u$ and $q_n^u$ can be understood as the substitutes for
$\rho_n\varphi_n$ and $\rho_n\varphi_{n-1}^*$, when changing $a_n\in\D$ by $u\in\T$
in the $n$-th step of (\ref{RR-OP}) and the reversed version of (\ref{RR-OP2})
respectively. In fact, using (\ref{RR-OP}) and (\ref{RR-OP2}) we see that these
polynomials are related by
\beq \label{u-v}
\kern-5pt
p_n^u=\frac{u-a_n}{\rho_n}q_n^v, \kern9pt
q_n^v=\frac{\overline a_n - \overline v}{\rho_n}p_n^u, \kern9pt
u=-v\frac{1-a_n\overline v}{1-\overline a_nv}, \kern9pt
v=-u\frac{1-a_n\overline u}{1-\overline a_nu}.
\eeq
They are called para-orthogonal polynomials of order $n$ associated with the measure
$\mu$ \cite{JoNjTh89}. It was known that they have simple zeros lying on $\T$,
which is in agreement with Theorem \ref{NORMALTRUNC}. Notice that the freedom in the
parameter $u\in\T$ means that we can arbitrarily fix in $\T$ one of the zeros of a
para-orthogonal polynomial with given order, the rest of the zeros being determined by
this choice.

\er

\medskip

The importance of Theorem \ref{NORMALTRUNC} is that it provides and, at the same time,
closes the possible ways of applying Proposition \ref{GENTRUNC} to the unitary matrix
$C(\bsa)$. We will identify any truncation of $C(\bsa)$ on $\ell^2_n$ ($\ell^{2\bot}_n$)
with its matrix representation with respect to $\{e_k\}_{k \leq n}$
($\{e_k\}_{k > n}$). The matrix form of the normal truncations of $C(\bsa)$ can be
obtained using the decomposition
\beq \label{DEC}
C(\bsa) = C(a_1,\dots,a_n) \oplus C(\bsa^{(n)};a_n),
\quad |a_n|=1, \quad \bsa^{(n)}:=(a_{n+j})_{j\geq1},
\eeq
where $C(a_1,\dots,a_n)$ is the principal matrix of $C(\bsa)$ of order $n$, and
\beq \label{DEC2}
\kern-5pt C(\bsa;u) :=
\cases{W C(\overline u\bsa) W^* & even $k$,
\smallskip \cr
(W C(\overline u\bsa) W^*)^t & odd $k$,}
\quad
W:=\pmatrix{ w &&&& \cr & \kern-3pt \overline w &&& \cr && \kern-3pt w && \cr
&&& \kern-3pt \overline w & \cr &&&& \kern-3pt \ddots },
\eeq
with $u=w^2$.
The factorization $C(\bsa)=C_o(\bsa)C_e(\bsa)$ is useful, where
\beq \label{FAC}
C_o(\bsa) := \pmatrix{
\Theta(a_1) & & & \cr
 & \kern-12pt\Theta(a_3) & & \cr
 & & \kern-12pt\Theta(a_5) & \cr
 & & & \kern-12pt\ddots
},
\quad
C_e(\bsa) := \pmatrix{
1 & & & \cr
 & \kern-5pt\Theta(a_2) & & \cr
 & & \kern-12pt\Theta(a_4) & \cr
 & & & \kern-12pt\ddots
},
\eeq
\vskip-5pt
$$
\Theta(a) := \pmatrix{
-a & \rho \cr
\rho & \overline a
},
\quad \rho:=\sqrt{1-|a|^2},
\quad |a|\leq1.
$$
Also, $C(a_1,\dots,a_n)=C_o(a_1,\dots,a_n)C_e(a_1,\dots,a_n)$, where
$C_o(a_1,\dots,a_n)$ and $C_e(a_1,\dots,a_n)$ are the principal submatrices of order $n$
of $C_o(\bsa)$ and $C_e(\bsa)$ respectively. All these properties hold for
$|a_n|\leq1$ \cite{MIN}.

If $(\chi_n)_{n\geq0}$ are the Laurent orthonormal polynomials related to a sequence
$\bsa$ of Schur parameters, we also have (see \cite{FIVE})
\beq \label{THETA}
\cases{
\pmatrix{\chi_{n-1}(z) & \chi_n(z)} = \pmatrix{\chi_{n-1*}(z) & \chi_{n*}(z)} \Theta(a_n),
& even $n$,
\cr
z\pmatrix{\chi_{n-1*}(z) & \chi_{n*}(z)} = \pmatrix{\chi_{n-1}(z) & \chi_n(z)} \Theta(a_n),
& odd $n$.
}
\eeq

\bc \label{NORMTRUNC-C}

For any sequence $\bsa$ in $\D$, the normal truncations of $C(\bsa)$ on $\ell^2_n$ and
$\ell^{2\bot}_n$ are unitary and they have respectively the form $C(a_1,\dots,a_{n-1},u)$
and $C(\bsa^{(n)};u)$, with $u\in\T$. For both kinds of truncations, the related
projections $Q_n(\bsa;u)$ satisfy $\|Q_n(\bsa;u)\|=\sqrt{1+|u-a_n|^2/\rho_n^2}$. For any
$u\in\T$ the spectrum of $C(a_1,\dots,a_{n-1},u)$ is simple and coincides with the zeros
of the para-orthogonal polynomial $p_n^u$ associated with the measure related to $\bsa$.

\ec

\bpr

Let $\mu$ be the measure on $\T$ whose sequence of Schur parameters is $\bsa$, and let
$(\chi_n)_{n\geq0}$ be the corresponding Laurent orthonormal polynomials. Since
$\spn\{\chi_0,\dots,\chi_{n-1}\}=\PP_{m,n}$, $m=-[\frac{n-1}{2}]$, the unitary
equivalence between $U^\mu$ and $C(\bsa)$ implies that the normal truncations of
$C(\bsa)$ on $\ell^2_n$ ($\ell^{2\bot}_n$) are given by the matrix representations of the
normal truncations of $U^\mu$ on $\PP_{m,n}$ ($\PP_{m,n}^\bot$), when using the basis
$\{\chi_j\}_{j<n}$ ($\{\chi_j\}_{j \geq n}$). So, it just remains to prove that
these representations have the matrix form given in the corollary.

Let $\hat\bsa$ be the sequence obtained from $\bsa$ when substituting $a_n$ by
$u\in\T$. Property (\ref{DEC}) implies that
$C(\hat\bsa) = C(a_1,\dots,a_{n-1},u) \oplus C(\bsa^{(n)};u)$.
If $\Delta:=C(\bsa)-C(\hat\bsa)$ and $\bsc:=(\chi_n)_{n\geq0}$,
(\ref{FAC}) and (\ref{THETA}) lead to
$$
\Delta^t\bsc(z) = z^m(b_np_n^u(z)+d_nq_n^u(z)),
\quad b_n\in\ell^2_n, \quad d_n\in\ell^{2\bot}_n.
$$
Since $C(\bsa)^t\bsc(z)=z\bsc(z)$, if $\bsc_n:=(\chi_j)_{j<n}$ and
$\bsc^{(n)}:=(\chi_j)_{j\geq n}$,
$$
\ba{l}
z\bsc_n(z) = C(a_1,\dots,a_{n-1},u)^t \bsc_n(z) + \beta_n z^m p_n^u(z),
\quad \beta_n\in\C^n,
\smallskip \\
z\bsc^{(n)}(z) = C(\bsa^{(n)};u)^t \bsc^{(n)}(z) + \delta_n z^m q_n^u(z),
\quad \delta_n\in\ell^2.
\ea
$$
This shows that $C(a_1,\dots,a_{n-1},u)$ ($C(\bsa^{(n)};u)$) is the matrix
representation of the normal truncation of $U^\mu$ on $\PP_{m,n}$ ($\PP_{m,n}^\bot$)
along $\spn\{z^mp_n^u\}\oplus\PP_{m,n+1}^\bot$ $(\spn\{z^mq_n^u\} \oplus \PP_{m+1,n-1})$
with respect to $\bsc_n$ ($\bsc^{(n)}$).

\epr

Notice that the normal truncations provided by Theorem \ref{NORMALTRUNC} and Corollary
\ref{NORMTRUNC-C} are always non-orthogonal. These truncations have a remarkable meaning.
Concerning the finite ones, the spectrum of the matrices $C(a_1,\dots,a_n,u)$, $u\in\T$,
provides the nodes of the Szeg\H o quadrature formulas \cite{JoNjTh89} for the measure
$\mu$ corresponding to $C(\bsa)$. In fact, these matrices are the five-diagonal
representations of the unitary multiplication operators related to the finitely supported
measures $\mu_{u,n}$ associated with such quadrature formulas \cite{MIN}. As for the
co-finite truncations, for any $u\in\T$, the matrices $C(\bsa^{(n)};u)$ are unitarily
equivalent to $C(\overline u\bsa^{(n)})$, which are the five-diagonal representations
corresponding to the family $\mu_u^{(n)}$ of Aleksandrov measures related to the
$n$-associated orthogonal polynomials.

When applied to the truncations of $C(\bsa)$ given in the previous corollary, Proposition
\ref{GENTRUNC} states that $\{\supp\,\mu\}'=\{\supp\,\mu_u^{(n)}\}'$ for any $u\in\T$,
and $\supp\,\mu \subset \{z\in\C : \lim_nd(z,\supp\,\mu_{u_n,n})=0\}$ for any sequence
$\bsu=(u_n)_{n\geq1}$ in $\T$. These relations were previously known (see \cite[Theorem
8]{Go02} for the result concerning $\mu_{u,n}$ and \cite{Pe96,GoNe01,Si04-1,Si04-2} for
results related to $\mu_u^{(n)}$). The relevance of our approach is that it proves that
the only possibility of applying Proposition \ref{GENTRUNC} to $C(\bsa)$ necessarily
leads to the measures $\mu_{u,n}$ and $\mu_u^{(n)}$.

%%%%%%%%%%%%%%%%%%%%%%%%%%%%%%%%%%%%%%%%%%%%%%%%%%%%%%%%%%%%%%%%%%%%%%%%%%%%%%%%%%%%%%%%%%%%%%
\section{Co-finite truncations of $C(\bsa)$ and the derived set of the support of the measure}
%%%%%%%%%%%%%%%%%%%%%%%%%%%%%%%%%%%%%%%%%%%%%%%%%%%%%%%%%%%%%%%%%%%%%%%%%%%%%%%%%%%%%%%%%%%%%%

If $\bsa$ is the sequence of Schur parameters for the measure $\mu$ on $\T$, our aim is
to relate $\{\supp\,\mu\}'=\spec_e(C(\bsa))$ and the asymptotic behaviour of $\bsa$ using
Proposition \ref{GENTRUNC}.2 and some results of operator theory. Concerning the
notation, if $T\in\frB(H)$ is normal, we denote by $E_T$ its spectral measure, so that
$$
(Tx,y) = \int_\C \lambda \; d(E_T(\lambda)x,y),
\quad \forall x,y \in H.
$$
In fact, the above expression can be understood as an integral over $\spec(T)$.

As for the results of operator theory that we will apply, we start with a
characterization of the essential spectrum of a normal operator and a lower bound for the
distance from a point to this essential spectrum. In what follows we use the notation
$D_\epsilon(z)$ for an open disk of centre $z\in\C$ and radius $\epsilon>0$.

\bp \label{ESS1}

Let $T\in\frB(H)$ be normal. A point $z\in\C$ lies on $\spec_e(T)$ if and only if
$\gamma(z-T;S)=0$ for any subspace $S \subset H$ with finite co-dimension. Moreover, given
an arbitrary subspace $S \subset H$ with finite co-dimension,
$d(z,\spec_e(T))\geq\gamma(z-T;S)$ for any $z\in\C$.

\ep

\bpr

A point $z\in\C$ belongs to $\spec_e(T)$ if and only if $\rank\,E_T(D_\epsilon(z))$ is
infinite for any $\epsilon>0$. Therefore, if $z\notin\spec_e(T)$,
$S_\epsilon:=E_T(\C\backslash D_\epsilon(z))H$ is a subspace with finite co-dimension for
some $\epsilon>0$. Moreover,
$$
\|(z-T)x\|^2 =
\int_{\C \backslash D_\epsilon(z)} |z-\lambda|^2 \, d(E_T(\lambda)x,x)
\geq \epsilon^2\|x\|^2, \quad \forall x \in S_\epsilon,
$$
which proves that $\gamma(z-T;S_\epsilon)\geq\epsilon>0$.

If, on the contrary, $z\in\spec_e(T)$, $S'_\epsilon:=E_T(D_\epsilon(z))H$ has infinite
dimension for any $\epsilon>0$ and, thus, given a subspace $S \subset H$ with finite
co-dimension there always exists a non-null vector $x_\epsilon \in S \cap S'_\epsilon$.
Then,
$$
\|(z-T)x_\epsilon\|^2 =
\int_{D_\epsilon(z)} |z-\lambda|^2 \, d(E_T(\lambda)x_\epsilon,x_\epsilon)
\leq \epsilon^2\|x_\epsilon\|^2.
$$
Since $\epsilon$ is arbitrary, we conclude that $\gamma(z-T;S)=0$.
Moreover, given $z\in\C$ and a subspace $S \subset H$ with finite co-dimension,
$$
\gamma(z-T;S)\leq|z-w|+\gamma(w-T;S)=|z-w|, \quad \forall w\in\spec_e(T),
$$
which proves that $d(z,\spec_e(T))\geq\gamma(z-T;S)$.

\epr

Concerning perturbative results, if $T_0,T\in\frB(H)$ and $T_0$ is normal, it is known
that $\spec(T)\subset\{z\in\C:d(z,\spec(T_0))\leq\|T-T_0\|\}$. The next result is the
analogue for the essential spectrum of normal operators.

\bp \label{ESS2}

If $T_0,T\in\frB(H)$ are normal, for any subspace $S \subset H$ with finite co-dimension
$$
\spec_e(T) \subset \{z\in\C : d(z,\spec_e(T_0))\leq\|T-T_0\|_S\}.
$$

\ep

\bpr

Suppose that $d(z,\spec_e(T_0))>\|T-T_0\|_S$. Consider a real number $\epsilon$ such that
$d(z,\spec_e(T_0))>\epsilon>\|T-T_0\|_S$. The subspace
$S_\epsilon:=E_{T_0}(\C \backslash D_\epsilon(z))H$ has finite co-dimension and, similarly
to the proof of the previous proposition, $\gamma(z-T_0;S_\epsilon)\geq\epsilon$.
$S'_\epsilon := S \cap S_\epsilon$ has also finite co-dimension and
$$
\gamma(z-T;S'_\epsilon) \geq \gamma(z-T_0;S'_\epsilon)-\|(T-T_0)\|_{S'_\epsilon}
\geq \gamma(z-T_0;S_\epsilon)-\|(T-T_0)\|_S.
$$
Therefore, $\gamma(z-T;S'_\epsilon)\geq\epsilon-\|T-T_0\|_S>0$ and, from Proposition \ref{ESS1},
$z\notin\spec_e(T)$.

\epr

When applying to $C(\bsa)$ and its co-finite normal truncations $C(\bsa^{(n)};u)$,
$u\in\T$, propositions \ref{GENTRUNC}.2, \ref{ESS1} and \ref{ESS2} give the following
result. For convenience, given an operator $T\in\frB(\ell^2)$, we use the notation
$\|T\|_n:=\|T\|_{\ell^{2\bot}_n}$, $\gamma_n(T):=\gamma(T;\ell^{2\bot}_n)$, $n\geq1$, and
$\|T\|_0:=\|T\|$, $\gamma_0(T):=\gamma(T)$.

\bt \label{ESS-C}

Let $\bsa$ be the sequence of Schur parameters of a measure $\mu$ on $\T$, $\bsb$ a
sequence in $\overline\D$, $u\in\T$, $m\geq0$ and $z\in\C$.

\be

\item
$\ds d(z,\{\supp\,\mu\}') \geq \, \sup_{n\geq0}\gamma_m(z-C(\overline u\bsa^{(n)}))$.

\item
$\ds \inf_{n\geq0}\|C(\overline u\bsa^{(n)})-C(\bsb^{(n)})\|_m < \, d(z,\spec_e(C(\bsb)))
\;\Rightarrow\; z\notin\{\supp\,\mu\}'$.

\ee

\et

\bpr

From Proposition \ref{GENTRUNC}.2,
$\spec_e(C(\overline u\bsa^{(n)}))=\spec_e(C(\bsa))=\{\supp\,\mu\}'$.
So, a direct application of Proposition \ref{ESS1} gives
$d(z,\{\supp\,\mu\}') \geq \gamma_m(z-C(\overline u\bsa^{(n)}))$
for any $m$, which proves the first statement.

Let us suppose now that
$\inf_{n\geq0}\|C(\overline u\bsa^{(n)})-C(\bsb^{(n)})\|_m < d(z,\spec_e(C(\bsb)))$.
Then, for some $n$,
$\|C(\overline u\bsa^{(n)})-C(\bsb^{(n)})\|_m < d(z,\spec_e(C(\bsb)))$.
Since Proposition \ref{GENTRUNC}.2 implies that
$\spec_e(C(\bsb))=\spec_e(C(\bsb^{(n)}))$, it follows from Proposition \ref{ESS2} that
$z\notin\spec_e(C(\overline u\bsa^{(n)}))$. Hence, using again Proposition
\ref{GENTRUNC}.2, we find that $z\notin\spec_e(C(\bsa))=\{\supp\,\mu\}'$.

\epr

For the application of the preceding propositions we have to obtain lower bounds for
$\gamma_m(C(\bsa)-C(\bsb))$ and upper bounds for $\|C(\bsa)-C(\bsb)\|_m$, $\bsa$ and
$\bsb$ being sequences in $\overline\D$. This is all we need since $C(\bsb)=z$ for
$b_n=(-z)^n$ with $z\in\T$. Notice that
$$
\left\|\pmatrix{-\alpha & \overline \beta \cr \beta & \overline \alpha}
\pmatrix{x \cr y}\right\| =
\left\|\pmatrix{\alpha \cr \beta}\right\|
\left\|\pmatrix{x \cr y}\right\|,
\qquad \pmatrix{\alpha \cr \beta},\pmatrix{x \cr y}\in\C^2.
$$
Hence, from the equality
$$
C(\bsa)-C(\bsb) =
\left(C_o(\bsa)-C_o(\bsb)\right)C_e(\bsa) + C_o(\bsb)(C_e(\bsa)-C_e(\bsb)),
$$
we get
\beq \label{NORM1}
\ba{l}
\gamma_m(C(\bsa)-C(\bsb)) \; \geq
\ds \inf_{\scriptsize \ba{c} j\geq m-1 \\ \odd\,(\even)\,j \ea}
\kern-8pt k(a_j,b_j) \;\;
- \sup_{\scriptsize \ba{c} j\geq m-1 \\ \even\,(\odd)\,j \ea}
\kern-8pt k(a_j,b_j),
\smallskip \\
\|C(\bsa)-C(\bsb)\|_m \, \leq
\ds \sup_{\scriptsize \ba{c} j\geq m-1 \\ \odd\,j \ea}
\kern-5pt k(a_j,b_j) \;\;
+ \sup_{\scriptsize \ba{c} j\geq m-1 \\ \even\,j \ea}
\kern-5pt k(a_j,b_j),
\ea
\eeq
where $k(x_1,x_2):=\gamma(\Theta(x_1)-\Theta(x_2))=\|\Theta(x_1)-\Theta(x_2)\|$,
that is,
\beq \label{k}
k(x_1,x_2)=\sqrt{|x_1-x_2|^2+|y_1-y_2|^2}, \quad |x_j|\leq1, \quad y_j:=\sqrt{1-|x_j|^2}.
\eeq

Equipped with these results, we can apply Theorem \ref{ESS-C} in different ways to get
information about the derived set of the support of a measure on $\T$ from the behaviour
of the related Schur parameters. The first example of this is the next theorem.
In what follows, given $z,w\in\C$, $w=e^{i\omega}z$, $\omega\in[0,2\pi)$, we denote
$(z,w):=\{e^{i\theta}z:\theta\in(0,\omega)\}$,
$[z,w]:=\{e^{i\theta}z:\theta\in[0,\omega]\}$.
Also, for any $z\in\C$ and $\alpha\in[0,\pi]$,
$\triangle_\alpha(z):=[e^{i\alpha}z,e^{-i\alpha}z]$,
$\Gamma_\alpha(z):=(e^{-i\alpha}z,e^{i\alpha}z)$ and
$\overline\Gamma_\alpha(z):=[e^{-i\alpha}z,e^{i\alpha}z]$.
Besides, for any sequence $\bsa$ in $\C$ and any $z\in\C$ we define
$\bsa(z):=(a_n(z))_{n\geq1}$ by $a_n(z):=\overline z^na_n$.

\bt \label{DERIVED1}

Let $\mu$ be a measure on $\T$ with a sequence $\bsa$ of Schur parameters. Assume that
for some $\lambda\in\T$ the limit points of the odd and even subsequences of
$\bsa(-\lambda)$ are separated by a band
$\cB(u,\alpha_1,\alpha_2):=\{z\in\C:\cos\alpha_2<\re(\overline uz)<\cos\alpha_1\}$,
$u\in\T$, $0\leq\alpha_1<\alpha_2\leq\pi$. Then,
$$
\{\supp\,\mu\}'\subset\triangle_\alpha(\lambda),
\qquad \sin\frac{\alpha}{2}=
\max\{\sin\frac{\alpha_2}{2}-\sin\frac{\alpha_1}{2},
\cos\frac{\alpha_1}{2}-\cos\frac{\alpha_2}{2}\}.
$$

\et

\bpr

Let $\bsb$ be defined by $b_n:=(-\lambda)^n$, $\lambda\in\T$. $C(\bsb^{(n)})$ is diagonal
with diagonal elements equal to $\lambda$, except the first one that is $(-1)^n\lambda^{n+1}$.
Therefore, for $m\geq1$,
$\gamma_m(\lambda-C(\overline u\bsa^{(n)}))=\gamma_m(C(\bsb^{(n)})-C(\overline u\bsa^{(n)}))$
and, from (\ref{NORM1}), we get
$$
\ba{l}
\kern-4pt
\ds\sup_{n\geq0}\gamma_m(\lambda-C(\overline u\bsa^{(n)})) \geq
\sup_{n\geq0} \kern-7pt \inf_{\scriptsize \ba{c} \odd\,(\even)\,j \\ j\geq n+m-1 \ea}
\kern-15pt k(\overline ua_j,b_j) -
\inf_{n\geq0} \kern-7pt \sup_{\scriptsize \ba{c} \odd\,(\even)\,j \\ j\geq n+m-1 \ea}
\kern-15pt k(\overline ua_j,b_j) =
\smallskip \\ \kern35pt
\ds = \kern-5pt \limi_{\odd\,(\even)\,n}
\kern-5pt \sqrt{2(1-\re(\overline ua_n(-\lambda)))} \,-
\kern-5pt \lims_{\even\,(\odd)\,n}
\kern-5pt \sqrt{2(1-\re(\overline ua_n(-\lambda)))}.
\ea
$$

Assume that the limit points of the even and odd subsequences of $\bsa(-\lambda)$ are
separated by the band $\cB(u,\alpha_1,\alpha_2)$. Then,
$$
\limi_{\even\,(\odd)\,n} \re(\overline ua_n(-\lambda)) \geq \cos\alpha_1,
\quad
\lims_{\odd\,(\even)\,n} \re(\overline ua_n(-\lambda)) \leq \cos\alpha_2,
$$
which gives
$$
\sup_{n\geq0}\gamma_m(\lambda-C(\overline u\bsa^{(n)})) \geq
2\sin\frac{\alpha_2}{2}-2\sin\frac{\alpha_1}{2}.
$$
Since the limit points of the even and odd subsequences of $\bsa(-\lambda)$ are
separated by the band $\cB(-u,\pi-\alpha_2,\pi-\alpha_1)$ too, we also get
$$
\sup_{n\geq0}\gamma_m(\lambda-C(-\overline u\bsa^{(n)})) \geq
2\cos\frac{\alpha_1}{2}-2\cos\frac{\alpha_2}{2}.
$$
So, if we define $\alpha\in(0,\pi]$ by
$\sin\frac{\alpha}{2}=
\max\{\sin\frac{\alpha_2}{2}-\sin\frac{\alpha_1}{2},
\cos\frac{\alpha_1}{2}-\cos\frac{\alpha_2}{2}\}$,
Proposition \ref{ESS-C}.1 gives
$d(\lambda,\{\supp\,\mu\}') \geq 2\sin\frac{\alpha}{2}$, that is,
$\{\supp\,\mu\}'\subset\triangle_\alpha(\lambda)$.

\epr

Given $\lambda\in\T$, the above theorem ensures that $\lambda\notin\{\supp\,\mu\}'$ if
the limit points of the odd and even subsequences of $\bsa(-\lambda)$ can be separated by
a straight line. In particular, $\lambda\notin\{\supp\,\mu\}'$ if the limit points of
$\bsa(\lambda)$ lie on an open half-plane whose boundary contains the origin, because
then, the limit points of the odd and even subsequences of $\bsa(-\lambda)$ are separated
by the straight line that limits such half-plane. In fact, we have the following
immediate consequence of Theorem \ref{DERIVED1}.

\bc \label{DERIVED11}

Let $\mu$ be a measure on $\T$ with a sequence $\bsa$ of Schur parameters. Assume that
for some $\lambda\in\T$ the limit points of $\bsa(\lambda)$ lie on
$\cD(u,\alpha_0):=\{z\in\C:\re(\overline uz)\geq\cos\alpha_0\}$, $u\in\T$,
$0\leq\alpha_0<\pi/2$. Then,
$$
\{\supp\,\mu\}'\subset\triangle_\alpha(\lambda),
\qquad \cos\frac{\alpha}{2}=\sqrt{\sin\alpha_0}.
$$

\ec

\bpr

Under the assumptions of the corollary, the limit points of the odd and even subsequences of
$\bsa(-\lambda)$ are separated by the band $\cB(u,\alpha_0,\pi-\alpha_0)$. So, the direct
application of Theorem \ref{DERIVED1} proves that
$\{\supp\,\mu\}'\subset\triangle_\alpha(\lambda)$
with $\alpha\in(0,\pi]$ given by
$\sin\frac{\alpha}{2}=\cos\frac{\alpha_0}{2}-\sin\frac{\alpha_0}{2}$,
that is,
$\cos^2\frac{\alpha}{2}=\sin\alpha_0$.

\epr

\br \label{remDER1}

If $\mu$ is the measure related to the sequence $\bsa$ of Schur parameters, the measure
obtained by rotating $\mu$ through an angle $\theta$ is associated with the sequence
$\bsa(e^{i\theta})$. Thus, Proposition \ref{DERIVED1} and Corollary \ref{DERIVED11} are
just the rotated versions of the following basic statements:
\be
\item If the limit points of the odd and even subsequences of $\bsa$ are separated by a band
obtained by a rotation of $\cos\alpha_1<\re(z)<\cos\alpha_2$, then
$\{\supp\,\mu\}'\subset\triangle_\alpha(-1)$ with
$\sin\frac{\alpha}{2}=
\max\{\sin\frac{\alpha_2}{2}-\sin\frac{\alpha_1}{2},
\cos\frac{\alpha_1}{2}-\cos\frac{\alpha_2}{2}\}$.
\item If the limit points of $\bsa$ lie on a domain obtained by a rotation of
$\re(z)\geq\cos\alpha_0>0$, then
$\{\supp\,\mu\}'\subset\triangle_\alpha(1)$ with $\cos\frac{\alpha}{2}=\sqrt{\sin\alpha_0}$.
\ee

\er

Let us show an example of application of the previous results. In what follows we denote by
$\frak L\{\bsa\}$ the set of limit points of a sequence $\bsa$ in $\C$.

\bex \label{EXDER1}
$\frak L\{\bsa(\lambda)\}=\{a,b\}$, $\lambda\in\T$, $a\neq b$, $a,b\neq0$.

\smallskip

\noindent Let $\bsa$ be the Schur parameters of a measure $\mu$ on $\T$. With the help of
the previous results we can get information about the case in which we just know that
$\bsa(\lambda)$ has two different subsequential limit points $a$, $b$, no matter from
which subsequence. Suppose that $0<|a|\leq|b|$ and let
$\frac{b-a}{|b-a|}=\frac{a}{|a|}e^{i\theta}$, $\theta\in(-\pi,\pi]$. Then,
$\{a,b\}\subset\cD(u,\alpha_0)$ where
$$
\ba{l}
\ds u=\frac{a}{|a|}, \quad \cos\alpha_0=|a|,
\quad {\rm if} \kern5pt |\theta|\leq\frac{\pi}{2},
\smallskip \\
\ds u=-\sg(\theta)\,i\frac{b-a}{|b-a|},
\quad \cos\alpha_0=|a|\sin|\theta|, \quad {\rm if} \kern5pt |\theta|>\frac{\pi}{2}.
\ea
$$
Therefore, using Corollary \ref{DERIVED11} we find that, if $\theta\neq\pi$ (which means
that $\frac{b}{|b|}\neq-\frac{a}{|a|}$),
$\{\supp\,\mu\}'\subset\triangle_\alpha(\lambda)$ with
$$
\cos\frac{\alpha}{2}=\cases{
\sqrt[4]{1-|a|^2} & if $|\theta|\leq\frac{\pi}{2}$,
\smallskip \cr
\sqrt[4]{1-|a|^2\sin^2\theta} & if $|\theta|>\frac{\pi}{2}$.}
$$

\eex

\medskip

The next results use the second statement of Theorem \ref{ESS-C}. This requires the
comparison of the matrix $C(\bsa)$ related to a measure on $\T$ with another matrix
$C(\bsb)$ with known essential spectrum. The simplest case where the essential spectrum
of $C(\bsb)$ is known is when it is a diagonal matrix, which means that $\bsb$ is a
sequence in $\T$. Applying Theorem \ref{ESS-C}.2 to the comparison between $C(\bsa)$ and
$C(\bsb)$ with a suitable choice for $\bsb$ in $\T$, we get the following result.

\bp \label{DERIVED2}

Let $a_n=|a_n|u_n \, (u_n\in\T)$ be the Schur parameters of a measure $\mu$ on $\T$, and
assume that $c(\bsa) := \min\{c_1(\bsa),c_2(\bsa)\} < 1$, where
$$
\ba{l} \ds c_1(\bsa) :=
\frac{1}{2}\lims_n\left(||a_{n+1}|-|a_n|| + \rho_n + \rho_{n+1}\right),
\smallskip \\
\ds c_2(\bsa) := \lims_{\odd\,n}\sqrt{\frac{1-|a_n|}{2}} +
\lims_{\even\,n}\sqrt{\frac{1-|a_n|}{2}}.
\ea
$$
Then,
$$
\{\supp\,\mu\}' \subset
\kern-10pt\bigcup_{\lambda\in\frak L\{\overline u_nu_{n+1}\}}\kern-15pt
\triangle_\alpha(\lambda), \qquad \cos\frac{\alpha}{2}=c(\bsa).
$$

\ep

\bpr

Let us define the sequence $\bsb$ by $b_n:=u_n$. Since $u_n\in\T$, $C(\bsb)$ is the
diagonal matrix
\beq \label{DIAGONAL}
C(\bsb) = -\pmatrix{
u_1 & & & \cr
 & \kern-9pt \overline u_1u_2 & & \cr
 & & \kern-9pt \overline u_2u_3 & \cr
 & & & \kern-9pt \ddots
},
\eeq
and, hence, $\spec_e(C(\bsb))=-\frak L\{\overline u_nu_{n+1}\}$. Using (\ref{NORM1})
we get
$$
\inf_{n\geq0}\|C(\bsa^{(n)})-C(\bsb^{(n)})\| \leq 2c_2(\bsa).
$$

We can find another upper bound for $\inf_{n\geq0}\|C(\bsa^{(n)})-C(\bsb^{(n)})\|$ in the
following way. The factorizations $C(\bsa)-C(\bsb) =
C_o(\bsa)\left(C_e(\bsa)-C_o^*(\bsa)C(\bsb)\right)$ and
$C_e(\bsa)-C_o^*(\bsa)C(\bsb)=A(\bsa)B(\bsb)$, where
$$
\ba{l}
A(\bsa) := \pmatrix{
1-|a_1| & \rho_1\overline u_1 & & & & \cr
\rho_1u_1 & |a_1|-|a_2| & \rho_2u_2 & & & \cr
 & \rho_2\overline u_2 & |a_2|-|a_3| & \rho_3\overline u_3 & & \cr
 & & \rho_3u_3 & |a_3|-|a_4| & \kern12pt \rho_4u_4 & \cr
 & & & \kern-60pt\ddots & \kern-45pt\ddots & \ddots
},
\smallskip \\
B(\bsb) := \pmatrix{
1 & & & & & \cr
 & \kern-3pt u_2 & & & & \cr
 & & \kern-3pt \overline u_2 & & & \cr
 & & & \kern-3pt u_4 & & \cr
 & & & & \kern-3pt \overline u_4 & \cr
 & & & & & \kern-5pt \ddots
}, \ea
$$
together with the fact that $A(\bsa)$ is unitarily equivalent (by a diagonal
transformation) to the Jacobi matrix
$$
J(\bsa) := \pmatrix{
1-|a_1| & \rho_1 & & & & \cr
\rho_1 & |a_1|-|a_2| & \rho_2 & & & \cr
 & \rho_2 & |a_2|-|a_3| & \rho_3 & & \cr
 & & \rho_3 & |a_3|-|a_4| & \kern15pt \rho_4 & \cr
 & & & \kern-60pt\ddots & \kern-45pt\ddots & \ddots
},
$$
show that $\|C(\bsa)-C(\bsb)\|=\|J(\bsa)\|$.

Since $J(\bsa)$ is a bounded Jacobi matrix, it defines a self-adjoint operator on
$\ell^2$. So, it follows from Proposition \ref{GENTRUNC}.1 that
$\spec(J(\bsa)) \subset \{z\in\C : \lim_nd(z,\spec(J(a_1,\dots,a_n)))=0\}$,
where $J(a_1,\dots,a_n)$ is the principal submatrix of $J(\bsa)$ of order $n$.
A direct application of Gershgorin theorem shows that
$$
\spec(J(a_1,\dots,a_n)) \subset
\{z\in\C : |z| \leq \max_{j=1}^n(||a_{j-1}|-|a_j||+\rho_{j-1}+\rho_j)\},
$$
where $a_0=1$ and $\rho_0=0$. Therefore,
$$
\|C(\bsa)-C(\bsb)\|=\|J(\bsa)\|=\kern-5pt\max_{\lambda\in\spec(J(\bsa))}|\lambda |
\leq \sup_{j\geq0}(||a_j|-|a_{j+1}||+\rho_j+\rho_{j+1}).
$$

Taking into account that
$|1-|a_{n+1}||+\rho_{n+1} \leq ||a_n|-|a_{n+1}||+\rho_n+\rho_{n+1}$,
a similar reasoning leads to
$$
\|C(\bsa^{(n)})-C(\bsb^{(n)})\|=\|J(\bsa^{(n)})\| \leq
\sup_{j \geq n}(||a_n|-|a_{n+1}||+\rho_n+\rho_{n+1}).
$$
From this inequality we get
$$
\inf_{n\geq0}\|C(\bsa^{(n)})-C(\bsb^{(n)})\| \leq 2c_1(\bsa).
$$

Summarizing,
$\inf_{n\geq0}\|C(\bsa^{(n)})-C(\bsb^{(n)})\| \leq 2c(\bsa)$.
Thus, Theorem \ref{ESS-C} implies that a point $z\in\T$ is outside $\{\supp\,\mu\}'$
if $d(z,-\frak L\{\overline u_nu_{n+1}\}) > 2c(\bsa)$, which can be satisfied only
if $c(\bsa)<1$. In such a case we can write $c(\bsa)=\cos\frac{\alpha}{2}$,
$\alpha\in(0,\pi]$, and
$$
\ba{l}
\{\supp\,\mu\}' \subset
\{z\in\T : d(z,-\frak L\{\overline u_nu_{n+1}\}) \leq 2\cos\frac{\alpha}{2}\}
= \kern-10pt\ds\bigcup_{\lambda\in\frak L\{\overline u_nu_{n+1}\}}\kern-15pt
\triangle_\alpha(\lambda).
\ea
$$

\epr

An immediate corollary of this proposition is a condition for the Schur parameters which
ensures that a certain arc of $\T$ is outside the derived set of the support of the measure.

\bc \label{DERIVED21}

Under the conditions of Proposition \ref{DERIVED2},
$$
\frak L\{\overline u_nu_{n+1}\} \subset \overline\Gamma_\zeta(\lambda),
\quad 0\leq\zeta<\alpha \quad\Rightarrow\quad
\{\supp\,\mu\}' \subset \triangle_{\alpha-\zeta}(\lambda).
$$

\ec

A remarkable consequence of Proposition \ref{DERIVED2} is obtained when studying
measures $\mu$ in the class $\lim_n|a_{n+1}/a_n|=l$. Notice that
$l\leq1$ because $\bsa$ is bounded. When $l<1$, $\lim_n|a_n|=0$ and, thus,
$\supp\,\mu=\T$ as a consequence of Weyl's theorem, since, for $b_n=0$, $C(\bsa)-C(\bsb)$
is compact. So, concerning $\supp\,\mu$, the only non trivial case is $l=1$.

The condition $\lim_n|a_{n+1}/a_n|=1$ covers the case $\lim_n|a_n|=1$, for which
$C(\bsa)$ differs in a compact perturbation from a diagonal matrix with diagonal elements
$-\overline a_{n-1}a_n$. Therefore, in this case,
$\{\supp\,\mu\}'=-\frak L\{a_{n+1}/a_n\}$
(see \cite{Go00a} for a similar argument using Hessenberg representations).

$\lim_n|a_{n+1}/a_n|=1$ is also verified when $\lim_na_n=a\in\D\backslash\{0\}$, which
implies that $C(\bsa)-C(\bsb)$ is compact for $b_n=a$. So, under this condition,
$\{\supp\,\mu\}'=\triangle_\alpha(1)$, $\sin\frac{\alpha}{2}=|a|$, as in the Geronimus
case corresponding to constant Schur parameters equal to $a$ \cite{Ge41,GoNeAs95}. A
bigger class is $\lim_na_{n+1}/a_n=\lambda\in\T$, $\lim_n|a_n|=r\in(0,1)$, which is known
as the L\'opez class \cite{BaLo99}. It is a particular case of $\lim_n|a_{n+1}/a_n|=1$
too. In the L\'opez class, $C(\bsa)$ is unitarily equivalent to a matrix obtained by a
compact perturbation of the matrix $C(\bsb)$ associated with the rotated Geronimus case
$b_n=\lambda^nr$ (see \cite[Chapter 4]{Si04-1}) and, therefore,
$\{\supp\,\mu\}'=\triangle_\alpha(\lambda)$, $\sin\frac{\alpha}{2}=r$.

All these results are known. We mention them to help understand to what extent the next
theorem is an extension of them. Notice that, not only the L\'opez class, but also the two
conditions that define this class are separately particular cases of
$\lim_n|a_{n+1}/a_n|=1$.

\bt \label{DERIVED22}

If $\bsa$ is the sequence of Schur parameters associated with a measure $\mu$ on $\T$,
$$
\ds \lim_n\left|\frac{a_{n+1}}{a_n}\right|=1
\quad\Rightarrow\quad
\{\supp\,\mu\}' \subset
\kern-7pt\bigcup_{\lambda\in\frak L\left\{\frac{a_{n+1}}{a_n}\right\}}\kern-12pt
\triangle_\alpha(\lambda),
\quad \sin\frac{\alpha}{2}=\limi_n|a_n|.
$$

\et

\bpr

The result is trivial when $\limi_n|a_n|=0$, so we just have to consider the case
$\limi_n|a_n|\neq0$. The result follows from Proposition \ref{DERIVED2}, taking into
account that, if $\lim_n|a_{n+1}/a_n|=1$ and $\lim_n|a_n|\neq0$,
$c(\bsa)=c_1(\bsa)=\lims_n\rho_n<1$ and
$\frak L\{\overline u_nu_{n+1}\} = \frak L\{a_{n+1}/a_n\}$.

\epr

Notice that the bounds provided by Theorem \ref{DERIVED22} give the exact location of
$\{\supp\,\mu\}'$ for the L\'opez class and also for the case $\lim_n|a_n|=1$.
A direct consequence of this theorem is a result for the class defined only by the first
of the L\'opez conditions.

\bc \label{DERIVED221}

If $\bsa$ is the sequence of Schur parameters associated with a measure $\mu$ on $\T$,
$$
\lim_n\frac{a_{n+1}}{a_n}=\lambda\in\T \quad\Rightarrow\quad
\{\supp\,\mu\}' \subset \triangle_\alpha(\lambda),
\quad \sin\frac{\alpha}{2}=\limi_n|a_n|.
$$

\ec

Theorem \ref{DERIVED22} can be used to supplement the conclusions of Theorem \ref{DERIVED1}
and Corollary \ref{DERIVED11}. Let us see an example.

\bex \label{EXDER2}
$\frak L\{\bsa(\lambda)\}=\{a,b\}$, $\lambda\in\T$, $a\neq b$, $|a|=|b|\neq0$.

\smallskip

\noindent When $\bsa(\lambda)$ has two different limit points $a$, $b$, using Corollary
\ref{DERIVED11} we got information about an arc centred at $\lambda$ which is free of
$\{\supp\,\mu\}'$. Theorem \ref{DERIVED22} helps us find other arcs outside
$\{\supp\,\mu\}'$ when $|a|=|b|\neq0$ since, in this case, $\lim_n|\frac{a_{n+1}}{a_n}|=1$.

Without loss of generality we can suppose $b=ae^{i\zeta}$, $\zeta\in(0,\pi]$, so,
$\frak L\{\frac{a_{n+1}}{a_n}\}\subset\{\lambda,\lambda e^{i\zeta},\lambda e^{-i\zeta}\}$.
Hence, Theorem \ref{DERIVED22} gives three possible arcs lying on
$\T\backslash\{\supp\,\mu\}'$, centred at $\lambda$ and $-\lambda e^{\pm i\frac{\zeta}{2}}$.
More precisely,
$$
\ba{l}
\ds \sin\frac{\zeta}{2}<|a| \kern3pt\Rightarrow\kern3pt
\{\supp\,\mu\}'\subset\triangle_{\alpha_1}(\lambda),
\kern4pt \alpha_1=\alpha-\zeta,
\smallskip \\
\ds \cos\frac{\zeta}{4}<|a| \kern3pt\Rightarrow\kern3pt
\{\supp\,\mu\}'\subset
\triangle_{\alpha_2}(-\lambda e^{i\frac{\zeta}{2}})
\cap \triangle_{\alpha_2}(-\lambda e^{-i\frac{\zeta}{2}}),
\kern4pt \alpha_2=\alpha+\frac{\zeta}{2}-\pi,
\ea
$$
where $\alpha\in(0,\pi)$ is given by $\sin\frac{\alpha}{2}=|a|$.

\eex

\medskip

Concerning Proposition \ref{DERIVED2}, when $\lim_n(|a_{n+1}|-|a_n|)=0$, $c(\bsa)=c_1(\bsa)$
since $\rho_n<\sqrt{2(1-|a_n|)}$. On the contrary, $c(\bsa)=c_2(\bsa)$ if
$\lim_n|a_{2n-1}|=1$ or $\lim_n|a_{2n}|=1$, due to the inequality $1-|a_n|+\rho_n >
\sqrt{2(1-|a_n|)}$. So, in the last case, $c(\bsa)=\lims_n\sqrt{(1-|a_n|)/2}$ and we get
the following corollary, which is a generalization of a result given in \cite[pg. 72]{Go00a}.

\bc \label{DERIVED23}

Let $a_n=|a_n|u_n \, (u_n\in\T)$ be the Schur parameters of a measure $\mu$ on $\T$.
If $\lim_n|a_{2n-1}|=1$ or $\lim_n|a_{2n}|=1$, then,
$$
\{\supp\,\mu\}' \subset
\kern-10pt\bigcup_{\lambda\in\frak L\{\overline u_nu_{n+1}\}}\kern-15pt
\triangle_\alpha(\lambda),
\qquad \cos\alpha=-\limi_n|a_n|.
$$

\ec

The rotated Geronimus case $C(\bsb)$, $b_n=\lambda^na$, $a\in\D\backslash\{0\}$,
$\lambda\in\T$, can also be used for the comparison in Theorem \ref{ESS-C}, since we know
that $\spec_e(C(\bsb))=\triangle_\alpha(\lambda)$, $\sin\frac{\alpha}{2}=|a|$. In fact,
the consequences of this comparison are a particularization of a more general result
concerning the comparison with the case $b_{2n-1}=\lambda^{2n-1}a_o$,
$b_{2n}=\lambda^{2n}a_e$, where $a_o,a_e\in\overline\D$. In this case it is known that
$\spec_e(C(\bsb))= \triangle_{\alpha_+}(\lambda)\cap\triangle_{\alpha_-}(-\lambda)$ where
$\alpha_\pm\in[0,\pi]$ are given by
\beq \label{PER}
\cos\alpha_\pm=\rho_o\rho_e\mp\re(\overline a_oa_e), \quad \rho_i:=\sqrt{1-|a_i|^2},
\quad i=o,e
\eeq
(see \cite{Ge44,PeSt96}). That is, $\spec_e(C(\bsb))$ has two connected components except
for the cases $a_o=\pm a_e$ which correspond to only one connected component.

\bp \label{DERIVED3}

Let $\bsa$ be the sequence of Schur parameters of a measure $\mu$ on $\T$ and, for
$a_o,a_e\in\overline\D$ and $\lambda\in\T$, let us define
$$
\ds s(\bsa) :=
\frac{1}{2}\left\{\lims_{\odd\,n}k(a_n(\lambda),a_o) +
\lims_{\even\,n}k(a_n(\lambda),a_e)\right\},
$$
with $k(\cdot,\cdot)$ given in (\ref{k}). Then, defining $\alpha_\pm$ as in (\ref{PER}),
$$
s(\bsa)<\sin\frac{\alpha_\pm}{2} \kern7pt\Rightarrow\kern7pt
\{\supp\,\mu\}' \subset \triangle_{\alpha_\pm-\zeta}(\pm\lambda),
\kern7pt \sin\frac{\zeta}{2}=s(\bsa), \kern7pt 0\leq\zeta<\alpha_\pm.
$$

\ep

\bpr

Let us consider the sequence $b_{2n-1}:=\lambda^{2n-1}a_o$, $b_{2n}:=\lambda^{2n}a_e$.
Using (\ref{NORM1}) we find that
$$
\inf_{n\geq0}\|C(\bsa^{(n)})-C(\bsb^{(n)})\| \leq 2s(\bsa).
$$
When $s(\bsa)<\sin\frac{\alpha_\pm}{2}$ we can write $s(\bsa)=\sin\frac{\zeta}{2}$,
$\zeta\in[0,\alpha_\pm)$, and, since
$\spec_e(C(\bsb))=\triangle_{\alpha_+}(\lambda)\cap\triangle_{\alpha_-}(-\lambda)$,
Theorem \ref{ESS-C} implies that
$$
\{\supp\,\mu\}' \subset \{z\in\T : d(z,\triangle_{\alpha_\pm}(\pm\lambda)) \leq
2\sin\frac{\zeta}{2}\} = \triangle_{\alpha_\pm-\zeta}(\pm\lambda).
$$

\epr

\br \label{RS}

With the notation of (\ref{k}), from $y_1^2-y_2^2=|x_2|^2-|x_1|^2$, we get
$$
|y_1-y_2|\leq\frac{|x_1|+|x_2|}{y_1+y_2}|x_1-x_2|,
$$
hence,
$$
k(x_1,x_2) \leq \frac{\sqrt{2(1+|x_1||x_2|+y_1y_2)}}{y_1+y_2}|x_1-x_2| \leq
\frac{2}{y_1+y_2}|x_1-x_2|.
$$
Therefore,
$$
k(a_n(\lambda),a_i) < \frac{2}{\rho_i}|a_n(\lambda)-a_i|,
$$
and the conclusions of Proposition \ref{DERIVED3} hold if
$$
\frac{1}{\rho_o}\lims_{\odd\,n}|a_n(\lambda)-a_o|
+\frac{1}{\rho_e}\lims_{\even\,n}|a_n(\lambda)-a_e|
<\sin\frac{\alpha_\pm}{2}.
$$

\er

\medskip

Another subclass of $\lim_n|a_{n+1}/a_n|=1$ is given by the second L\'opez condition,
$\lim_n|a_n|=r\in(0,1)$. In this subclass, Proposition \ref{DERIVED3} supplements Theorem
\ref{DERIVED22} with the result that we present below. Notice that $\lim_n|a_n|=r$ if and
only if, for $\lambda\in\T$, the limit points of $\bsa(\lambda)$ lie on the circle
$\{z\in\C : |z|=r\}$. The key idea is that the knowledge of the arcs of this circle in
which the limit points of $\bsa(\lambda)$ lie, gives information about the arc of $\T$
around $\lambda$ that is free of $\{\supp\,\mu\}'$. We will state a more general result
that deals with the case $\lim_n|a_{2n-1}(\lambda)|=r_o$, $\lim_n|a_{2n}(\lambda)|=r_e$.

\bt \label{DERIVED31}

Let $\bsa$ be the sequence of Schur parameters of a measure $\mu$ on $\T$.
Assume that for some $\lambda\in\T$ the limit points of the odd and even subsequences
of $\bsa(\lambda)$ lie on $\overline\Gamma_{\xi_o}(a_o)$ and $\overline\Gamma_{\xi_e}(a_e)$
respectively. Then, if $\alpha_\pm$ is given in (\ref{PER}),
$$
s \kern-1pt :=|a_o|\sin\frac{\xi_o}{2}+|a_e|\sin\frac{\xi_e}{2}<\sin\frac{\alpha_\pm}{2}
\kern3pt\Rightarrow\kern3pt \{\supp\,\mu\}'
\subset \triangle_{\alpha_\pm-\zeta}(\pm\lambda), \kern4pt \sin\frac{\zeta}{2}=s.
$$

\et

\bpr

Since $\lim_n|a_{2n-1}|=|a_o|$ and $\lim_n|a_{2n}|=|a_e|$,
$$
s(\bsa)=\frac{1}{2}\left\{\lims_{\odd\,n}|a_n(\lambda)-a_o|
+\lims_{\even\,n}|a_n(\lambda)-a_e|\right\}.
$$
The statement follows from Proposition \ref{DERIVED3} and the fact that, when the
limit points of $(a_n(\lambda))_{n\geq1}$ lie on $\overline\Gamma_\xi(a)$,
$$
\lims_n|a_n(\lambda)-a| \leq 2|a|\sin\frac{\xi}{2}.
$$

\epr

\br \label{GER}

The class $\lim_n|a_n|=r\in(0,1)$ corresponds to the case $|a_o|=|a_e|\neq0,1$.
Then, $\{a_o,a_e\}=\{a,ae^{i\omega}\}$ with $|a|=r$ and $\omega\in[0,\pi]$,
so, $\sin\frac{\alpha_+}{2}=r\cos\frac{\omega}{2}$ and
$\sin\frac{\alpha_-}{2}=r\sin\frac{\omega}{2}$.
Hence, in this case, the previous Theorem gives
$$
\eta:=\sin\frac{\xi_o}{2}+\sin\frac{\xi_e}{2}<
\cases{\cos\frac{\omega}{2}\cr\sin\frac{\omega}{2}}
\kern6pt\Rightarrow\kern6pt
\{\supp\,\mu\}'\subset\triangle_{\alpha_\pm-\zeta}(\pm\lambda),
\kern7pt \sin\frac{\zeta}{2}=\eta r.
$$

\er

%%%%%%%%%%%%%%%%%%%%%%%%%%%%%%%%%%%%%%%%%%%%%%%%%%%%%%%%%%%%%%%%%%%%%%%%%%%%%%%%%%%%%
\section{Finite truncations of $C(\bsa)$ and the support of the measure}
%%%%%%%%%%%%%%%%%%%%%%%%%%%%%%%%%%%%%%%%%%%%%%%%%%%%%%%%%%%%%%%%%%%%%%%%%%%%%%%%%%%%%

Given an operator on an infinite-dimensional Hilbert space, the search for finite
truncations whose spectra asymptotically approach the spectrum of the full operator is
an old and non trivial problem. For our purposes, the relevant question is, if, given a
sequence $\bsa$ in $\D$, there exist sequences of finite truncations of the unitary
matrix $C(\bsa)$ such that their spectra approximate to the spectrum of $C(\bsa)$, that is,
to the support of the related measure on $\T$ (for the analogous problem concerning Jacobi
matrices and measures on the real line see \cite{BaLoMaTo98, BaLoTo95, IfPa01b}). We will
see that the normal truncations of $C(\bsa)$ on $\ell^2_n$ give a positive answer to this
question.

At this point we have to remember the definitions of $\lims_n$, $\limi_n$ and $\lim_n$
for sequences of subsets of $\C$, in the sense of Hahn \cite{Ha48} or Kuratowski
\cite{Ku61}.

\bd \label{LIMITPOINTS}

Given a sequence $\bsE=(E_n)_{n\geq1}$, $E_n\subset\C$,
$$
\ba{l}
\ds \lims_nE_n := \{\lambda\in\C : \limi_n d(\lambda,E_n)=0\},
\smallskip \\
\ds \limi_nE_n := \{\lambda\in\C : \lim_n d(\lambda,E_n)=0\},
\smallskip \\
\ds E=\lim_nE_n \;{\rm iff}\; \limi_nE_n=E=\lims_nE_n.
\ea
$$
The points in $\lims_nE_n$ are called the (weak) limit points of $\bsE$, while the points
in $\limi_nE_n$ are called the strong limit points of $\bsE$.

\ed

$\limi_nE_n$ and $\lims_nE_n$ are closed sets such that $\limi_nE_n \subset \lims_nE_n$.
The points in $\limi_nE_n$ are the limits of the convergent sequences
$(\lambda_n)_{n\geq1}$ with $\lambda_n \in E_n$, $\forall n\geq1$, while $\lims_nE_n$
contains the strong limit points of all the subsequences of $\bsE$. We have the following
relations.

\bl \label{DIST-LIMITPOINTS}

For any $z\in\C$,
$$
d(z,\lims_nE_n) = \limi_n d(z,E_n) \leq \lims_n d(z,E_n) \leq d(z,\limi_nE_n).
$$

\el

\bpr

If $\lambda \in \limi_nE_n$, there exists $(\lambda_n)_{n\geq1}$, $\lambda_n \in E_n$,
$\forall n\geq1$, such that $\lambda=\lim_n\lambda_n$. Therefore,
$d(z,\lambda) = \lim_n d(z,\lambda_n) \geq \lims_n d(z,E_n)$. Hence,
$d(z,\limi_nE_n) \geq \lims_n d(z,E_n)$.

A similar argument shows that $d(z,\lims_nE_n) \geq \limi_n d(z,E_n)$. So, in the case
$\limi_nd(z,E_n) = \infty$ the relation is proved. Otherwise, let $(E_n)_{n\in\cI}$
be a subsequence of $\bsE$ such that $\limi_nd(z,E_n)=\lim_{n\in\cI}d(z,E_n)$.
Then, $\limi_nd(z,E_n)=\lim_{n\in\cI}d(z,\lambda_n)$, $\lambda_n \in E_n$,
$\forall n\in\cI$. Since $(\lambda_n)_{n\in\cI}$ must be bounded, it has a convergent
subsequence $(\lambda_n)_{n\in\cJ}$. Therefore,
$\lambda=\lim_{n\in\cJ}\lambda_n \in \lims_nE_n$ and
$d(z,\lims_nE_n) \leq d(z,\lambda) = \limi_n d(z,E_n)$.

\epr

Given a sequence $(T_n)_{n\geq1}$ of truncations of an operator $T$, we are interested in
the limit and strong limit points of the related spectra, that is, $\lims_n\spec(T_n)$
and $\limi_n\spec(T_n)$. With this notation, Proposition \ref{GENTRUNC}.1 says that, if
$T$ is a bounded normal band operator on $\ell^2$ and $(T_n)_{n\geq1}$ is a bounded
sequence of normal truncations of $T$, $T_n$ being a truncation on $\ell^2_n$, then
$\spec(T)\subset\limi_n\spec(T_n)$.

Concerning the relation between the limit points of the spectra for different sequences
of truncations, we have the following result.

\bp \label{COMPTRUNC}

For $n\geq1$, let $T_n$, $T'_n$ be bounded truncations of a given operator on the same
subspace. If the truncations $T_n$ are normal, then
$$
\lims_n\spec(T'_n) \subset \{z\in\C : d(z,\lims_n\spec(T_n)) \leq \lims_n\|T'_n-T_n\|\}.
$$

\ep

\bpr

Let $\lambda \in \lims_n\spec(T'_n)$. There exists $(\lambda_n)_{n\in\cI}$,
$\cI\subset\N$, with $\lambda_n\in\spec(T'_n)$, $\forall n\in\cI$, and
$\lambda=\lim_{n\in\cI}\lambda_n$. Since $T_n$ is normal,
$d(\lambda_n,\spec(T_n)) \leq \|T'_n-T_n\|$ and, so,
$d(\lambda,\spec(T_n)) \leq |\lambda-\lambda_n|+\|T'_n-T_n\|$.
Therefore, using Lemma \ref{DIST-LIMITPOINTS} we get
$d(\lambda,\lims_n\spec(T_n)) =
\limi_n d(\lambda,\spec(T_n)) \leq \lims_n\|T'_n-T_n\|$.

\epr

Given a sequence $\bsa=(a_n)_{n\geq1}$ in $\D$ and a sequence $\bsu=(u_n)_{n\geq1}$ in
$\T$, we can consider the corresponding sequence $(C(a_1,\dots,a_{n-1},u_n))_{n\geq1}$ of
finite unitary truncations of $C(\bsa)$. Our aim is to study the relation between the
limit and strong limit points of the spectra of these truncations and the spectrum of
$C(\bsa)$. The spectrum of $C(a_1,\dots,a_{n-1},u_n)$ is the set of zeros of the
para-orthogonal polynomial $p^{u_n}_n$ associated with the measure related to $\bsa$.
This means that, in fact, we are going to study the connection between the support of a
measure on $\T$ and the limit and strong limit points of the zeros of sequences
$(p^{u_n}_n)_{n\geq1}$ of para-orthogonal polynomials associated with this measure. Some
previous results in this direction can be found in \cite{MPOP, Go02}. For convenience, in
what follows we use the notation
$$
\Sigma_n(\bsa;\bsu):=\spec(C(a_1,\dots,a_{n-1},u_n))=\{z\in\C : p_n^{u_n}(z)=0\},
$$
for any sequence $\bsa$ in $\D$ and $\bsu$ in $\T$. The results achieved till now have
the following consequences.

\bt \label{SLP-SPEC-C}

If $\bsa$ is the sequence of Schur parameters associated with a measure $\mu$ on $\T$ and
$\bsu$ is an arbitrary sequence in $\T$,
$$
\supp\,\mu \subset \limi_n\Sigma_n(\bsa;\bsu).
$$
Moreover, for any other sequence $\bsu'$ in $\T$,
$$
\lims_n\Sigma_n(\bsa;\bsu') \subset
\{z\in\T : d(z,\lims_n\Sigma_n(\bsa;\bsu)) \leq \lims_n|u'_n-u_n|\},
$$
and, so,
$$
\lim_n(u'_n-u_n)=0 \quad\Rightarrow\quad
\lims_n\Sigma_n(\bsa;\bsu') = \lims_n\Sigma_n(\bsa;\bsu).
$$

\et

\bpr

The first statement is a direct consequence of Proposition \ref{GENTRUNC}.1.

The second one follows from Proposition \ref{COMPTRUNC}, taking into account that
$\|C(a_1,\dots,a_{n-1},u'_n)-C(a_1,\dots,a_{n-1},u_n)\|=|u'_n-u_n|$, as can be easily
proved using
$C(a_1,\dots,a_{n-1},u_n)=C_o(a_1,\dots,a_{n-1},u_n)C_e(a_1,\dots,a_{n-1},u_n)$.

\epr

The fact that the strong limit points of the zeros of para-orthogonal polynomials include
the support of the orthogonality measure was already proved in \cite[Theorem 8]{Go02},
which deals with sequences of para-orthogonal polynomials with a fixed zero. We have
obtained the result as a particular case of a more general statement of operator theory.
This result will be improved later (see theorems \ref{SLP-C-1}, \ref{SLP-C-2} and
Corollary \ref{SLP-C-21}) although we can not always expect a strict equality between
those strong limit points and the support of the measure, due to the freedom in one of
the zeros for the para-orthogonal polynomials of a given order (see Remark \ref{POP}).

Concerning the weak limit points of the zeros of para-orthogonal polynomials, it was also
shown in \cite[Examples 9 and 10]{Go02} that some of them can lie outside the support of
the measure, even if we fix for all the para-orthogonal polynomials a common zero inside
the support of the measure. However, we can get some information about the location of
these limit points, which will be useful for the study of the convergence of rational
approximants for the Carath\'eodory function of the measure (see Section 5). The next
theorem is an example of this kind of results. If there is a limit point outside the
support of the measure, this theorem establishes how far it can be from the derived set
of this support. The proof, which follows the ideas of \cite[Theorem 2.3]{IfPa01b}
relating to orthogonal truncations of self-adjoint operators, needs the following lemmas.

\bl \label{NORM-LP}

Let $T_0,T\in\frB(H)$ be normal and such that $T-T_0$ is compact. If $T_n:=T[Q_n]$ is a
finite normal truncation of $T$ for $n\geq1$ and $\hat T_n \to T$,
$$
\sup_{\lambda\in\lims_n\spec(T_n)\backslash\spec(T)}\kern-22pt|\lambda-z| \leq
\lims_n\|Q_n\| \kern-3pt \sup_{\lambda\in\spec(T_0)}\kern-3pt|\lambda-z|,
\quad \forall z\in\C.
$$

\el

\bpr

Let $\lambda\in\lims_n\spec(T_n)\backslash\spec(T)$. Then, there exists
$(\lambda_n)_{n\in\cI}$, $\cI\subset\N$, such that $\lim_{n\in\cI}\lambda_n=\lambda$
and $\lambda_n$ is an eigenvalue of the finite truncation $T_n$. Let $x_n$ be a
unitary eigenvector of $T_n$ with eigenvalue $\lambda_n$. Since $T_n$ is normal,
we can suppose that $x_n$ is also an eigenvector of $T^*_n$ with eigenvalue
$\overline\lambda_n$. So, given an arbitrary $y\in H$,
$$
(x_n,(\lambda-T)y) = ((\overline \lambda-\hat T_n^*)x_n,y) + (x_n,(\hat T_n-T)y),
$$
which gives
$$
|(x_n,(\lambda-T)y)| \leq |\lambda-\lambda_n|\|y\| + \|(\hat T_n-T)y\|.
$$
$(x_n)_{n\in\cI}$ is bounded, thus, there exists a subsequence $(x_n)_{n\in\cJ}$ weakly
converging to some $x\in H$. Taking limits in the above inequality for $n\in\cJ$ we get
$((\overline \lambda-T^*)x,y)=0$, $\forall y\in H$, that is, $(\overline \lambda-T^*)x=0$.
Since $T$ is normal and $\lambda\notin\spec(T)$, $\overline\lambda$ is not an eigenvalue of
$T^*$, thus, $x=0$.

Let $T_{0n}:=T_0[Q_n]$. Then, for any $z\in\C$, we can write
$$
(\lambda-z)x_n = (T_{0n}-z)x_n + (T_n-T_{0n})x_n + (\lambda-\lambda_n)x_n,
$$
and, hence,
$$
|\lambda-z| \leq \|Q_n\|\|T_0-z\| + \|Q_n\|\|(T-T_0)x_n\| + |\lambda-\lambda_n|.
$$
The fact that $(x_n)_{n\in\cJ}$ weakly converges to $0$ and $T-T_0$ is compact implies
that $\lim_{n\in\cJ}\|(T-T_0)x_n\|=0$. We can suppose $\lims_n\|Q_n\|<\infty$, otherwise
the inequality of the theorem is trivial. Then, taking limits in the last inequality for
$n\in\cJ$ we obtain
$$
|\lambda-z| \leq \lims_n\|Q_n\|\|T_0-z\|.
$$
$T_0$ is normal, hence, $\|T_0-z\|=\sup_{\lambda\in\spec(T_0)}|\lambda-z|$.
So, the theorem is proved since $\lambda$ was an arbitrary point in
$\lims_n\spec(T_n)\backslash\spec(T)$.

\epr

The previous result for normal operators has the following consequence in the special
case of unitary operators.

\bl \label{UNIT-LP}

Let $U$ be a unitary operator on $H$ and $U_n:=U[Q_n]$ be a finite unitary truncation of
$U$ for $n\geq1$ such that $\hat U_n \to U$. Define $\alpha_0\in[0,\pi]$ by
$$
\cos\frac{\alpha_0}{2} = \limi_n \frac{1}{\|Q_n\|}.
$$
Then, if $\spec_e(U)\subset\triangle_\alpha(w)$,
$$
\alpha>\alpha_0 \quad\Rightarrow\quad
\lims_n\spec(U_n)\backslash\spec(U)\subset\triangle_\beta(w), \quad
\cos\frac{\beta}{2}=\frac{\cos\frac{\alpha}{2}}{\cos\frac{\alpha_0}{2}}.
$$

\el

\bpr

Let us suppose that $\spec_e(U)\subset\triangle_\alpha(w)$. Then, for any
$\epsilon\in(0,\alpha)$, $S_\epsilon:=E_U(\Gamma_\epsilon(w))H$ has finite dimension.
Thus, the unitary operator
$$
U^\epsilon :=
U_{\upharpoonright S_\epsilon^\bot} \oplus (-w 1_{\upharpoonright S_\epsilon})
$$
differs from $U$ in a finite rank perturbation and, so, $U-U^\epsilon$ is compact.
Moreover, $\spec(U^\epsilon)\subset\triangle_\epsilon(w)$. Hence, Lemma \ref{NORM-LP}
gives
$$
\sup_{\lambda\in\lims_n\spec(U_n)\backslash\spec(U)}\kern-22pt|\lambda+w| \leq
\lims_n\|Q_n\| \kern-3pt \sup_{\lambda\in\spec(U^\epsilon)}\kern-3pt|\lambda+w| \leq
2\frac{\cos\frac{\epsilon}{2}}{\cos\frac{\alpha_0}{2}}, \quad \forall
\epsilon\in(0,\alpha).
$$
If $\alpha>\alpha_0$, then
$\ds\frac{\cos\frac{\alpha}{2}}{\cos\frac{\alpha_0}{2}}=\cos\frac{\beta}{2}$,
$\beta\in(0,\pi]$. Thus, from the above inequality,
$$
\sup_{\lambda\in\lims_n\spec(U_n)\backslash\spec(U)}\kern-22pt|\lambda+w| \leq
2\cos\frac{\beta}{2},
$$
which proves the result.

\epr

Now we can get the announced result about the limit points of the zeros of
para-orthogonal polynomials.

\bt \label{LP-C-1}

Let $\bsa$ be the sequence of Schur parameters of a measure $\mu$ on $\T$,
$\{\Gamma_{\alpha_j}(w_j)\}_{j=1}^N$ $(N\in\N\cup\{\infty\})$ being the connected
components of $\T\backslash\{\supp\,\mu\}'$. Let $\bsu$ be a sequence in $\T$ and
define $\alpha_0\in[0,\pi]$ by
$$
\cos\frac{\alpha_0}{2} = \limi_n \frac{1}{\sqrt{1+\frac{|u_n-a_n|^2}{\rho_n^2}}}.
$$
Then, for any $j=1,\dots,N$,
$$
\alpha_j>\alpha_0 \kern6pt\Rightarrow\kern6pt
\lims_n\Sigma_n(\bsa;\bsu)\cap\Gamma_{\beta_j}(w_j) =
\supp\,\mu \cap \Gamma_{\beta_j}(w_j), \kern6pt
\cos\frac{\beta_j}{2}=\frac{\cos\frac{\alpha_j}{2}}{\cos\frac{\alpha_0}{2}}.
$$

\et

\bpr

Since $\supp\,\mu \subset \limi_n\Sigma_n(\bsa;\bsu) \subset \lims_n\Sigma_n(\bsa;\bsu)$,
we just have to prove that
$\lims_n\Sigma_n(\bsa;\bsu)\backslash\supp\,\mu\subset\triangle_{\beta_j}(w_j)$ for each
$j$ such that $\alpha_j>\alpha_0$. Remember that $C(a_1,\dots,a_{n-1},u)$ is the
unitary truncation of $C(\bsa)$ on $\ell^2_n$ associated with the projection
$Q_n(\bsa;u)$ whose norm is given in Corollary \ref{NORMTRUNC-C}. From Proposition
\ref{GENTRUNC}.1, $\hat C(a_1,\dots,a_{n-1},u) \to C(\bsa)$. So, if $\alpha_j>\alpha_0$,
a direct application of Lemma \ref{UNIT-LP} with $U=C(\bsa)$ and $Q_n=Q_n(\bsa;u_n)$
gives $\lims_n\Sigma_n(\bsa;\bsu)\backslash\supp\,\mu\subset\triangle_{\beta_j}(w_j)$.

\epr

The previous theorem says that $\lims_n\Sigma_n(\bsa;\bsu)$ can differ from $\supp\,\mu$
only in $\Gamma_{\alpha_j}(w_j)\backslash\Gamma_{\beta_j}(w_j)$ if $\alpha_j>\alpha_0$, or
in $\Gamma_{\alpha_j}(w_j)$ if $\alpha_j\leq\alpha_0$.

If $a_n\neq0$ for any $n$ big enough, the best choice for the sequence $\bsu$ in Theorem
\ref{LP-C-1} is $u_n=\frac{a_n}{|a_n|}$. Then,
$$
\cos\alpha_0=\limi_n|a_n|,
$$
so, $\alpha_0\leq\frac{\pi}{2}$. Taking into account Theorem \ref{SLP-SPEC-C}, Theorem
\ref{LP-C-1} also works with the above value of $\alpha_0$ if
$\lim_n(u_n-\frac{a_n}{|a_n|})=0$. In particular, if $\lim_n|a_n|=1$ we can choose $\bsu$
such that $\lim_n(u_n-a_n)=0$, which gives $\alpha_0=0$ and, hence,
$\supp\,\mu=\limi_n\Sigma_n(\bsa;\bsu)=\lims_n\Sigma_n(\bsa;\bsu)$. So we get the following
consequence of Theorem \ref{LP-C-1}.

\bc \label{LP-C-11}

If $\bsa$ is the sequence of Schur parameters associated with a measure $\mu$ on $\T$
and $\bsu$ is a sequence in $\T$,
$$
\lim_n(u_n-a_n)=0 \quad\Rightarrow\quad \lim_n\Sigma_n(\bsa;\bsu) = \supp\,\mu.
$$

\ec

As we pointed out, $\{\supp\,\mu\}'=-\frak L\{a_{n+1}/a_n\}$ when $\lim_n|a_n|=1$, so, if
$\lim_n(u_n-a_n)=0$, $\lim_n\Sigma_n(\bsa;\bsu)$ coincides with $-\frak L\{a_{n+1}/a_n\}$
plus, at most, a countable set that can accumulate only on $-\frak L\{a_{n+1}/a_n\}$.

\bex \label{EX11}
Rotated asymptotically 2-periodic Schur parameters.

\smallskip

\noindent Let us suppose that $\lim_na_{2n-1}(\lambda)=a_o$ and $\lim_na_{2n}(\lambda)=a_e$
for some $\lambda\in\T$. We know that
$\T\backslash\{\supp\,\mu\}'=\Gamma_{\alpha_+}(\lambda)\cup\Gamma_{\alpha_-}(-\lambda)$
where $\alpha_\pm$ is given in (\ref{PER}). If $u_n=\frac{a_n}{|a_n|}$, then
$\cos\alpha_0=\min\{|a_o|,|a_e|\}$ and
$$
\alpha_\pm>\alpha_0 \kern3pt\Leftrightarrow\kern3pt
\rho_o\rho_e\mp\re(\overline a_oa_e)<\min\{|a_o|,|a_e|\}.
$$
Under these conditions,
$\lims_n\Sigma_n(\bsa;\bsu)\cap\Gamma_{\beta_\pm}(\pm\lambda) =
\supp\,\mu \cap \Gamma_{\beta_\pm}(\pm\lambda)$
where
$$
\cos\frac{\beta_\pm}{2}=
\sqrt{\frac{1+\rho_o\rho_e\mp\re(\overline a_oa_e)}{1+\min\{|a_o|,|a_e|\}}}.
$$

\eex

\bex \label{EX12}
$\frak L\{\bsa(\lambda)\}=\{a,b\}$, $\lambda\in\T$, $a\neq b$, $|a|=|b|\neq0$.

\smallskip

\noindent Following Example \ref{EXDER2}, we can suppose $b=ae^{i\zeta}$,
$\zeta\in(0,\pi]$, and, then
$$
\ba{l}
\ds \Gamma_{\alpha_1}(\lambda)\subset\T\backslash\{\supp\,\mu\}',
\quad \alpha_1=\alpha-\zeta, \quad
{\rm if} \kern4pt \sin\frac{\zeta}{2}<|a|,
\smallskip \\
\ds
\Gamma_{\alpha_2}(-\lambda e^{\pm i\frac{\zeta}{2}})
\subset\T\backslash\{\supp\,\mu\}',
\quad \alpha_2=\alpha+\frac{\zeta}{2}-\pi, \quad
{\rm if} \kern4pt \cos\frac{\zeta}{4}<|a|,
\ea
$$
where $\alpha\in(0,\pi)$ is given by $\sin\frac{\alpha}{2}=|a|$.
Let $u_n=\frac{a_n}{|a_n|}$, so that $\cos\alpha_0=|a|$, that is,
$\alpha_0=\alpha-\zeta_0$ with
$$
\sin\frac{\zeta_0}{2}=(2|a|-1)\sqrt{\frac{1+|a|}{2}},
\quad -\frac{\pi}{2}<\zeta_0<\pi.
$$
Hence, we get
$$
\ba{l}
\ds \alpha_1>\alpha_0 \kern3pt\Leftrightarrow\kern3pt \zeta<\zeta_0
\kern3pt\Leftrightarrow\kern3pt \sin\frac{\zeta}{2}<(2|a|-1)\sqrt{\frac{1+|a|}{2}},
\\
\ds \alpha_2>\alpha_0 \kern3pt\Leftrightarrow\kern3pt \pi-\frac{\zeta}{2}<\zeta_0
\kern3pt\Leftrightarrow\kern3pt \cos\frac{\zeta}{4}<(2|a|-1)\sqrt{\frac{1+|a|}{2}}.
\ea
$$
Under each of these conditions,
$\lims_n\Sigma_n(\bsa;\bsu)\cap\Gamma_{\beta_1}(\lambda) =
\supp\,\mu \cap \Gamma_{\beta_1}(\lambda)$
and
$\lims_n\Sigma_n(\bsa;\bsu)\cap\Gamma_{\beta_2}(-\lambda e^{\pm i\frac{\zeta}{2}}) =
\supp\,\mu \cap \Gamma_{\beta_2}(-\lambda e^{\pm i\frac{\zeta}{2}})$
respectively, where
$$
\cos\frac{\beta_1}{2}=
\frac{|a|\sin\frac{\zeta}{2}+\rho\cos\frac{\zeta}{2}}{\sqrt{\frac{1+|a|}{2}}}, \qquad
\cos\frac{\beta_2}{2}=
\frac{|a|\cos\frac{\zeta}{4}+\rho\sin\frac{\zeta}{4}}{\sqrt{\frac{1+|a|}{2}}},
$$
being $\rho=\sqrt{1-|a|^2}$.

\eex

\medskip

We can go further in the analysis of $\lims_n\Sigma_n(\bsa;\bsu)$ and
$\limi_n\Sigma_n(\bsa;\bsu)$ using the analytic properties of the para-orthogonal
polynomials. The following result for the corresponding zeros was proved in \cite{MPOP}
and \cite{Go02}.

\medskip

\noindent{\bf Theorem A.} \cite[Corollary 2]{MPOP} \cite[Theorem 2]{Go02} {\it Given a
measure $\mu$ on $\T$, the closure $\overline\Gamma$ of any arc
$\Gamma\subset\T\backslash\supp\,\mu$ contains at most one zero of the para-orthogonal
polynomial $p_n^u$ related to $\mu$ for any $u\in\T$ and $n\in\N$.}

\medskip

With this property and Theorem \ref{SLP-SPEC-C} we can achieve the following result.

\bt \label{LP-C-2}

Let $\mu$ be a measure on $\T$ with a sequence $\bsa$ of Schur parameters, and let $\bsu$
be a sequence in $\T$. Consider a connected component $\Gamma$ of
$\T\backslash\{\supp\,\mu\}'$ and $w\in\Gamma$.
\be
\item
$\ds w \in \Sigma_n(\bsa;\bsu) \;\; \forall n \geq 1 \;\Rightarrow\;
\lims_n\Sigma_n(\bsa;\bsu)\cap\Gamma=(\supp\,\mu\cap\Gamma)\cup\{w\}$.
\item
$\ds w \in \limi_n\Sigma_n(\bsa;\bsu)\backslash\supp\,\mu \;\Rightarrow\;
\lims_n\Sigma_n(\bsa;\bsu)\cap\Gamma=(\supp\,\mu\cap\Gamma)\cup\{w\}$.
\ee

\et

\bpr

$\Gamma\cap\supp\,\mu$ is at most a countable set which can accumulate only at
the endpoints of $\overline\Gamma$. Consider one of the two connected components of
$\Gamma\backslash\{w\}$, let us say $\Gamma_+=(w,w_+)$. Let
$\Gamma_+\backslash\supp\,\mu=(w,w_1)\cup(w_1,w_2)\cup\cdots$ be the decomposition in
connected components. From Theorem A, it is clear that $(w_j,w_{j+1})$ has, at most, one
point in $\Sigma_n(\bsa;\bsu)$ for each $j,n\geq1$.

Assume that $w \in \Sigma_n(\bsa;\bsu)$, $\forall n\geq1$. Then, Theorem A implies that
$(w,w_1) \cap \Sigma_n(\bsa;\bsu) = \emptyset$, $\forall n\geq1$. Hence, for $j\geq1$,
since $w_j\in\supp\,\mu \subset \limi_n\Sigma_n(\bsa;\bsu)$, it is necessary that
$(w_j,w_{j+1}) \cap \Sigma_n(\bsa;\bsu) = \{\lambda_j^{(n)}\}$ for $n$ greater than
certain index $n_j$, and $\lim_n\lambda_j^{(n)}=w_j$. So, we conclude that
$\Gamma_+ \cap \lims_n\Sigma_n(\bsa;\bsu) = \{w_1,w_2,\dots\} = \Gamma_+ \cap \supp\,\mu$.
A similar analysis for the other connected component of $\Gamma\backslash\{w\}$
finally gives
$\Gamma \cap \lims_n\Sigma_n(\bsa;\bsu) = (\Gamma \cap \supp\,\mu) \cup \{w\}$.

Let us suppose now that $w \in \limi_n\Sigma_n(\bsa;\bsu)\backslash\supp\,\mu$. Let
$\Gamma_-=(w_-,w)$ be the other connected component of $\Gamma\backslash\{w\}$,
$\Gamma_-\backslash\supp\,\mu = (w_{-1},w_0)\cup(w_{-2},w_{-1})\cup\cdots$ being the
decomposition in connected components. From Theorem A and the condition for $w$ we
conclude that $(w_{-1},w_1) \cap \Sigma_n(\bsa;\bsu) = \{\lambda_0^{(n)}\}$ for $n$
greater than certain index $n_0$, and $\lim_n\lambda_0^{(n)}=w$. From here on, similar
arguments to the previous case prove that
$\Gamma \cap \lims_n\Sigma_n(\bsa;\bsu) = (\Gamma \cap \supp\,\mu) \cup \{w\}$.

\epr

Theorem \ref{LP-C-2} says that, if we choose a sequence of para-orthogonal polynomials
with a fixed zero $w$ outside $\{\supp\,\mu\}'$, or with a zero converging to a point $w$
outside $\supp\,\mu$, then, in the connected component of $\T\backslash\{\supp\,\mu\}'$
where $w$ lies, the limit points of the zeros of the para-orthogonal polynomials coincide
with $\supp\,\mu$ up to, at most, the point $w$.

\br \label{FIX}

The above result can be read as a statement about the zeros of the para-orthogonal
polynomials $p_n^{u_n}(z)=z\varphi_{n-1}(z)+u_n\varphi_{n-1}^*(z)$. The requirement for
fixing a common zero $w$ for all these polynomials means that we have to choose
$u_n=-w\varphi_{n-1}(w)/\varphi_{n-1}^*(w)$ for $n\geq1$. However, if we
consider Theorem \ref{LP-C-2} as a result about the eigenvalues of the unitary matrices
$C(a_1,\dots,a_{n-1},u_n)$, the interesting question is: how to get from the sequence
$\bsa$ of Schur parameters the sequence $\bsu=\bsu^w$ that fixes a common zero $w$ for
all the polynomials $p_n^{u_n}$? Since the polynomial $p_n^{u_n}$ is proportional to
$q_n^{v_n}=\varphi_n^*-\overline v_n\varphi_n$ with
$v_n=-u_n\frac{1-a_n\overline u_n}{1-\overline a_nu_n}$ (see Remark \ref{POP}),
we find that such a sequence $\bsu=\bsu^w$ must satisfy $u_{n+1}\overline v_n=-w$
and, thus,
\beq \label{FIXeq}
\ba{l}
u_1^w=-w,
\\
\ds u_{n+1}^w=wu_n^w\frac{1-a_n\overline u_n^w}{1-\overline a_nu_n^w},
\quad \forall n\geq1.
\ea
\eeq
This recurrence answers the question.

\er

\medskip

Corollary \ref{LP-C-11} stated that, for the family of measures whose Schur parameters
tend to the unit circle, it is possible to choose a sequence of para-orthogonal
polynomials whose zeros exactly converge to the support of the measure. Theorem
\ref{LP-C-2} gives another class of measures where this essentially happens, as the
following corollary shows.

\bc \label{LP-C-21}

Let $\mu$ be a measure on $\T$ with a sequence $\bsa$ of Schur parameters, and let $\bsu$
be a sequence in $\T$. Assume that $\{\supp\,\mu\}'$ is connected.
\be
\item
$\ds w \in \Sigma_n(\bsa;\bsu)\backslash\{\supp\,\mu\}' \;\; \forall n \geq 1 \;\Rightarrow\;
\lim_n\Sigma_n(\bsa;\bsu) = \supp\,\mu\cup\{w\}$.
\item
$\ds w \in \limi_n\Sigma_n(\bsa;\bsu)\backslash\supp\,\mu \;\Rightarrow\;
\lim_n\Sigma_n(\bsa;\bsu) = \supp\,\mu\cup\{w\}$.
\ee

\ec

The results of Section 3 provide very general situations where Theorem \ref{LP-C-2} can
be applied. Concerning the more stringent result of Corollary \ref{LP-C-21}, the
following example gives a remarkable situation where it holds.

\bex \label{EX2} The L\'opez class.

\smallskip

\noindent If $\lim_n\frac{a_{n+1}}{a_n}=\lambda\in\T$ and $\lim_n|a_n|=r\in(0,1)$, we
know that $\{\supp\,\mu\}'=\triangle_\alpha(\lambda)$ with $\alpha\in(0,\pi)$ given by
$\sin\frac{\alpha}{2}=r$. Therefore, for any $w\in\Gamma_\alpha(\lambda)$,
$\lim_n\Sigma_n(\bsa;\bsu^w) = \supp\,\mu\cup\{w\}$.

\eex

\smallskip

Theorem \ref{LP-C-2}.2 has a consequence about the strong limit points of the zeros of
para-orthogonal polynomials. It implies that $\limi_n\Sigma_n(\bsa;\bsu)$ can differ from
$\supp\,\mu$ in, at most, one point in each connected component of
$\T\backslash\{\supp\,\mu\}'$. This result will be improved in theorems \ref{SLP-C-1},
\ref{SLP-C-2} and Corollary \ref{SLP-C-21}, which use some results of \cite{Kh02}. In
this work, S. Khrushchev defines the so-called class of Markoff measures on the unit
circle (Mar$(\T)$), which includes, as a particular case, all the measures whose support
does not cover the unit circle \cite[pg. 268]{Kh02}. For a measure $\mu$ in this class he
proves some results that give information about the asymptotics of
$\varphi_n/\varphi_n^*$ in $\T\backslash\supp\,\mu$, $(\varphi_n)_{n\geq0}$ being the
orthonormal polynomials in $L^2_\mu$. As we will see, this is a key tool to control the
strong limit points of the zeros of the para-orthogonal polynomials. Let us summarize the
referred results in \cite{Kh02}.

\medskip

\noindent{\bf Theorem B.} \cite[Lemma 8.4.1]{Kh02} {\it Let $\mu\in{\rm Mar}(\T)$. Then
there exists a positive number $\delta(\mu)$ such that
$$
\sup_{|z|\leq1/2} \left|\frac{\varphi_n(z)}{\varphi_n^*(z)}\right| > \delta(\mu) > 0,
\quad \forall n\geq1.
$$
}

\medskip

\noindent{\bf Theorem C.} \cite[Corollary 8.6]{Kh02}
{\it Let $B$ be a Blaschke product with zeros $\{z_j\}$ such that $|B(z_0)|>\delta>0$ for
some $z_0$, $|z_0|\leq1/2$. If $|z-z_j|\geq\epsilon>0$, then
$|B(z)|>c=c(\delta,\epsilon)>0$.}

\medskip

The interest of this last result is that $\varphi_n/\varphi_n^*$ is always a Blaschke
product. With the above tools we get a result for the asymptotic behaviour of the
sequence $(\varphi_n^*/\varphi_n)_{n\geq0}$ on $\C\backslash\co(\supp\,\mu)$, where
$\co(E)$ means the convex hull of $E\subset\C$.

\bp \label{BLAS}

Let $\mu$ be a measure on $\T$ and let $g_n:=\varphi_n^*/\varphi_n$ for $n\geq0$,
$(\varphi_n)_{n\geq0}$ being the orthonormal polynomials in $L^2_\mu$. Then, any
subsequence of $(g_n)_{n\geq0}$ has a subsequence which uniformly converges on
compact subsets of $\C\backslash\co(\supp\,\mu)$.

\ep

\bpr

Let us suppose first that $\supp\,\mu\neq\T$, and let $\cK$ be a compact subset of
$\C\backslash\co(\supp\,\mu)$. Since $\co(\supp\,\mu)$ is also compact, the distance
between $\cK$ and $\co(\supp\,\mu)$ must be a positive number $\epsilon$. The zeros
$\{z_j^{(n)}\}_{j=1}^n$ of the polynomial $\varphi_n$ lie on $\co(\supp\,\mu)$, so,
$|z-z_j^{(n)}|\geq\epsilon>0$, $j=1,\dots,n$, for all $z\in\cK$ and $n\in\N$. Therefore,
Theorem B and C imply that $|1/g_n(z)|>c(\delta(\mu),\epsilon)>0$, $\forall z \in \cK$,
$\forall n\in\N$, and, thus, $(g_n)_{n\geq0}$ is uniformly bounded on $\cK$. That is,
$(g_n)_{n\geq0}$ is uniformly bounded on any compact subset of
$\C\backslash\co(\supp\,\mu)$. This is also true for the case $\supp\,\mu=\T$ since
$|g_n(z)|\leq1$ for $|z|\geq1$. Therefore, $(g_n)_{n\geq0}$ is always a normal family in
$\C\backslash\co(\supp\,\mu)$, which proves the proposition.

\epr

Now we are ready to prove the main results about the strong limit points of the zeros of
para-orthogonal polynomials. The following set will be important in the next discussions.

\bd \label{LL}

Given a sequence $\bsE=(E_n)_{n\geq1}$, $E_n\subset\C$, we define $\limss_nE_n$ as the
set of points $\lambda\in\C$ such that, for some infinite set $\cI\subset\N$,
$$
\lim_{n\in\cI}d(\lambda,E_n)=\lim_{n\in\cI}d(\lambda,E_{n+1})=0.
$$
We call $\limss_nE_n$ the set of double limit points of the sequence $\bsE$.

\ed

Obviously, $\limss_nE_n$ is a closed set such that
$\limi_nE_n\subset\limss_nE_n\subset\lims_nE_n$.

\bt \label{SLP-C-1}

Let $\bsa$ be the sequence of Schur parameters of a measure $\mu$ on $\T$ and let $\bsu$
be a sequence in $\T$. Then, $\limi_n\Sigma_n(\bsa;\bsu)$ coincides with $\supp\,\mu$
except, at most, at one point. If this point exists,
$\limss_n\Sigma_n(\bsa;\bsu)=\limi_n\Sigma_n(\bsa;\bsu)$ and, thus,
$\limss_n\Sigma_n(\bsa;\bsu)$ equals $\supp\,\mu$ up to such a point.

\et

\bpr

It is enough to prove that the conditions
$w\in\limi_n\Sigma_n(\bsa;\bsu)\backslash\supp\,\mu$ and
$z\in\limss_n\Sigma_n(\bsa;\bsu)\backslash\supp\,\mu$ imply $z=w$. Let
$w\in\limi_n\Sigma_n(\bsa;\bsu)$ and $z\in\limss_n\Sigma_n(\bsa;\bsu)$. There exist two
sequences $(w_n)_{n\geq1}$ and $(z_n)_{n\geq1}$ with $w_n,z_n\in\Sigma_n(\bsa;\bsu)$,
$\forall n\geq1$, such that $\lim_nw_n=w$ and $\lim_{n\in\cI}z_n=\lim_{n\in\cI}z_{n+1}=z$
for some infinite set $\cI\subset\N$. Let $\mu$ be the measure whose sequence of Schur
parameters is $\bsa$, and let $g_n=\varphi_n^*/\varphi_n$ for $n\geq0$, $(\varphi_n)_{n\geq0}$
being the orthonormal polynomials in $L^2_\mu$. Since $w_n$ and $z_n$ are zeros of the same
para-orthogonal polynomial $p_n^{u_n}$, we get
\beq \label{g1}
\overline z_{n+1} g_n(z_{n+1}) = \overline w_{n+1} g_n(w_{n+1}), \quad n\geq0.
\eeq
Taking into account that $p_n^{u_n}$ is proportional to
$q_n^{v_n}=\varphi_n^*-\overline v_n\varphi_n$ for some $v_n\in\T$ (see Remark \ref{POP}),
we also find that
\beq \label{g2}
g_n(z_n) = g_n(w_n), \quad n\geq1.
\eeq

Assume that $w,z\notin\supp\,\mu$ and let $\cK$ be a compact subset of
$\T\backslash\supp\,\mu$ containing two open arcs centred at $w$ and $z$ respectively.
$w_n,z_n\in\cK$ for any $n$ big enough. On the other hand, Proposition \ref{BLAS} ensures the
uniform convergence on $\cK$ of a subsequence $(g_n)_{n\in\cJ}$, $\cJ\subset\cI$. If $g$
is the uniform limit of this subsequence,
$\lim_{n\in\cJ}g_n(w_n)=\lim_{n\in\cJ}g_n(w_{n+1})=g(w)$ and
$\lim_{n\in\cJ}g_n(z_n)=\lim_{n\in\cJ}g_n(z_{n+1})=g(z)$.
Taking limits for $n\in\cJ$ in (\ref{g1}) and (\ref{g2}) we conclude that $z=w$ since
$g(z)\neq0$ (in fact, $|g(z)|=1$).

\epr

It is clear that we can control the possible strong limit point of the zeros that lies
outside the support of the measure choosing a sequence of para-orthogonal polynomials
with a fixed zero outside this support. More surprising is that the choice of a fixed
zero in the support always gives an exact equality between the strong limit points and
the support of the measure. This is a consequence of the next theorem, which delimits the
possible double limit points outside the support of the measure. In fact, it provides for
any measure sequences of para-orthogonal polynomials that ensure the strict equality
between the double limit points of the zeros and the support of the measure.

\bt \label{SLP-C-2}

If $\bsa$ is the sequence of Schur parameters of a measure $\mu$ on $\T$ and $\bsu$ is a
sequence in $\T$,
$$
\limss_n \Sigma_n(\bsa;\bsu)\backslash\supp\,\mu \subset
\frak L \left\{ \frac{u_{n+1}}{u_n}\frac{1-\overline a_nu_n}{1-a_n\overline u_n} \right\}.
$$

\et

\bpr

Let $z\in\limss_n\Sigma_n(\bsa;\bsu)$. There exists a sequence $(z_n)_{n\geq1}$ such that
$z_n\in\Sigma_n(\bsa;\bsu)$, $\forall n\geq1$, and
$\lim_{n\in\cI}z_n=\lim_{n\in\cI}z_{n+1}=z$ for some infinite set $\cI\subset\N$. Since
$z_n$ is a zero of $p_n^{u_n}$ and $q_n^{v_n}$ with $v_n=-u_n\frac{1-a_n\overline
u_n}{1-\overline a_nu_n}$ (see Remark \ref{POP}), we find that
$u_{n+1}=-z_{n+1}\overline{g_n(z_{n+1})}$ and $v_n=\overline{g_n(z_n)}$ for $n\geq1$,
which gives \beq \label{uu} \frac{u_{n+1}}{u_n}\frac{1-\overline a_nu_n}{1-a_n\overline
u_n} = z_{n+1}\frac{g_n(z_n)}{g_n(z_{n+1})}, \quad n\geq1. \eeq

Let us suppose that $z\notin\supp\,\mu$. Using again Proposition \ref{BLAS} we find that
a subsequence $(g_n)_{n\in\cJ}$, $\cJ\subset\cI$, uniformly converges to a function $g$
on a compact subset of $\T\backslash\supp\,\mu$ containing an open arc centred at $z$.
Hence, $\lim_{n\in\cJ}g_n(z_n)=\lim_{n\in\cJ}g_n(z_{n+1})=g(z)$ and, taking limits for
$n\in\cJ$ in (\ref{uu}), it follows that
$z \in \frak L \left\{
\frac{u_{n+1}}{u_n}\frac{1-\overline a_nu_n}{1-a_n\overline u_n} \right\}$.

\epr

As a first consequence of the previous theorem we find infinitely many sequences of
para-orthogonal polynomials $(p_n^{u_n})_{n\geq1}$, the double limit points of whose
zeros coincide exactly with $\supp\,\mu$. They are those defined by sequences $\bsu$ such
that
$\frak L
\left\{ \frac{u_{n+1}}{u_n}\frac{1-\overline a_nu_n}{1-a_n\overline u_n} \right\}
\subset \supp\,\mu$.

An interesting choice for $\bsu$ is given by the phases of $\bsa$, that is,
$u_n=\frac{a_n}{|a_n|}$ if $a_n\neq0$ and $u_n$ arbitrarily chosen in $\T$ otherwise.
Then, the previous theorem states that
$\limss_n\Sigma_n(\bsa;\bsu)\backslash\supp\,\mu\subset\frak L \{\overline u_nu_{n+1}\}$.

If we are interested in locating the possible strong limit point outside $\supp\,\mu$ at
a certain place $w\in\T$, we can choose $\bsu$ so that
$\lim_n \frac{u_{n+1}}{u_n}\frac{1-\overline a_nu_n}{1-a_n\overline u_n} = w$.
Then, $\limss_n\Sigma_n(\bsa;\bsu)\subset\supp\,\mu\cup\{w\}$. So, if $w\in\supp\,\mu$,
$\limss_n\Sigma_n(\bsa;\bsu)=\supp\,\mu$.

A particular choice of $\bsu$ which ensures
$\limss_n\Sigma_n(\bsa;\bsu)\subset\supp\,\mu\cup\{w\}$
is given by the recurrence
\beq \label{uw}
u_{n+1}=wu_n\frac{1-a_n\overline u_n}{1-\overline a_nu_n},
\quad n\geq1,
\eeq
with an initial condition $u_1=u$, $u$ an arbitrary point in $\T$. Let us study the
form of the related para-orthogonal polynomials $p_n^{u_n}$.

Without loss of generality we write $u_n=-wr_{n-1}(w)/r_{n-1}^*(w)$, with $r_n$ a
polynomial of degree $n$. (\ref{uw}) is equivalent to
$r_n(w)/r_n^*(w)=s_n(w)/s_n^*(w)$, $s_n(w)=wr_{n-1}(w)+a_nr_{n-1}^*(w)$. So,
$r_n(w)=\lambda_n(wr_{n-1}(w)+a_nr_{n-1}^*(w))$, $\lambda_n\in\R\backslash\{0\}$.
This equation has two independent solutions: $(\delta_n\varphi_n(w))_{n\geq0}$,
$(i\delta_n\psi_n(w))_{n\geq0}$, where
$\delta_n=\lambda_1\rho_1\cdots\lambda_n\rho_n$, $(\varphi_n)_{n\geq0}$ are the
orthonormal polynomials related to the Schur parameters $\bsa$, and
$(\psi_n)_{n\geq0}$ are the orthonormal second kind polynomials, associated with the
Schur parameters $-\bsa$ \cite{Sz75,Ge61}.

Therefore, the sequence $\bsu$ satisfies (\ref{uw}) if and only if $p_n^{u_n}(z)$
is proportional to
$p_{n-1}^*(w) z\varphi_{n-1}(z) - wp_{n-1}(w) \varphi_{n-1}^*(z)$, where
$p_n=c_1\varphi_n+ic_2\psi_n$, $(c_1,c_2)\in\R^2\backslash\{(0,0)\}$.
Since $\rho_np_n(w)=wp_{n-1}(w)+a_np_{n-1}^*(w)$, this is equivalent to saying that
$p_n^{u_n}(z)$ is proportional to $p_n^*(w) \varphi_n(z) - p_n(w) \varphi_n^*(z)$.
Notice that the initial condition $u=-w$ means that $p_0$ is real, which corresponds
to the case $c_2=0$, giving the sequence of para-orthogonal polynomials with a fixed
zero at $w$ studied in Remark \ref{FIX}.

Summarizing, as a consequence of Theorem \ref{SLP-C-2} we obtain the following result.

\bc \label{SLP-C-21}

Let $\mu$ be a measure on $\T$, $(\varphi_n)_{n\geq0}$ be the orthonormal polynomials in
$L^2_\mu$ and $(\psi_n)_{n\geq0}$ be the related orthonormal second kind polynomials.
Any sequence $(P_n)_{n\geq1}$ of para-orthogonal polynomials given by
$$
\ba{c}
P_n(z) := p_n^*(w) \varphi_n(z) - p_n(w) \varphi_n^*(z),
\medskip \\
p_n := c_1\varphi_n + ic_2\psi_n,
\quad (c_1,c_2)\in\R^2\backslash\{(0,0)\}, \quad w\in\T,
\ea
$$
has the property that the double limit points of the corresponding zeros coincide with
$\supp\,\mu$ except, at most, at the point $w$. Hence, if $w\in\supp\,\mu$, the double
limit points of the zeros exactly coincide with $\supp\,\mu$. When $c_2=0$, all the
para-orthogonal polynomials have a common zero at $w$ and the double and strong limit
points coincide with $\supp\,\mu\cup\{w\}$.

\ec

The para-orthogonal polynomials $(P_n)_{n\geq1}$ given in the previous corollary
appeared previously in \cite{MPOP}, where it was proved that they have other interesting
properties concerning the interlacing of zeros: for all $n$, $P_n$ and $P_{n+1}$ have
interlacing zeros in $\T\backslash\{w\}$ \cite[Theorem 1]{MPOP}.

In \cite{Go02}, L. Golinskii conjectured that the strong limit points (which he called the
strong attracting points) of the zeros of para-orthogonal polynomials with a fixed zero
lying on the support of the measure, must coincide with this support. Theorems
\ref{SLP-C-1}, \ref{SLP-C-2} and Corollary \ref{SLP-C-21}, not only confirm this
conjecture, but go even further in two senses: the achieved results cover any sequence of
para-orthogonal polynomials, not only the case of a fixed zero on the support of the
measure; even in this case, the results are stronger than the one conjectured by L.
Golinskii since we have proved the equality between the support of the measure and the
double limit points of the zeros of the para-orthogonal polynomials.

The previous results give a method for approximating the support of a measure $\mu$
on the unit circle starting from its sequence $\bsa$ of Schur parameters, based
on the computation of the eigenvalues of the finite unitary matrices
$C(a_1,\dots,a_{n-1},u_n)$ for a sequence $\bsu$ in $\T$. A recommendable choice is the
sequence $\bsu=\bsu^w$ given in (\ref{FIXeq}) that fixes a common eigenvalue $w$ for all
the finite matrices, because it permits us to control the only possible strong limit point
that, according to Theorem \ref{SLP-C-1}, can lie outside $\supp\,\mu$. In this case
Theorem \ref{SLP-C-2} proves that the double limit points coincide with
$\supp\,\mu\cup\{w\}$, so the computation of the eigenvalues for pairs of consecutive
matrices can be used to eliminate those weak limit points that are spurious points of
$\supp\,\mu$.

The following figures show some examples of the previous method of approximation. They
represent $\Sigma_n(\bsa;\bsu)$ for some choices of $\bsa$, $\bsu$ and
$n=50,51,200,201,1000,1001$. The computations have been made applying the double
precision routines of MATLAB to the calculation of the eigenvalues of
$C(a_1,\dots,a_{n-1},u_n)$. We have to remark that the computations can be also made
using a Hessenberg matrix unitarily equivalent to $C(a_1,\dots,a_{n-1},u_n)$
\cite{AmGrRe91,MIN}. However, although the time of computation of eigenvalues is only a
little smaller with the five-diagonal representation (using the standard routines), the
computational cost of building the matrix is much bigger in the Hessenberg case, growing
very much faster as $n$ increases.

The first three figures correspond to different choices of $\bsu$ in the case of constant
Schur parameters $a_n=\frac{1}{2}$, where $\supp\,\mu=\triangle_{\frac{\pi}{3}}(1)$.
Corollary \ref{LP-C-21} proves that $\lim_n\Sigma_n(\bsa;\bsu^1)=\supp\,\mu\cup\{1\}$, as
can be seen in Figure 1. According to Corollary \ref{SLP-C-21},
$\limss_n\Sigma_n(\bsa;\bsu^{-1})=\supp\,\mu$. This is in agreement with Figure 2, where
we see that $\lim_n\Sigma_{2n-1}(\bsa;\bsu^{-1})=\supp\,\mu$ but
$\lim_n\Sigma_{2n}(\bsa;\bsu^{-1})=\supp\,\mu\cup\{1\}$, so, $1$ is a weak but not a
double limit point. Such behaviour was predicted by L. Golinskii in \cite[Example
10]{Go02}. As we have seen throughout the paper, another interesting choice is
$u_n=\frac{a_n}{|a_n|}$ which, used in Figure 3, seems to give
$\lim_n\Sigma_n(\bsa;\bsu)=\supp\,\mu$, although Theorem \ref{SLP-C-2} says that
$\limss_n\Sigma_n(\bsa;\bsu)$ could differ from $\supp\,\mu$ in (at most) the point 1.

The next three figures deal with the 2-periodic Schur parameters $a_{2n-1}=\frac{1}{4}$,
$a_{2n}=\frac{3}{4}$, whose measure satisfies
$\{\supp\,\mu\}'=\triangle_{\alpha_+}(1)\cap\triangle_{\alpha_-}(-1)$ with
$\alpha_+\approx0.35\pi$ and $\alpha_-\approx0.19\pi$. According to those figures there
are no isolated mass points. From Corollary \ref{SLP-C-21},
$\limss_n\Sigma_n(\bsa;\bsu^1)=\supp\,\mu\cup\{1\}$, while Theorem \ref{LP-C-2} states
that, apart from the point $1$, there cannot be weak limit points in the gap around $1$.
These results agree with Figure 4 which shows that, in this case, $-1$ is the only weak
limit point that is not a double limit point. The choice $\bsu=\bsu^i$ fixes a common
eigenvalue at $\supp\,\mu$ and, thus, Corollary \ref{SLP-C-21} implies that
$\limss_n\Sigma_n(\bsa;\bsu^i)=\supp\,\mu$. However, we observe in Figure 5 that this
choice yields a more chaotic behaviour of the eigenvalues in the gaps, which could give weak
limit points outside $\supp\,\mu$. As in the case $a_n=\frac{1}{2}$, Figure 6 seems to
indicate that $\limss_n\Sigma_n(\bsa;\bsu)=\supp\,\mu$ also for $u_n=\frac{a_n}{|a_n|}$,
where Theorem \ref{SLP-C-2} predicts that
$\limss_n\Sigma_n(\bsa;\bsu)\subset\supp\,\mu\cup\{1\}$. However, contrary to the case of
constant Schur parameters, a weak limit point, $-1$, appears now outside $\supp\,\mu$.

It is interesting to compare the first example, $a_n=\frac{1}{2}$, with the situation of
a random sequence $\bsa$ lying on $\re(z)\geq\frac{1}{2}$. Figure 7, which represents
this last case for $\bsu=\bsu^{-1}$, agrees with Corollary \ref{DERIVED11}, which
predicts that $\{\supp\,\mu\}'\subset\triangle_\alpha(1)$, $\alpha\approx0.24\pi$.

The second example of 2-periodic Schur parameters, $a_{2n-1}=\frac{1}{4}$,
$a_{2n}=\frac{3}{4}$, can be compared with Figures 8 and 9. Figure 8 deals with the
choice $\bsu=\bsu^1$ for a sequence $\bsa$ whose odd and even subsequences are randomly
located on $\re(z)\leq\frac{1}{4}$ and $\re(z)\geq\frac{3}{4}$ respectively. This figure
confirms Theorem \ref{DERIVED1}, which implies that
$\{\supp\,\mu\}'\subset\triangle_\alpha(-1)$, $\alpha=\frac{\pi}{6}$. In Figure 9 the odd
subsequence of $\bsa$ is randomly chosen on the semicircle
$\overline\Gamma_\frac{\pi}{2}(\frac{1}{4})$ and $a_{2n}=\frac{3}{4}$ with $\bsu=\bsu^i$.
The figure is compatible with theorems \ref{DERIVED31} and \ref{DERIVED1} which give
$\{\supp\,\mu\}'\subset\triangle_\alpha(-1)\cap\triangle_\beta(1)$,
$\alpha=\frac{\pi}{6}$, $\beta\approx0.24\pi$.

Finally, figures 10, 11 and 12 correspond to different sequences of
Schur parameters having two different limit points
$re^{\pm i\frac{\pi}{3}}$ with equal modulus $r=\sin(\frac{3\pi}{8})$.
In figures 10 and 11, the subsequences of $\bsa$ converging to such limit
points are chosen so that $(a_{n+1}/a_n)_{n\geq1}$ has three limit points
$1,e^{\pm i\frac{2\pi}{3}}$. Then, Theorem \ref{DERIVED22} predicts that
$\{\supp\,\mu\}'$ is included in three arcs centred at $-1,e^{\pm i\frac{\pi}{3}}$
with angular radius $\frac{\pi}{4}$. As Example \ref{EXDER2}
shows, this means that $\{\supp\,\mu\}'$ has at least three gaps centred
at $1$ and $e^{\pm i\frac{2\pi}{3}}$ with an angular radius greater than or equal to
$\alpha=\frac{\pi}{12}$. In fact, from Example \ref{EXDER1} we see
that Corollary \ref{DERIVED11} ensures that the radius of the gap around $1$
is not less than $\beta\approx0.22\pi$ and, hence,
$\{\supp\,\mu\}'\subset
\triangle_\alpha(e^{i\frac{2\pi}{3}})\cap\triangle_\alpha(e^{-i\frac{2\pi}{3}})
\cap\triangle_\beta(1)$.
This result agrees with figures 10 and 11. The comparison with Figure 12
is of interest. It represents the case of 2-periodic Schur parameters with the
same limit points as in figures 10 and 11. In this case $(a_{n+1}/a_n)_{n\geq1}$
has only two limit points $e^{\pm i\frac{2\pi}{3}}$, hence, the arc around 1 is
now free of $\{\supp\,\mu\}'$. In fact, we know that
$\{\supp\,\mu\}'=\triangle_{\alpha_-}(-1)\cap\triangle_{\alpha_+}(+1)$,
$\alpha_-\approx0.59\pi$, $\alpha_+\approx0.31\pi$. Notice also the similarity
between figures 10, 11 and 12 concerning the isolated mass point close to
$e^{i\frac{2\pi}{3}}$.

\begin{figure}
\includegraphics{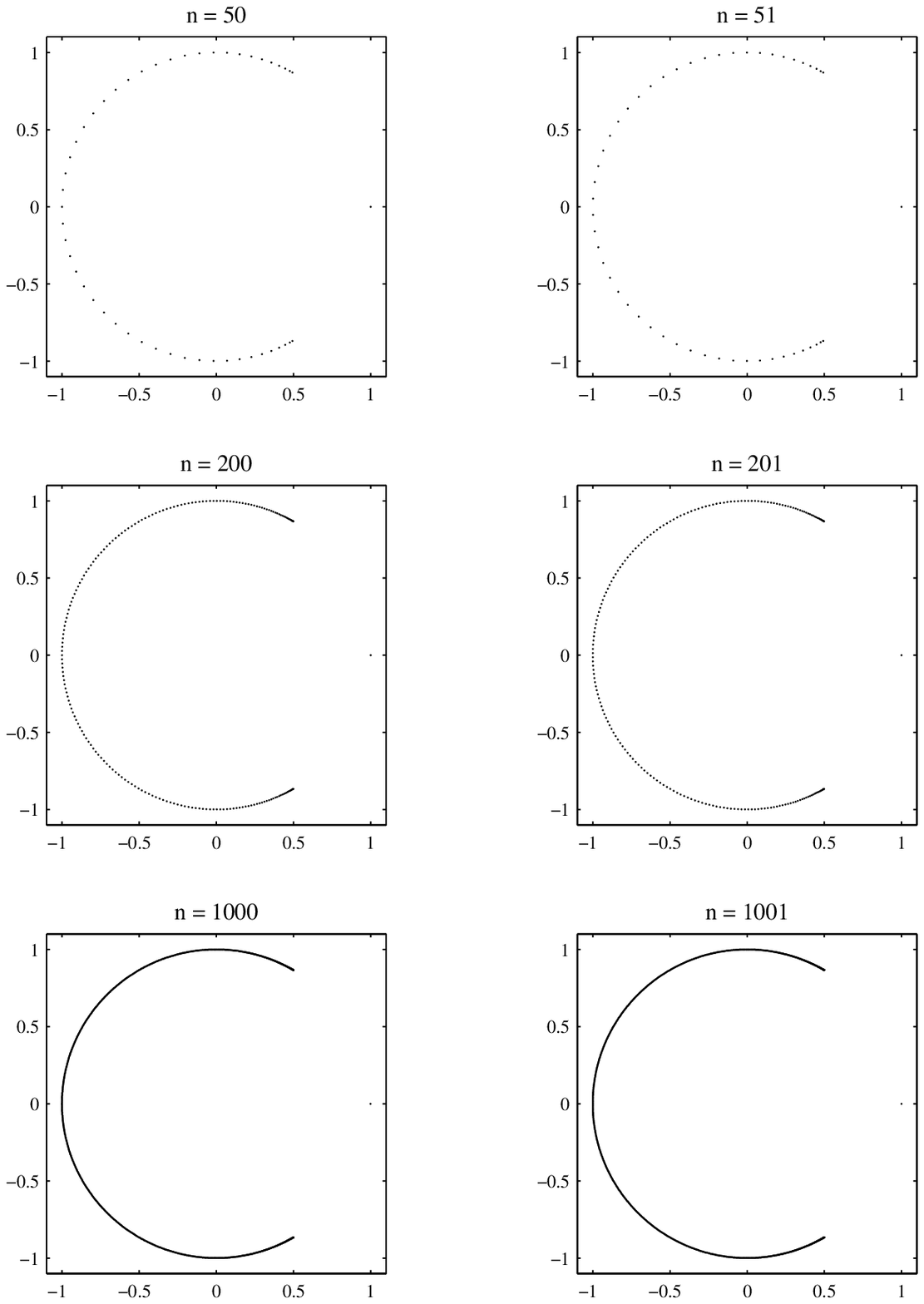}
\caption{$\Sigma_n(\bsa;\bsu^1)$ for $a_n=\frac{1}{2}$.}
\end{figure}

\begin{figure}
\includegraphics{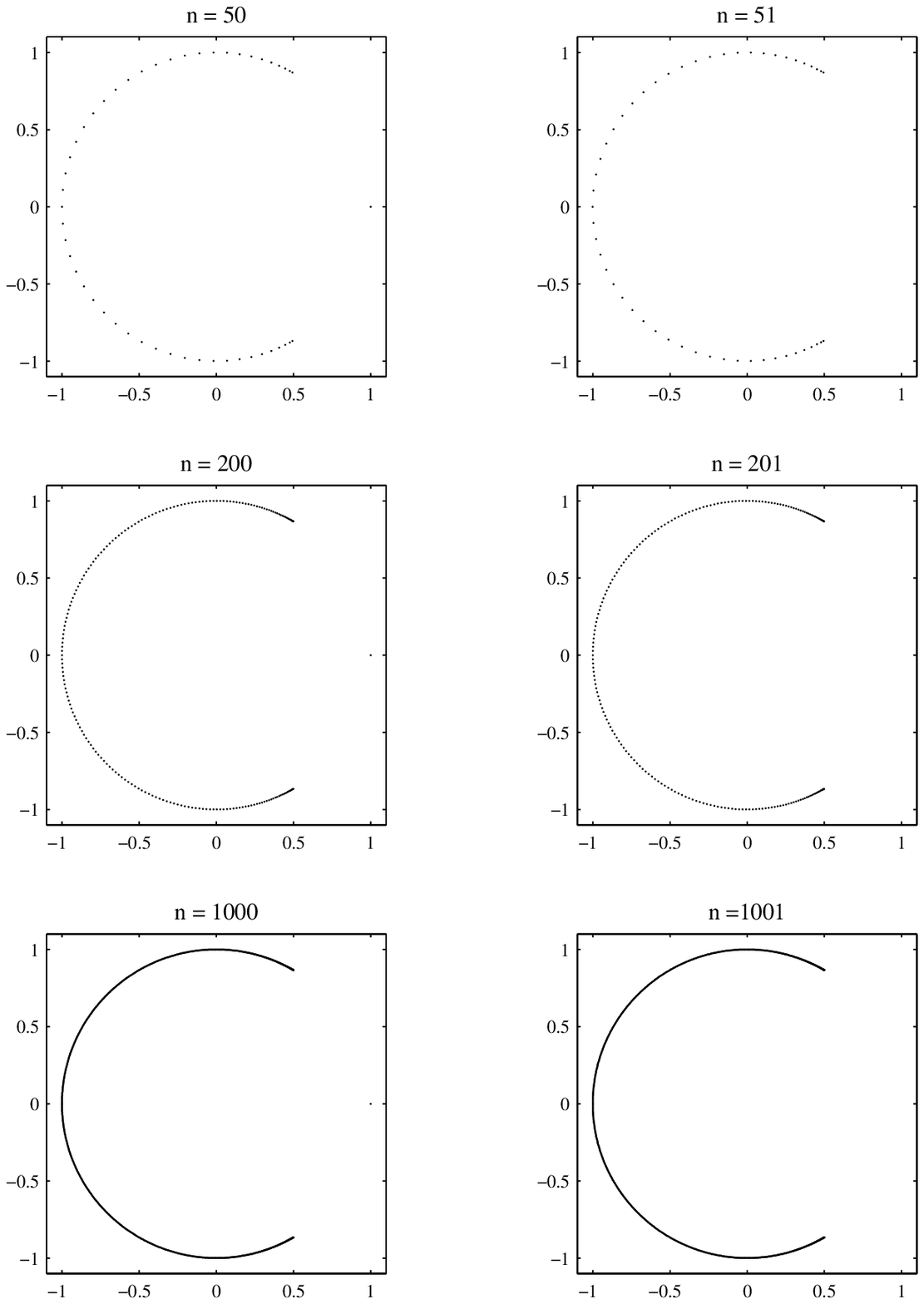}
\caption{$\Sigma_n(\bsa;\bsu^{-1})$ for $a_n=\frac{1}{2}$.}
\end{figure}

\begin{figure}
\includegraphics{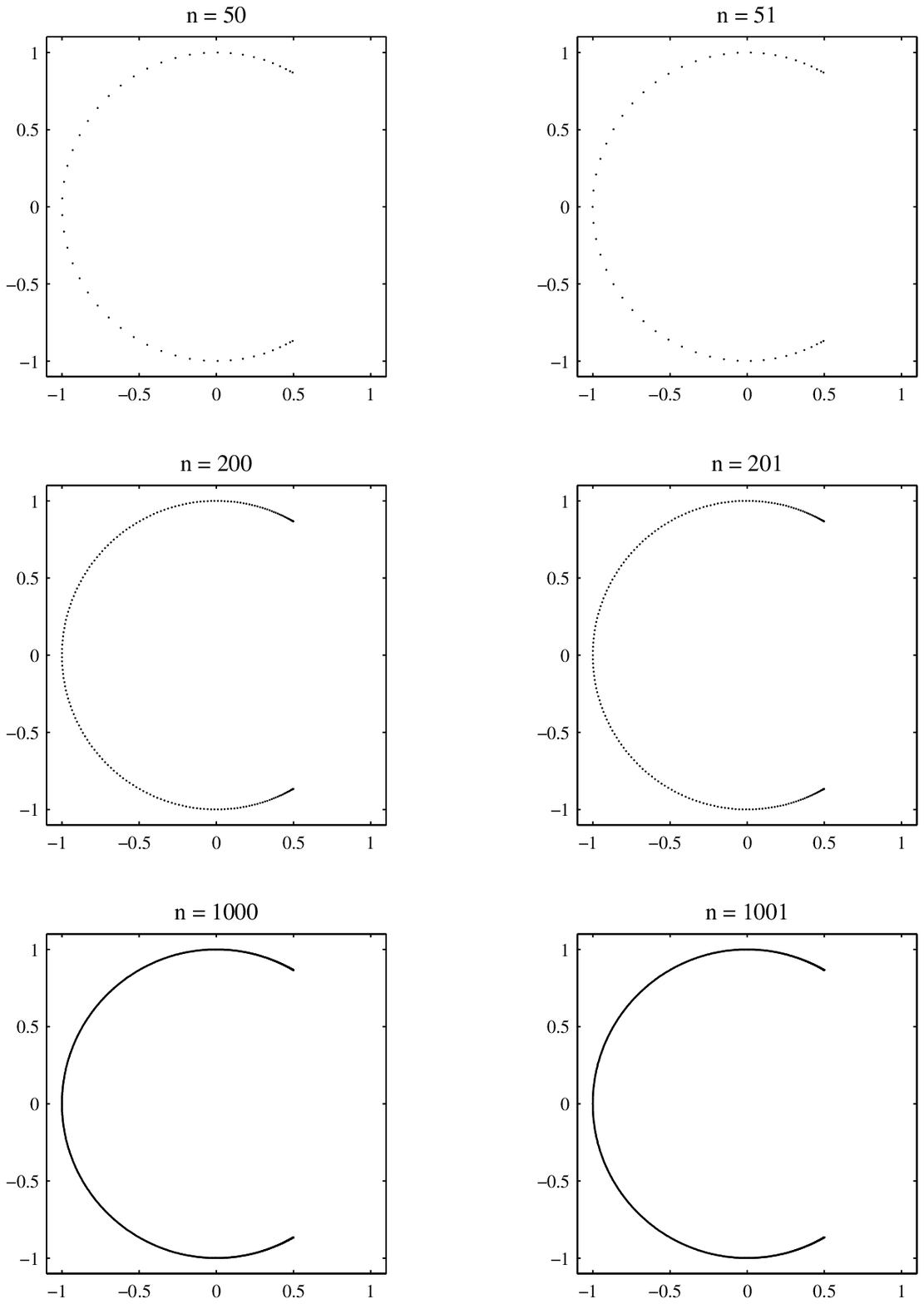}
\caption{$\Sigma_n(\bsa;\bsu)$ for $a_n=\frac{1}{2}$ and $u_n=\frac{a_n}{|a_n|}$.}
\end{figure}

\begin{figure}
\includegraphics{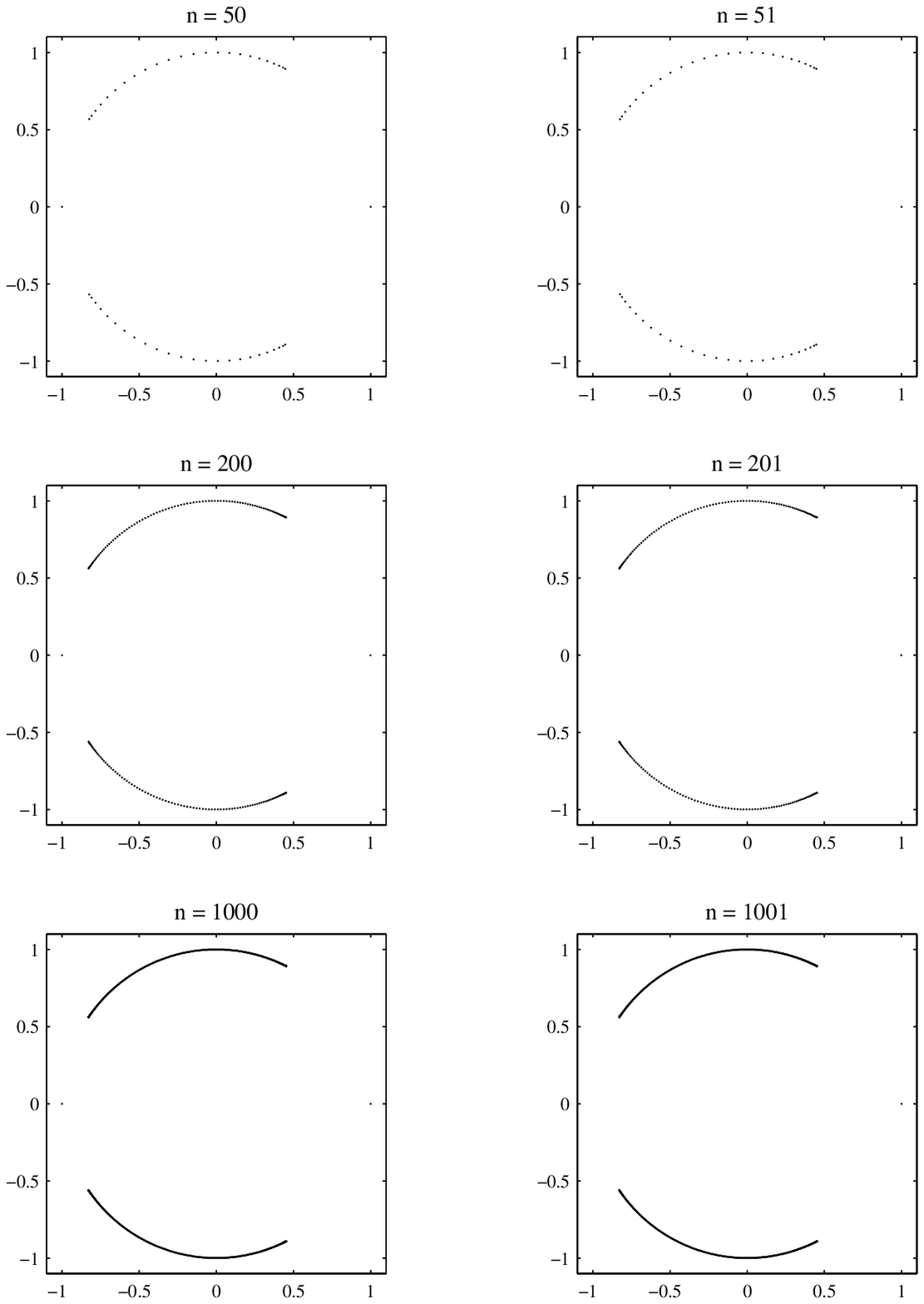}
\caption{$\Sigma_n(\bsa;\bsu^1)$ for $a_{2n-1}=\frac{1}{4}$ and $a_{2n}=\frac{3}{4}$.}
\end{figure}

\begin{figure}
\includegraphics{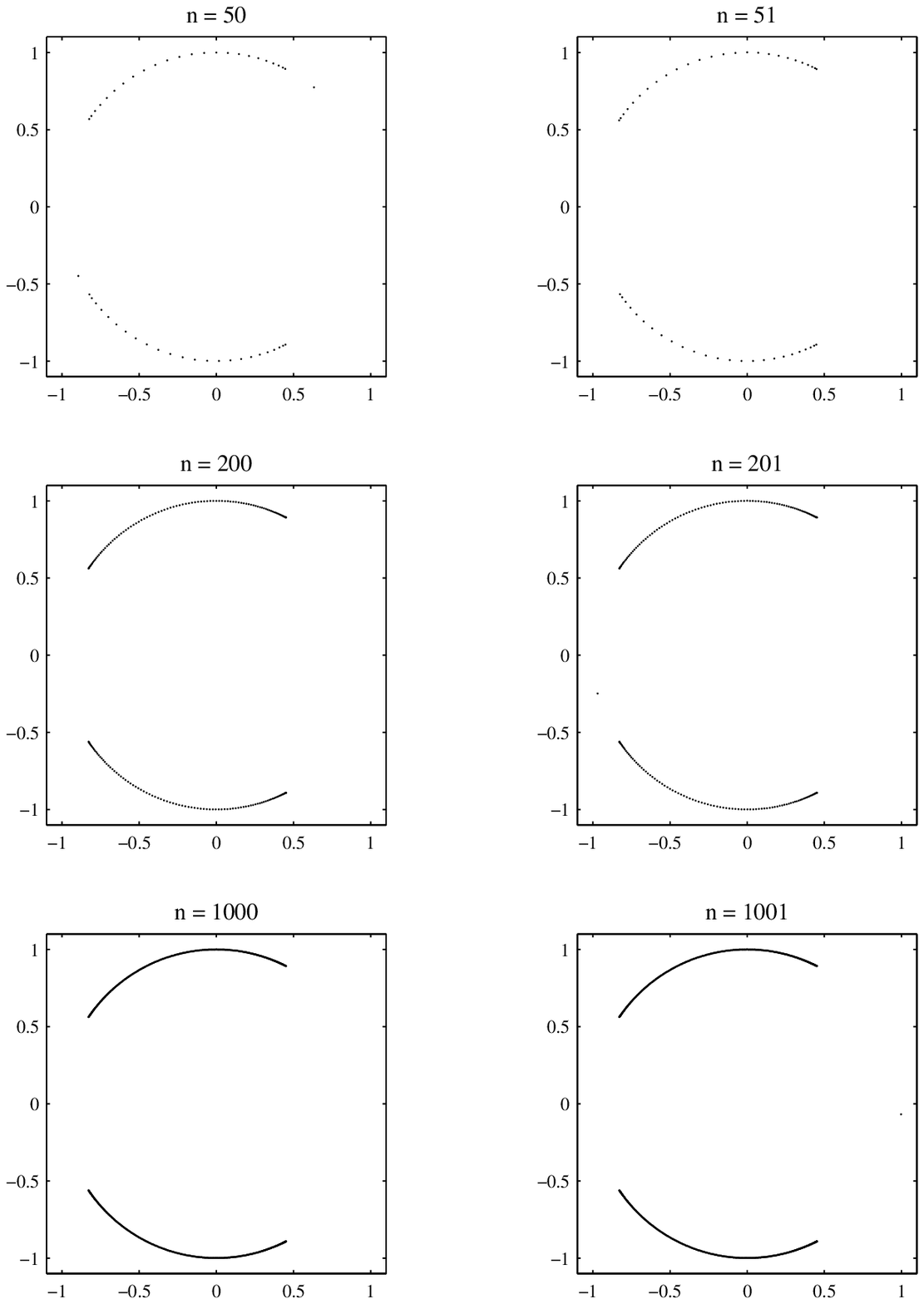}
\caption{$\Sigma_n(\bsa;\bsu^i)$ for $a_{2n-1}=\frac{1}{4}$ and $a_{2n}=\frac{3}{4}$.}
\end{figure}

\begin{figure}
\includegraphics{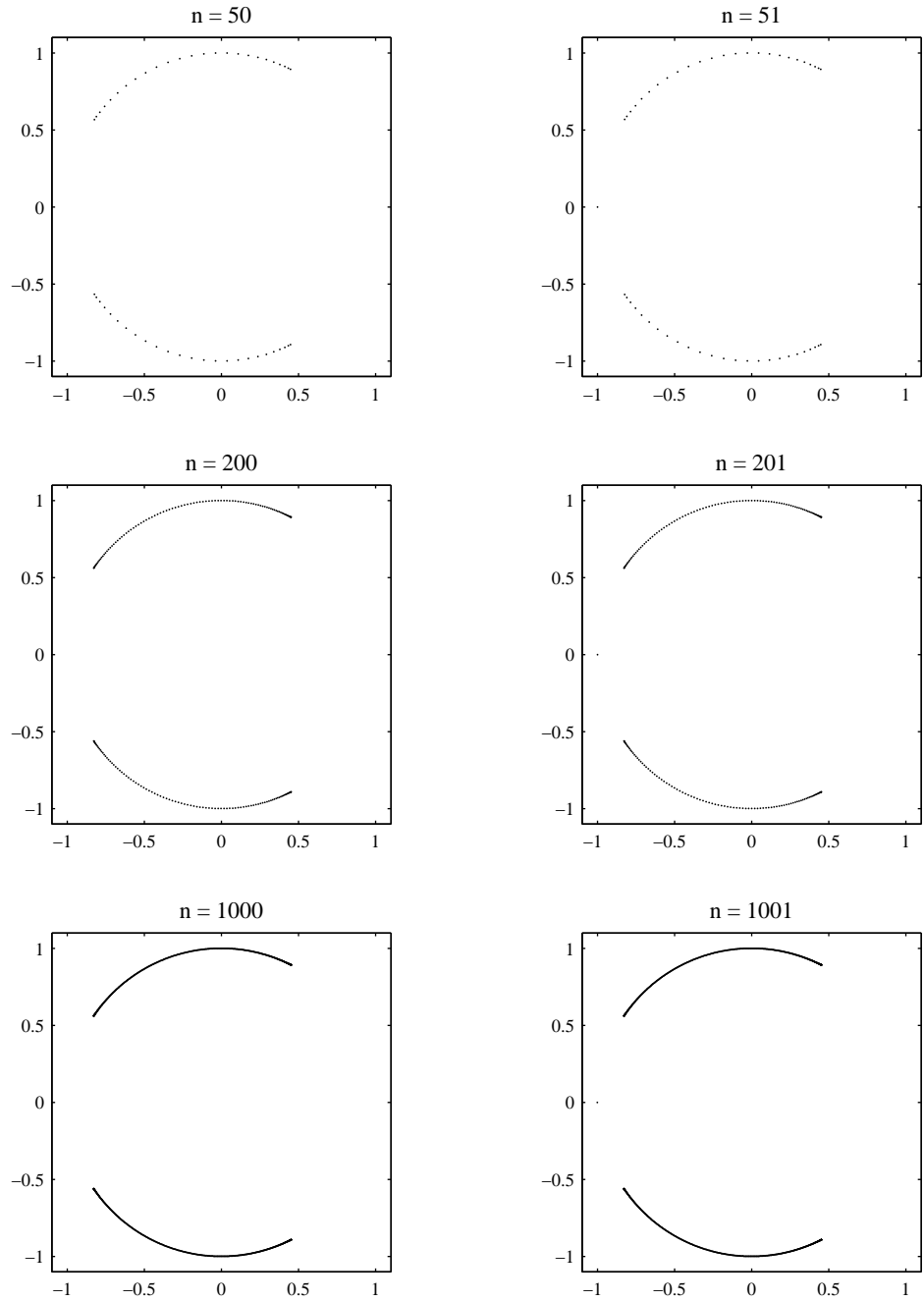}
\caption{$\Sigma_n(\bsa;\bsu)$ for $a_{2n-1}=\frac{1}{4}$, $a_{2n}=\frac{3}{4}$ and
$u_n=\frac{a_n}{|a_n|}$.}
\end{figure}

\begin{figure}
\includegraphics{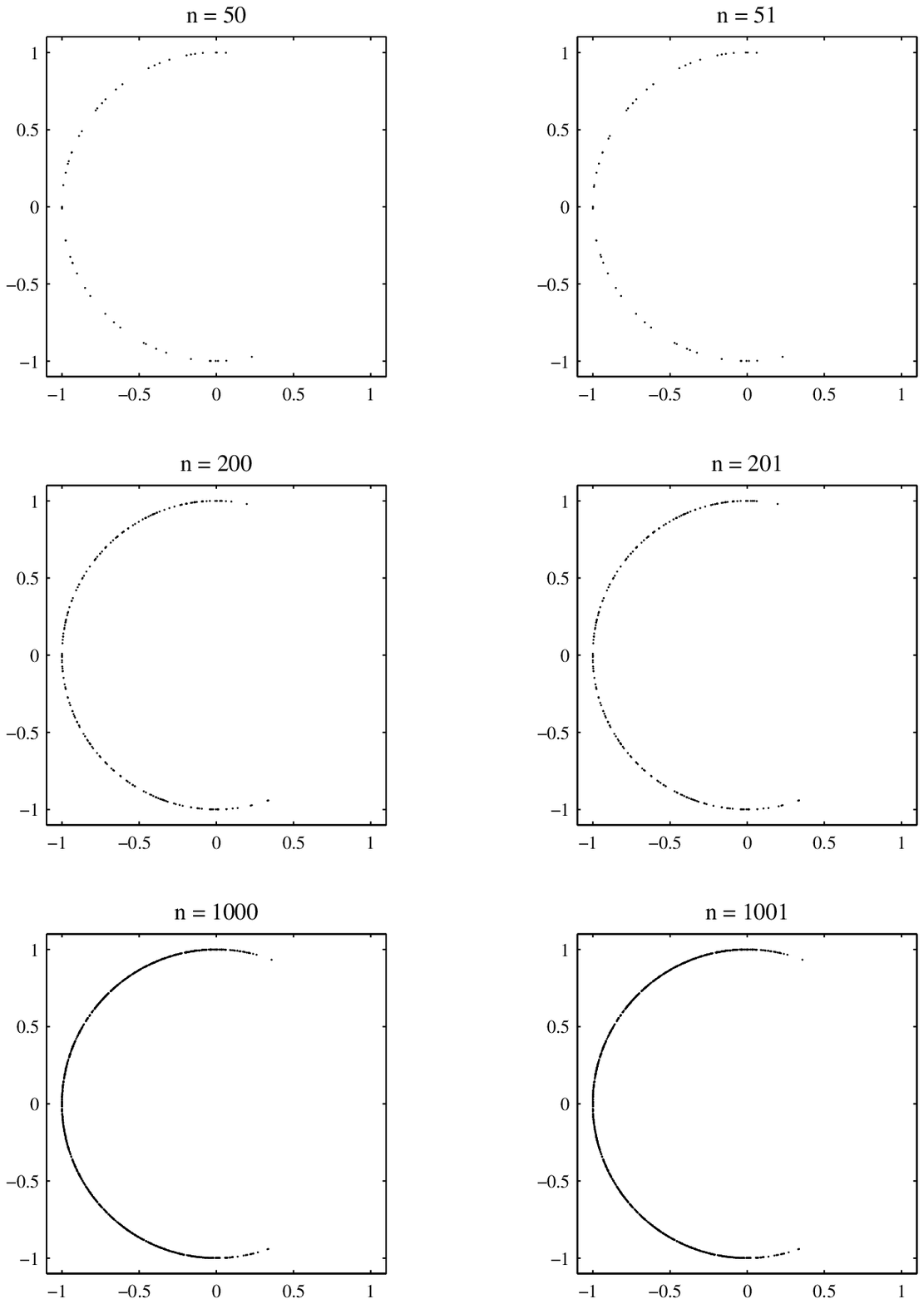}
\caption{$\Sigma_n(\bsa;\bsu^{-1})$ for $a_n$ randomly distributed on
$\re(z)\geq\frac{1}{2}$.}
\end{figure}

\begin{figure}
\includegraphics{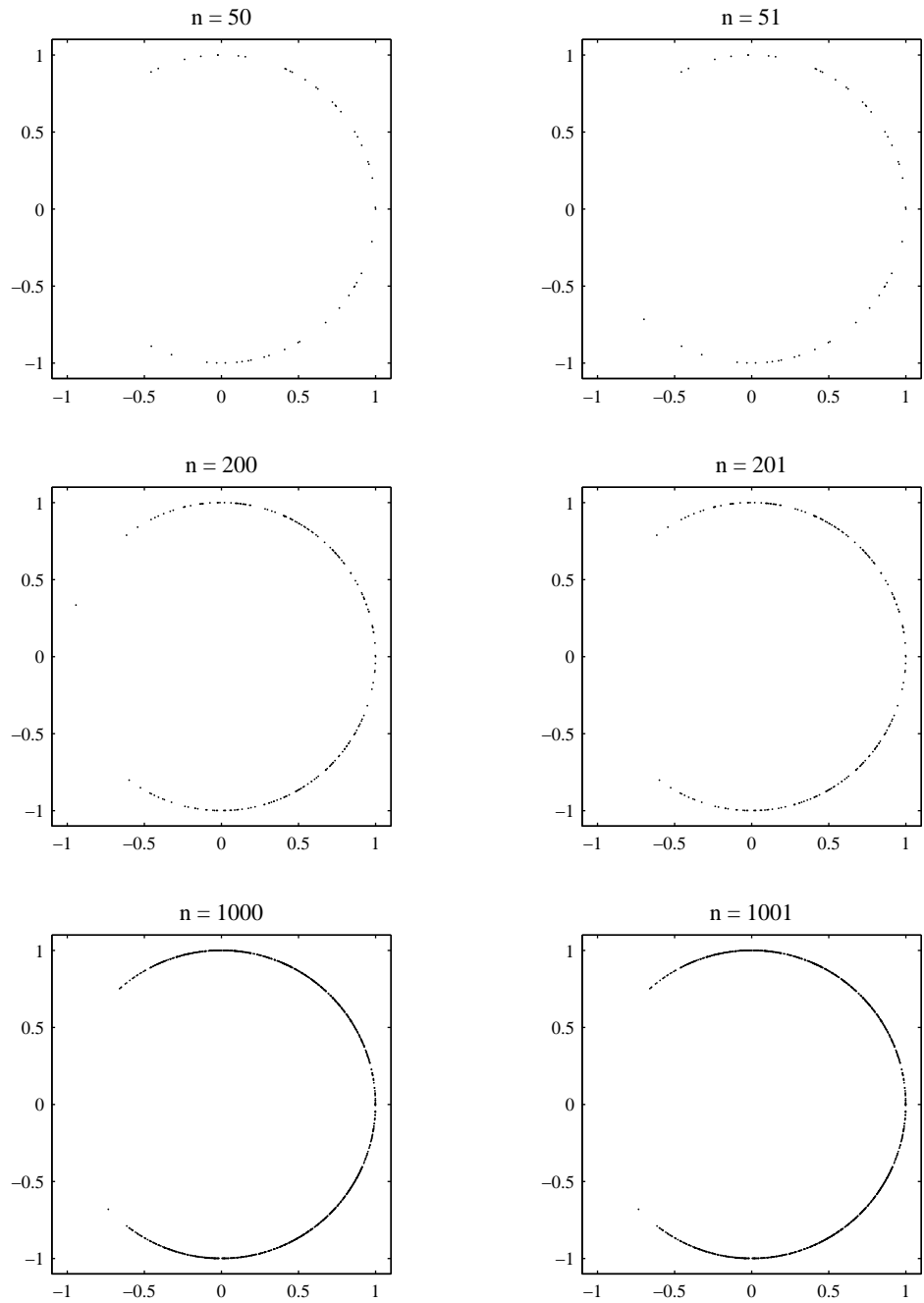}
\caption{$\Sigma_n(\bsa;\bsu^1)$ for $a_{2n-1}$ and $a_{2n}$ randomly distributed on
$\re(z)\leq\frac{1}{4}$ and $\re(z)\geq\frac{3}{4}$ respectively.}
\end{figure}

\begin{figure}
\includegraphics{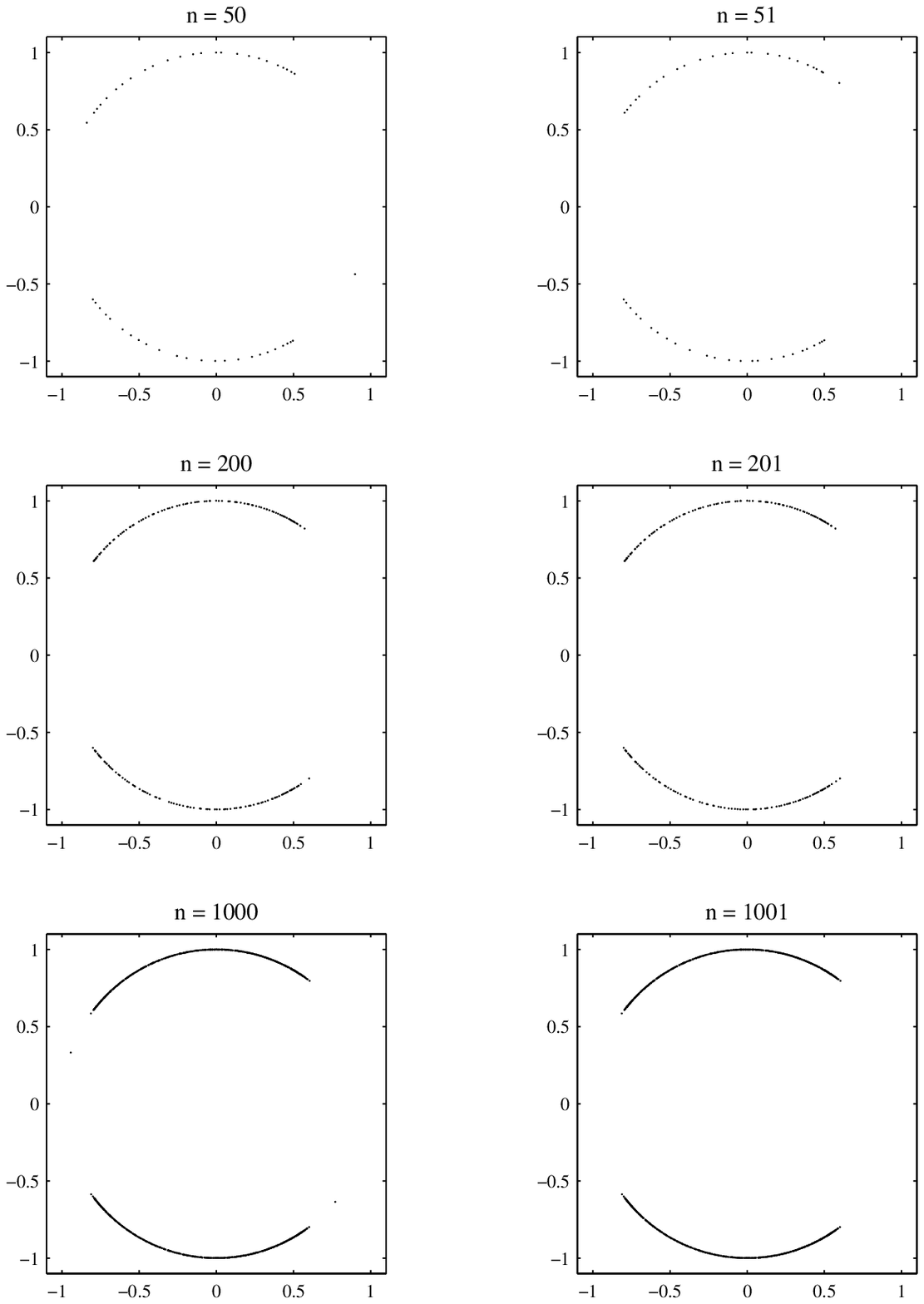}
\caption{$\Sigma_n(\bsa;\bsu^i)$ for $a_{2n-1}$ randomly distributed on
$\overline\Gamma_{\frac{\pi}{2}}(\frac{1}{4})$ and $a_{2n}=\frac{3}{4}$.}
\end{figure}

\begin{figure}
\includegraphics{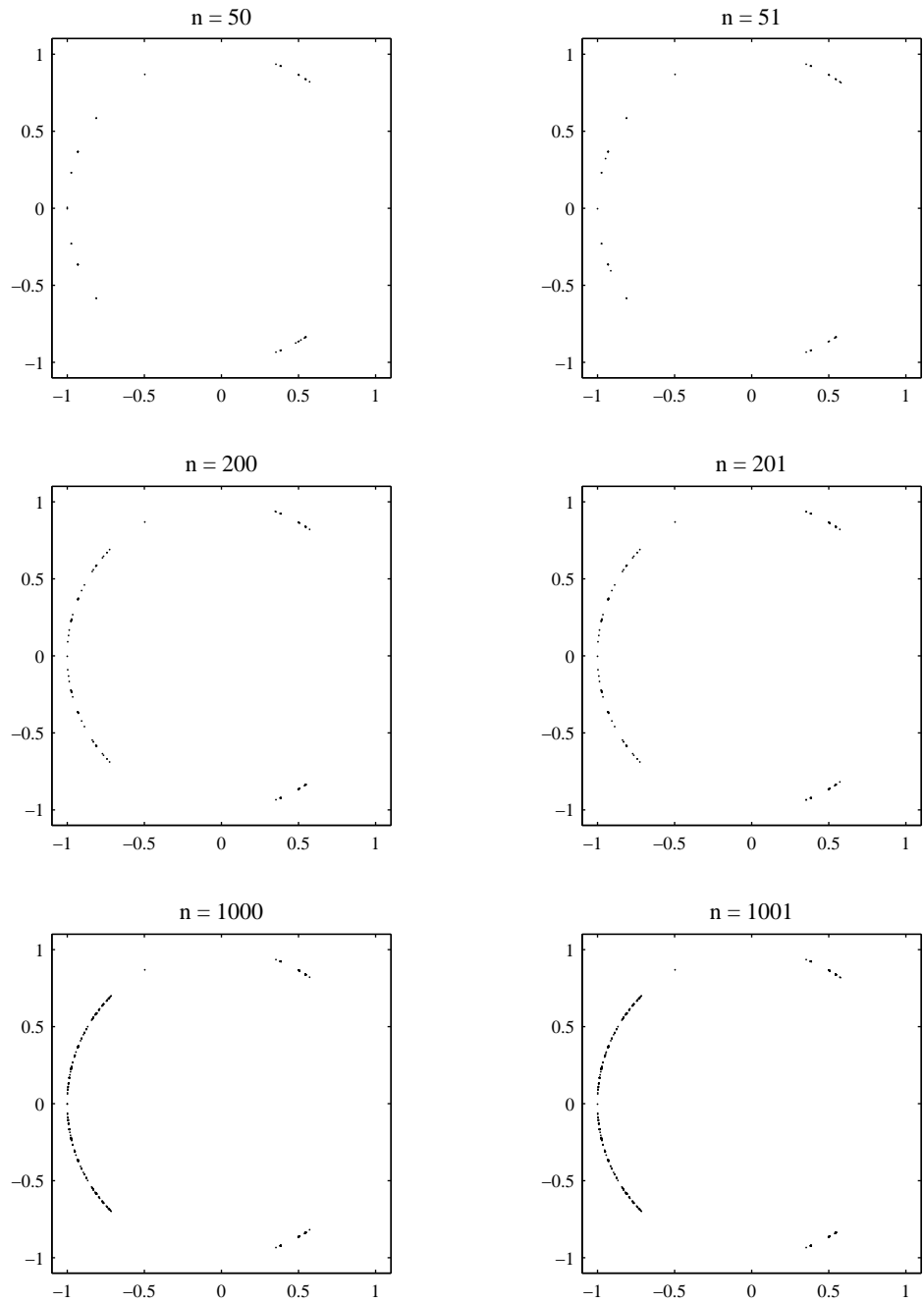}
\caption{$\Sigma_n(\bsa;\bsu^w)$ for $a_n=\sin(\frac{3\pi}{8})e^{\pm i\frac{\pi}{3}}$ if
$n$ is prime/not prime and $w=e^{i\frac{\pi}{3}}$.}
\end{figure}

\begin{figure}
\includegraphics{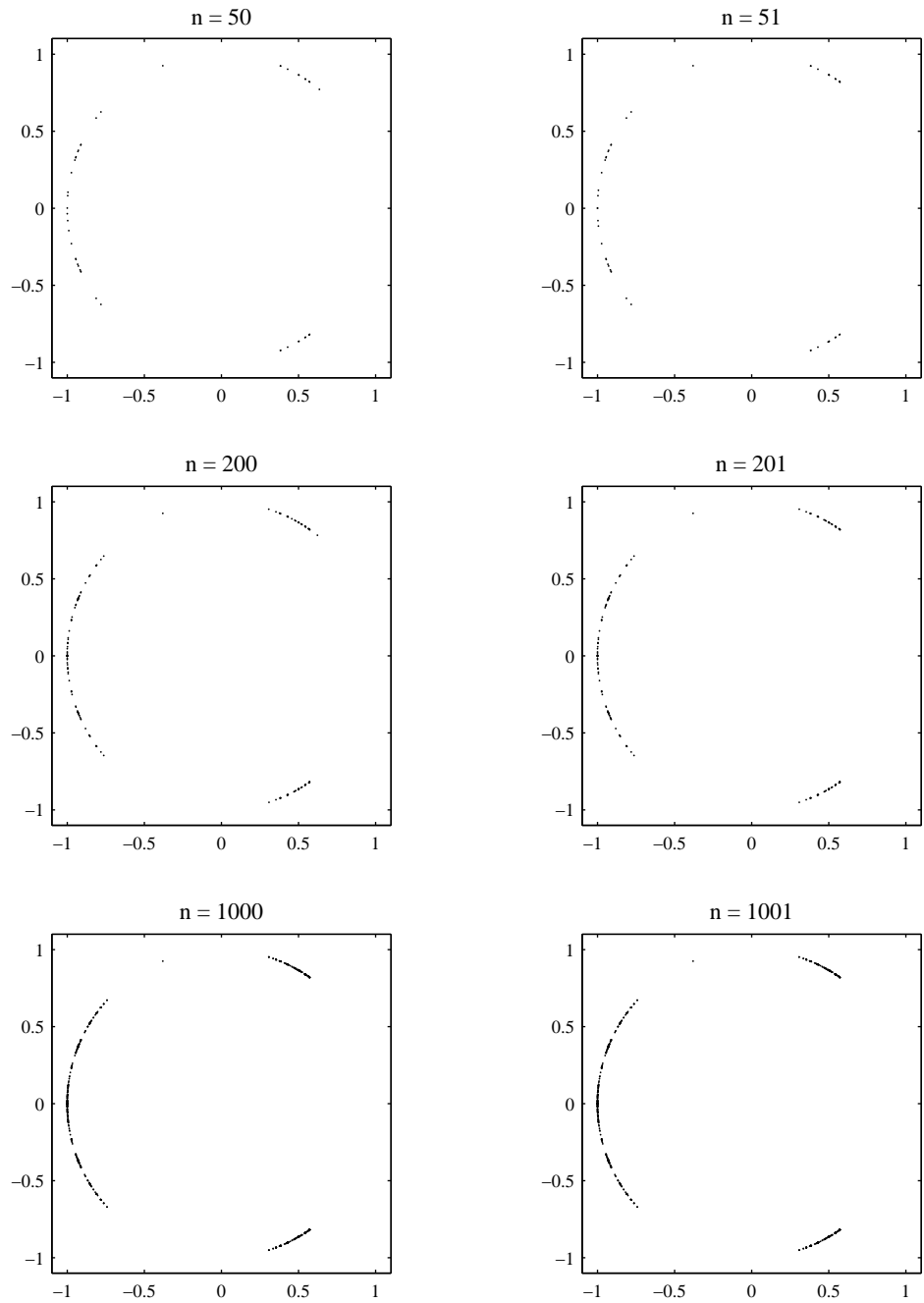}
\caption{$\Sigma_n(\bsa;\bsu^w)$ for $a_n$ randomly distributed on
$\{\sin(\frac{3\pi}{8})e^{\pm i\frac{\pi}{3}}\}$ and $w=e^{i\frac{\pi}{3}}$.}
\end{figure}

\begin{figure}
\includegraphics{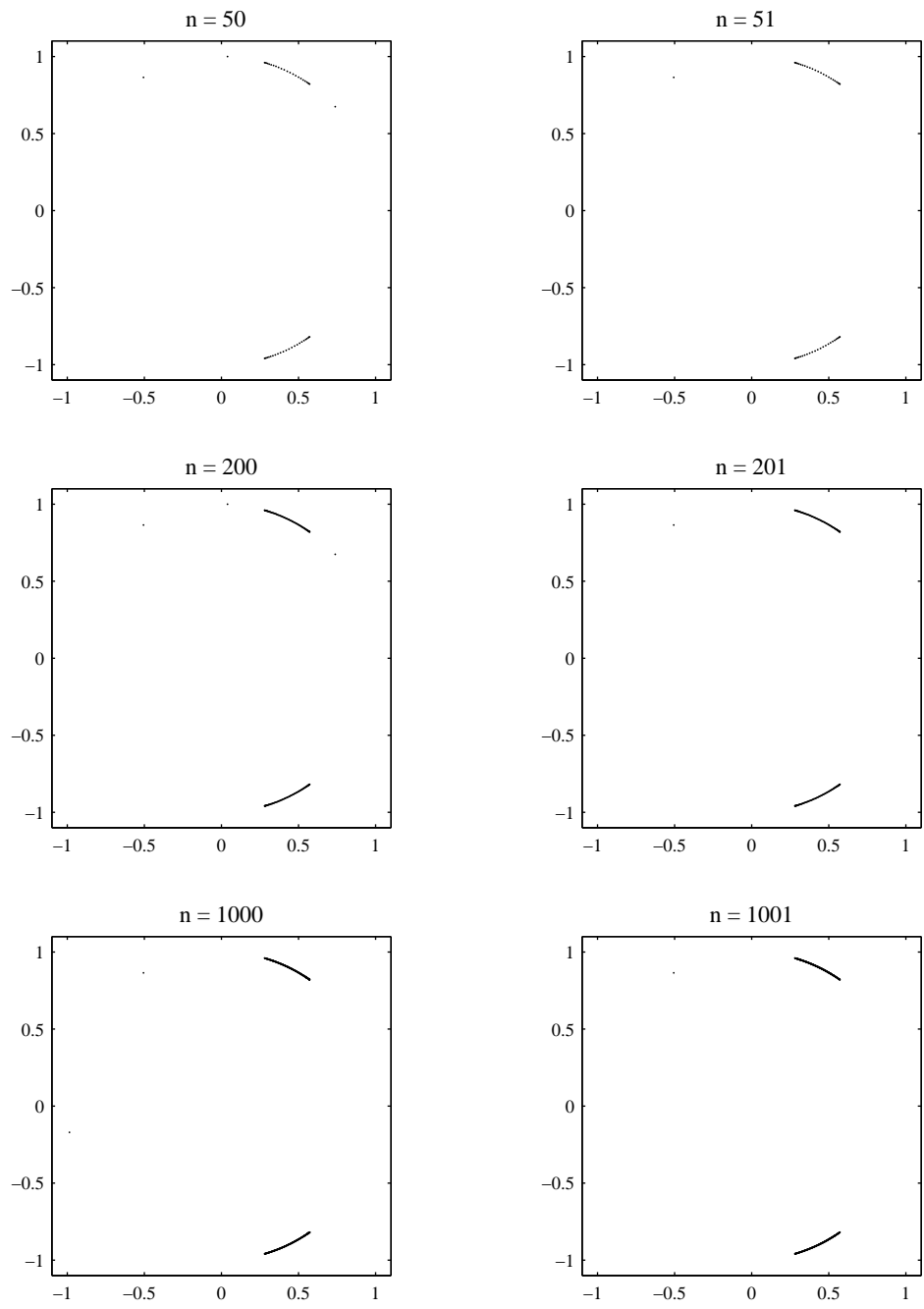}
\caption{$\Sigma_n(\bsa;\bsu^w)$ for
$a_n=\sin(\frac{3\pi}{8})e^{\pm i\frac{\pi}{3}}$ if $n$ is even/odd
and $w=e^{i\frac{\pi}{3}}$.}
\end{figure}

\vfil\eject

%%%%%%%%%%%%%%%%%%%%%%%%%%%%%%%%%%%%%%%%%%%%%%%%%%%%%%%%%%%%%%%%%%%%%%%%%%%%%%%%%%%%%
\section{Applications to the continued fractions}
%%%%%%%%%%%%%%%%%%%%%%%%%%%%%%%%%%%%%%%%%%%%%%%%%%%%%%%%%%%%%%%%%%%%%%%%%%%%%%%%%%%%%

In this section we will show some applications of the previous results to the study of
rational approximants of Carath\'eodory functions. In what follows,
$f_n(z) \rightrightarrows f(z)$, $z\in\Omega$, means that the sequence $(f_n)_{n\geq1}$
uniformly converges to $f$ on compact subsets of $\Omega$.

It is known that the monic orthogonal polynomials $(\Phi_n)_{n\geq0}$ corresponding to a
measure $\mu$ on $\T$ and the related monic second kind polynomials $(\Psi_n)_{n\geq0}$
provide rational approximants for the associated Carath\'eodory function
$$
F_\mu(z) := \int_\T \frac{\lambda+z}{\lambda-z} \, d\mu(\lambda).
$$
More precisely, $\Psi_n^*(z)/\Phi_n^*(z) \rightrightarrows F_\mu(z)$, $z\in\D$, and
$-\Psi_n(z)/\Phi_n(z) \rightrightarrows F_\mu(z)$, $z\in\C\backslash\overline\D$
\cite{JoNjTh89}. When $\supp\,\mu\neq\T$, it is possible to enlarge the above domains of
convergence through a careful analysis of the asymptotic behaviour of the zeros of the
orthonormal polynomials \cite[Section 9]{Kh02}.

All these can be read as results about the convergence of continued fractions. Remember
that, given a continued fraction
$$
K := \beta_0+\frac{\alpha_1}{\beta_1+\frac{\ds\alpha_2}{\ds\beta_2+
{\mathop{}\limits_{\ddots}}}},
$$
the related $w$-modified $n$-th approximant is
$$
K_n^w := \beta_0+\frac{\alpha_1}{\beta_1+\frac{\ds\alpha_2}{\ds\beta_2+
{\mathop{}\limits_{\ddots\mathop{}\limits_{\ds+\frac{\alpha_n}{\beta_n+w}}}}}}.
$$
In particular, $K_n:=K_n^0$ is called the $n$-th approximant of $K$. It is known that
(see \cite{Wa48})
\beq \label{CF}
K_n^w=\frac{A_n+wA_{n-1}}{B_n+wB_{n-1}}, \quad n\geq0,
\eeq
where $A_n$ and $B_n$ are given by the same recurrence
$$
X_n=\beta_nX_{n-1}+\alpha_nX_{n-2}, \quad X_n=A_n,B_n, \quad n\geq0,
$$
but with different initial conditions $A_0=\beta_0$, $A_{-1}=1$ and $B_0=1$,
$B_{-1}=0$.

The recurrences for $\Phi_n$ and $\Psi_n$ show that $A_{2n}=\Psi_n^*(z)$,
$B_{2n}=\Phi_n^*(z)$, $A_{2n+1}=-z\Psi_n(z)$ and $B_{2n+1}=z\Phi_n(z)$ for the
continued fraction  \cite{JoNjTh89}
\beq \label{PC}
K(\bsa;z) :=
1+\frac{-2z}{\ds z+\frac{\ds 1}{\ds \overline a_1+\frac{\rho_1^2z}{a_1z+\frac{\ds 1}
{\ds \overline a_2+\frac{\rho_2^2z}{a_2z+\mathop{}\limits_{\ddots}}}}}},
\eeq
where $\bsa$ is the sequence of Schur parameters of $\mu$. Hence,
$\Psi_n^*(z)/\Phi_n^*(z)$ and $-\Psi_n(z)/\Phi_n(z)$ are respectively the $2n$-th
approximant $K_{2n}(a_1,\dots,a_n;z)$ and the $2n+1$-th approximant
$K_{2n+1}(a_1,\dots,a_n;z)$ of $K(\bsa;z)$.

It is also clear from (\ref{CF}) that for any $u\in\T$ \cite{JoNjTh89}
\beq \label{MODAPPROX}
\kern-10pt
K_{2n}(a_1,\dots,a_{n-1},u;z)=K_{2n-1}^u(a_1,\dots,a_{n-1};z)=
-\frac{\Psi_n^u(z)}{\Phi_n^u(z)}, \kern7pt n\geq1,
\eeq
where $\Phi_n^u(z):=z\Phi_{n-1}(z)+u\Phi_{n-1}^*(z)$ and
$\Psi_n^u(z):=z\Psi_{n-1}(z)-u\Psi_{n-1}^*(z)$. Therefore, these modified approximants
are quotients of para-orthogonal polynomials. In fact, given a sequence $\bsu$ in $\T$,
the convergence properties for the modified approximants $-\Psi_n^{u_n}/\Phi_n^{u_n}$
are, in general, better than for the standard ones, since it is known that
$-\Psi_n^{u_n}(z)/\Phi_n^{u_n}(z) \rightrightarrows F_\mu(z)$, $z\in\C\backslash\T$
\cite{JoNjTh89}. The aim of this section is to find information about the convergence of
these modified approximants on the unit circle.

Closely related to the concept of Carath\'eodory function is the notion of resolvent
$R_z(T):=(z-T)^{-1}$ of an operator $T\in\frB(H)$, which is again a bounded
operator on $H$ for $z\in\C\backslash\spec(T)$. Moreover, when $T$ is normal,
$\|R_z(T)\|=1/d(z,\spec(T))$ for $z\in\C\backslash\spec(T)$. The Carath\'eodory function
of a measure $\mu$ on $\T$ with Schur parameters sequence $\bsa$ is related to the
resolvent $R_z(\bsa):=R_z(C(\bsa))$, which is a bounded operator on $\ell^2$ for
$z\in\C\backslash\supp\,\mu$. In fact,
$$
\int_\T\lambda^n\,d(E_{C(\bsa)}(\lambda)e_1,e_1)=(C(\bsa)^ne_1,e_1)=
({U^\mu}^n1,1)=\int_\T\lambda^n\,d\mu(\lambda), \kern5pt \forall n\in\Z,
$$
and, thus, $d\mu(\lambda)=d(E_{C(\bsa)}(\lambda)e_1,e_1)$. Therefore,
$$
F_\mu(z)=\int_\T\frac{\lambda+z}{\lambda-z}\,d(E_{C(\bsa)}(\lambda)e_1,e_1)=
1-2z(R_z(\bsa)e_1,e_1).
$$

Also, for any $u\in\T$, the modified approximant $-\Psi_n^u/\Phi_n^u$
is related to the resolvent $R_z(a_1,\dots,a_{n-1},u):=R_z(C(a_1,\dots,a_{n-1},u))$,
which defines an operator on $\ell^2_n$ for $z$ outside the spectrum of
$C(a_1,\dots,a_{n-1},u)$. More precisely, if
$f_j:= 1-2z(R_z(a_1,\dots,a_{n-1},u)e_1,e_j)$ for $j=1,\dots,n$,
the vector $\boldsymbol{f}:=\sum_{j=1}^nf_je_j$ satisfies
$$
(C(a_1,\dots,a_{n-1},u)-z)\boldsymbol{f}=(C(a_1,\dots,a_{n-1},u)+z)e_1.
$$
Just solving this system for $f_1$ we get
$$
f_1 = -\frac{\det(z-VC(-a_1,\dots,-a_{n-1},-u)V^*)}
{\det(z-C(a_1,\dots,a_{n-1},u))},
$$
where $V$ is the linear operator on $\ell^2_n$ defined by $Ve_j=(-1)^je_j$,
$j=1,\dots,n$. From Corollary \ref{NORMTRUNC-C},
$\Phi_n^u(z)=\det(z-C(a_1,\dots,a_{n-1},u))$, so, we finally get
$$
-\frac{\Psi_n^u(z)}{\Phi_n^u(z)} = 1-2z(R_z(a_1,\dots,a_{n-1},u)e_1,e_1).
$$

As a consequence of the previous discussion, given a sequence $\bsu$ in $\T$, the weak
convergence of $(\hat R_z(a_1,\dots,a_{n-1},u_n))_{n\geq1}$ to $R_z(\bsa)$ implies the
convergence of $(-\Psi_n^{u_n}(z)/\Phi_n^{u_n}(z))_{n\geq1}$ to $F_\mu(z)$. In the case
of self-adjoint band operators, the convergence of the resolvents of finite orthogonal
truncations was analyzed in \cite{BaLoMaTo99} and \cite{IfPa01a}, in connection with its
interest for the Jacobi fractions. An extension of the ideas in \cite{BaLoMaTo99}
and \cite{IfPa01a} gives the following result.

\bp \label{CONV-RES}

Let $T\in\frB(\ell^2)$ be a normal band operator. If $T_n$ is a normal truncation of $T$
on $\ell^2_n$ for $n\geq1$ and $(\|T_n\|)_{n\geq1}$ is bounded, for all $x\in\ell^2$,
$$
\hat R_z(T_n)x \rightrightarrows R_z(T)x, \quad z\in\C\backslash\lims_n\spec(T_n).
$$
Moreover, each $z\in\lims_n\spec(T_n)\backslash\limi_n\spec(T_n)$ has a neighbourhood
where the above uniform convergence holds at least for a subsequence of $(T_n)_{n\geq1}$.

\ep

\bpr

Let $z\in\C\backslash\limi_n\spec(T_n)$. There exist $\delta>0$ and a subsequence
$(T_n)_{n\in\cI}$ such that $d(z,\spec(T_n))\geq\delta$, $\forall n\in\cI$. Hence,
$D_\delta(z)\subset\C\backslash\limi_n\spec(T_n)$ and, from Proposition
\ref{GENTRUNC}.1, $D_\delta(z)\subset\C\backslash\spec(T)$. Therefore, $R_w(T)\in\frB(\ell^2)$
for $w \in D_\delta(z)$. Also, $R_w(T_n)$ exists for $n\in\cI$ and $w \in D_\delta(z)$.
Moreover, since $T_n$ is normal, $\|R_w(T_n)\|=1/d(w,\spec(T_n))\leq1/(\delta-|w-z|)$
for $n\in\cI$. Hence, $(\|R_w(T_n)\|)_{n\in\cI}$ is uniformly bounded with respect to $w$
on compact subsets of $D_\delta(z)$.

Let $P_n$ be the orthogonal projection on $\ell^2_n$ and $w \in D_\delta(z)$.
From the identities
$\hat R_w(T_n)=P_nR_w(\hat T_n)$ and $R_w(\hat T_n)-R_w(T)=R_w(\hat T_n)(\hat T_n-T)R_w(T)$
we get
$$
\hat R_w(T_n)-R_w(T) = \hat R_w(T_n)(\hat T_n-T)R_w(T) + (P_n-1)R_w(T).
$$
Proposition \ref{GENTRUNC}.1 states that $\hat T_n \to T$. Since $P_n\to1$ and
$(\|R_w(T_n)\|)_{n\in\cI}$ is bounded we conclude that
$\hat R_w(T_n) \mathop{\to} \limits_{n\in\cI} R_w(T)$.
Moreover, the equality
\vskip-7pt
$$
\ba{l}
\hat R_{w'}(T_n) - R_{w'}(T) = \hat R_w(T_n) - R_w(T) \, +
\smallskip \\ \kern100pt
+ \, (w'-w)(R_w(T) R_{w'}(T) - \hat R_{w'}(T_n) \hat R_w(T_n))
\ea
$$
shows that, given $x \in H$ and $\epsilon>0$, there is a disk centred at $w$ such that
$\|\hat R_{w'}(T_n)x-R_{w'}(T)x\|<\epsilon$ for $w'$ lying on such a disk and $n$ big enough.
Then, standard arguments prove that
$\hat R_w(T_n)x \mathop{\rightrightarrows} \limits_{n\in\cI} R_w(T)x$, $w \in D_\delta(z)$.

In the preceding discussion, if $z\notin\lims_n\spec(T_n)$, the subsequence
$(T_n)_{n\in\cI}$ can be chosen such that $\cI=\{n\in\N:n\geq N\}$, $N\in\N$, and, so,
the uniform convergence $\hat R_w(T_n)x \rightrightarrows R_w(T)x$, $w \in D_\delta(z)$,
holds for the full sequence. Therefore, the convergence is uniform on compact subsets of
$\C\backslash\lims_n\spec(T_n)$.

\epr

From the preceding proposition a result for the resolvent of $C(\bsa)$ immediately
follows.

\bt \label{CONV-RES-C}

Given a sequence $\bsa$ in $\D$ and a sequence $\bsu$ in $\T$, for all $x\in\ell^2$,
$$
\hat R_z(a_1,\dots,a_{n-1},u_n)x \rightrightarrows R_z(\bsa)x,
\quad z\in\C\backslash\lims_n\Sigma_n(\bsa;\bsu).
$$
Moreover, each $z\in\lims_n\Sigma_n(\bsa;\bsu)\backslash\supp\,\mu$ (up to, at most,
one point) has a neighbourhood where the above uniform convergence holds at least for
a subsequence.

\et

\bpr

Apply Proposition \ref{CONV-RES-C} to $C(\bsa)$ and its finite unitary truncations
$C(a_1,\dots,a_n,u_n)$, taking into account that, from Theorem \ref{SLP-C-1},
$\limi_n\Sigma_n(\bsa;\bsu)$ coincides with $\supp\,\mu$ up to, at most, at one point.

\epr

Since the strong convergence of operators implies the weak convergence, we get a
conclusion for the convergence of the modified approximants (\ref{MODAPPROX}).

\bc \label{CONV-F}

If $\bsa$ is the sequence of Schur parameters of a measure $\mu$ on $\T$ and $\bsu$ is a
sequence in $\T$,
$$
K_{2n}(a_1,\dots,a_{n-1},u_n;z) \rightrightarrows F_\mu(z),
\quad z\in\C\backslash\lims_n\Sigma_n(\bsa;\bsu).
$$
Moreover, each $z\in\lims_n\Sigma_n(\bsa;\bsu)\backslash\supp\,\mu$ (up to, at most, one
point) has a neighbourhood where the above uniform convergence holds at least for a
subsequence.

\ec

This corollary says that $(\Psi_n^{u_n}/\Phi_n^{u_n})_{n\geq1}$ converges to $F_\mu$,
not only outside the unit circle, but also at the points in the unit circle that are not
limit points of the zeros of the para-orthogonal polynomials $(\Phi_n^{u_n})_{n\geq1}$.
The results of the previous sections can now be used to get information about
the convergence of the sequence $(\Psi_n^{u_n}/\Phi_n^{u_n})_{n\geq1}$, as the following
examples show.

\bex \label{EX3}

Let $\bsa$ be the sequence of Schur parameters of a measure $\mu$ on $\T$ and let $\bsu$
be a sequence in $\T$.

\be

\item
Schur parameters converging to the unit circle.

If $\lim_n|a_n|=1$ and $u_n=\frac{a_n}{|a_n|}$, Corollary \ref{LP-C-11} gives
$$
K_{2n}(a_1,\dots,a_{n-1},u_n;z) \rightrightarrows F_\mu(z),
\quad z\in\C\backslash\supp\,\mu.
$$

\item
Rotated asymptotically 2-periodic Schur parameters.

Let $\lim_na_{2n-1}(\lambda)=a_o$, $\lim_na_{2n}(\lambda)=a_e$, $\lambda\in\T$.
Example \ref{EX11} shows that, if $u_n=\frac{a_n}{|a_n|}$, the conditions
$\rho_o\rho_e\mp\re(\overline a_oa_e)<\min\{|a_o|,|a_e|\}$ respectively imply that
$$
K_{2n}(a_1,\dots,a_{n-1},u_n;z) \rightrightarrows F_\mu(z), \quad
z\in\C\backslash(\supp\,\mu\cup\triangle_{\beta_\pm}(\pm\lambda)),
$$
where $\beta_\pm\in(0,\pi]$ are given by
$\cos\frac{\beta_\pm}{2}=
\sqrt{\frac{1+\rho_o\rho_e\mp\re(\overline a_oa_e)}{1+\min\{|a_o|,|a_e|\}}}$.

\item
The limit points of the odd and even subsequences of $\bsa(-\lambda)$, $\lambda\in\T$,
separated by a band.

If $\cB(u,\alpha_1,\alpha_2)$, $u\in\T$, $0\leq\alpha_1<\alpha_2\leq\pi$, is such a band,
theorems \ref{DERIVED1} and \ref{LP-C-2} prove that
$$
K_{2n}(a_1,\dots,a_{n-1},u_n^w;z) \rightrightarrows F_\mu(z), \quad
z\in\C\backslash(\supp\,\mu\cup\triangle_{\alpha}(\lambda)\cup\{w\}),
$$
where $\alpha\in(0,\pi]$ is given by
$
\sin\frac{\alpha}{2}=
\max\{\sin\frac{\alpha_2}{2}-\sin\frac{\alpha_1}{2},
\cos\frac{\alpha_1}{2}-\cos\frac{\alpha_2}{2}\}
$
and $w$ is arbitrarily chosen in $\Gamma_\alpha(\lambda)$.

\item
$(\frac{a_{n+1}}{a_n})_{n\geq1}$ converging to the unit circle.

Let us suppose that $\frak L\{\frac{a_{n+1}}{a_n}\}\subset\overline\Gamma_\zeta(\lambda)$,
$\lambda\in\T$, $\zeta\in[0,\pi)$. Theorems \ref{DERIVED22}, \ref{LP-C-2} and Corollary
\ref{DERIVED21} imply that, if $\ds\sin\frac{\zeta}{2}<\limi_n|a_n|$, then
$$
K_{2n}(a_1,\dots,a_{n-1},u_n^w;z) \rightrightarrows F_\mu(z),
\kern5pt
z\in\C\backslash(\supp\,\mu\cup\triangle_{\alpha-\zeta}(\lambda)\cup\{w\}),
$$
where $\alpha\in[0,\pi]$ is given by $\sin\frac{\alpha}{2}=\limi_n|a_n|$ and $w$ is
any point in $\Gamma_{\alpha-\zeta}(\lambda)$.

\hskip10pt In particular, in the L\'opez class
$\lim_n\frac{a_{n+1}}{a_n}=\lambda\in\T$, $\lim_n|a_n|\in(0,1)$, we have $\zeta=0$ and
$\{\supp\,\mu\}'=\triangle_\alpha(\lambda)$, hence,
$$
K_{2n}(a_1,\dots,a_{n-1},u_n^w;z) \rightrightarrows F_\mu(z), \quad
z\in\C\backslash(\supp\,\mu\cup\{w\}),
$$
if we choose $w\in\Gamma_\alpha(\lambda)$.

\ee

\eex

\bigskip
\noindent{\bf Acknowledgements}
\smallskip

\noindent The work of the authors was supported by Project E-12/25 of DGA
(Diputaci\'on General de Arag\'on) and by Ibercaja under grant IBE2002-CIEN-07.

\end{document}